%% file: main.tex
\documentclass[final,3p,times]{elsarticle}
\usepackage{graphicx}
\usepackage{amsmath, bm}
\usepackage{amssymb}
\usepackage{mathtools}
\usepackage{caption}
\usepackage{subcaption}
\usepackage{natbib}
\usepackage{xcolor}
\usepackage[hidelinks]{hyperref}
\usepackage{fix-cm}

\usepackage{tikz}
\usepackage{tikz-3dplot}
\usepackage{pgfplots}
\usetikzlibrary{
        angles,
        quotes,
        arrows.meta,
        backgrounds,
        calc,
        decorations.markings
}
\usepgfplotslibrary{external, groupplots}
\pgfplotsset{width=10cm,compat=1.9}

\pgfdeclareplotmark{o*}{
    \pgfpathcircle{\pgfpointorigin}{2pt}\pgfusepath{stroke}
    \pgfpathcircle{\pgfpointorigin}{1pt}\pgfusepath{fill}
}

\input{macros}

\input{colors}

\journal{Computational Physics}

\begin{document}

\begin{frontmatter}





\title{Beyond first-order accuracy in continuous-forcing immersed boundary methods, and their well-conditioned projection-based solution}





\author{Diederik Beckers\corref{cor1}} 
\ead{beckers@caltech.edu}
\cortext[cor1]{Corresponding author}

\affiliation{organization={Lynn Booth and Kent Kresa Department of Aerospace, California Institute of Technology},
            addressline={1200 East California Boulevard}, 
            city={Pasadena},
            postcode={91125}, 
            state={CA},
            country={USA}}

\author{H. Jane Bae} 

\affiliation{organization={Lynn Booth and Kent Kresa Department of Aerospace, California Institute of Technology},
            addressline={1200 East California Boulevard}, 
            city={Pasadena},
            postcode={91125}, 
            state={CA},
            country={USA}}

\author{Andres Goza} 

\affiliation{organization={Grainger College of Engineering, Department of Aerospace Engineering, University of Illinois Urbana-Champaign},
            addressline={104 S Wright Street}, 
            city={Urbana},
            postcode={61801}, 
            state={IL},
            country={USA}}

\begin{abstract}
    We introduce a refined immersed boundary (IB) methodology that is better-than-first-order accurate in practice, while preserving key properties of ``continuous-forcing'' IB approaches that retain a singular source term in the governing equations. Our method leverages a smoothed indicator (Heaviside) function, following ideas from multiphase flow and immersed layers formulations, to recast the IB solution as a composite of distinct interior and exterior fields. We demonstrate that, when cast through this composite-solution lens, prior continuous-forcing IB methods can be seen as neglecting terms in the governing and constraint equations that restrict the solution to first-order accuracy. We incorporate these terms to systematically improve accuracy without the need for heuristic corrections. In canonical Poisson problems, we empirically demonstrate second-order convergence, and in incompressible Navier--Stokes simulations the method achieves slightly sub-second-order performance. While our present study focuses on these cases, the framework suggests a path towards second-order accuracy or higher, with further extensions. This perspective reframes accuracy limitations typically attributed to IB schemes. Although continuous-forcing IB methods are often reported to be only first-order accurate, we show that neither smoothing nor interface interpolation inherently restricts attainable order. Moreover, we naturally incorporate this higher-order formulation into a projection-based solution process. The resulting algorithm simultaneously mitigates the spurious surface stresses produced by ill-conditioned linear systems and reduces sensitivity to geometric resolution, addressing both conditioning and accuracy concerns within a unified approach.
\end{abstract}


\end{frontmatter}

\section{Introduction}

Immersed boundary (IB) methods are used to solve partial differential equations in a physical domain $\domain \subset \R^3$ while enforcing interface conditions on an IB: an arbitrarily shaped surface $\surf \subset \domain$ lying strictly inside the boundaries of the computational domain. The key feature of these methods is that the discretization of $\domain$ need not conform to $\surf$, enabling the use of solvers optimized for Cartesian grids on simple geometries and avoiding expensive mesh (re-)generation, which is particularly advantageous when $\surf$ moves in time.

In the original IB method of~\citet{peskin_flow_1972}, the effect of $\surf$ is accounted for through a spatially singular term in the governing equation that transforms a surface forcing distribution $\scaforcedist(\surfc)$, where $\surfc$ is a set of surface coordinates on $\surf$, to the solution domain,
\begin{equation}
    \scaforce(\x) = \int_\surf \scaforcedist(\surfc) \delta \big(\x - \X(\surfc) \big) \dsurf , \quad \x \in \Omega,
    \label{eq:cont_reg}
\end{equation}
where $\delta(\x)$ is the three-dimensional Dirac delta function (DDF) and $\X(\surfc) = \big(X(\surfc),Y(\surfc),Z(\surfc)\big) \in \surf$ is the Lagrangian coordinate of a given point on the IB.
Similarly, interface conditions at the IB are specified using the sifting property of the Dirac delta function to evaluate the solution at the IB,
\begin{equation}
    \pot(\surfc) = \int_\domain \pot(\x) \ddf \big(\x - \X(\surfc) \big) \dvol.
    \label{eq:cont_interp}
\end{equation}

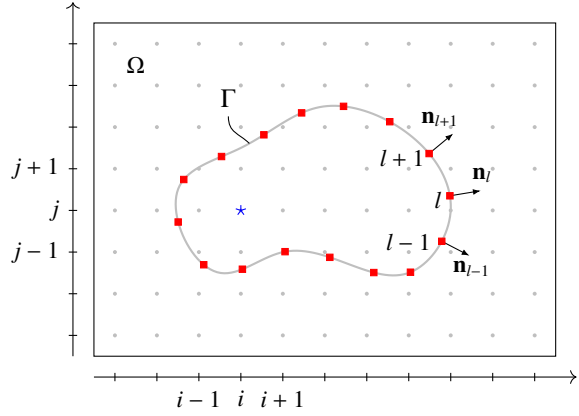
\begin{figure}
    \centering
  \tikzsetnextfilename{figures/2D_surface_discretization}%
  \input{figures/2D_surface_discretization.tikz}%

    \caption{Two-dimensional example of a surface $\surf$ immersed in a rectangular domain $\domain$. $\surf$ is discretized with a set of IB points (red squares) and $\domain$ is discretized with a Cartesian grid (gray dots). The position vector of the grid point indicated by the blue star is $\xgridvec[(i,j)]$.}
    \label{fig:domain_example}
\end{figure}

In this approach, $\surf$ is discretized with Lagrangian markers---material points tracked in time. Since these points generally do not coincide with the Eulerian grid used to discretize $\domain$, the singular DDF is replaced by a regularized (i.e., smoothed) DDF. For a point $\x = (x,y,z)$ in three dimensions and a grid that is uniform in each direction with spacings $\dxvec = (\dx, \dy, \dz) \in \R^3$, this regularized DDF is typically defined as the tensor product 
\begin{equation}
    \ddf_\dxvec(\x) = \ddf_\dx (x) \ddf_\dy (y) \ddf_\dz (z),
    \label{eq:ddf}
\end{equation}
where $\ddf_\dx (x) = \phi(x / \dx) / \dx$, with $\ddf_\dy (y)$ and $\ddf_\dz (z)$ defined analogously. Here $\phi$ is a smooth, compactly supported kernel function spanning only a few grid cells.
This construction ensures that Lagrangian forces are spread to nearby Eulerian grid points, and $\phi$ is chosen such that the regularized DDF satisfies $\sum_{i,j,k} \delta_\dxvec(\xgridvec[(i,j,k)] - \X) = 1$ for all $\X$. We use a non-italic sans-serif font to denote discretized Eulerian grid quantities and use the subscript ${(i,j,k)}$ (with parentheses) to index them on a Cartesian grid. For example, $\xgridvec[(i,j,k)]$ denotes the coordinates of an Eulerian grid point on a three-dimensional Cartesian grid, indexed by $i,j,k \in \Z$, and Fig.~\ref{fig:domain_example} demonstrates the discretization for a two-dimensional example.

We discretize $\surf$ with a set of $\numpts$ markers and use the single subscript $l$ (without parentheses) to denote a quantity associated to the $l$-th marker. Discrete surface quantities are denoted by a non-italic serif font. For example, $\xsurfvec[l] \coloneqq \X(\surfc[l])$ represents the Lagrangian coordinate of the $l$-th marker, and $\surfarea[l]$ is its associated surface area. Then we can discretize the right-hand side of Eq.~\eqref{eq:cont_reg} as a Riemann sum,
\begin{equation}
    \sum_{l} \scaforcedist(\surfc[l]) \ddf_\dxvec ( \x - \xsurfvec[l] ) \surfarea[l],
    \label{eq:peskin_reg}
\end{equation}
which is commonly referred to as the \emph{regularization} of the IB forcing. Similarly, the right-hand side of Eq.~\eqref{eq:cont_interp} becomes
\begin{equation}
    \sum_{i,j,k} \pot(\xgridvec[(i,j,k)]) \ddf_\dxvec ( \xgridvec[(i,j,k)] - \X ) \dx\dy\dz,
    \label{eq:peskin_interp}
\end{equation}
corresponding to a discrete convolution-type \emph{interpolation} using $\ddf_\dxvec$ as the kernel. Throughout the remainder of the paper, we refer to this discretization technique---the approximation of continuous convolution with a DDF kernel by a Riemann sum and regularized DDF---as the \emph{continuous-forcing discretization}. We similarly call methods that use this approach \emph{continuous-forcing IB methods}, following the classification of~\citet{mittal_immersed_2005}.

Both regularization and interpolation introduce errors, and the overall convergence rate depends on the type of IB problem. For regularization, \citet{beyer_analysis_1992} analyzed the IB method applied to the 1D heat equation with singular forcing. They show that, when the IB force is prescribed, the accuracy of the solution away from the interface is determined by the number of discrete moment conditions satisfied by the regularized DDF. \citet{tornberg_numerical_2004} extended this analysis to multidimensional Poisson problems and demonstrated that the same one-dimensional moment conditions control the accuracy of the multidimensional tensor-product DDF. Specifically, if a regularized, one-dimensional DDF satisfies the $m$ discrete moment conditions in grid direction $i$,
\begin{equation}
    \sum_i (\xgrid[(i)] - X)^r \ddf_\dx (\xgrid[(i)] - X) =
    \begin{cases}
        1, & r = 0,\\
        0, & 1 \le r \le m,
    \end{cases}
    \label{eq:moment_conditions}
\end{equation}
for all $X$, then, when used with a $m$-th order discretization of the underlying PDE, the IB method with prescribed forcing converges with order $m$ away from the interface and with first-order accuracy near the interface. This reduction to first order near the interface arises because the singular forcing induces jumps in the solution or its derivative, which the regularized IB forcing necessarily spreads over several grid points around the interface. As a result, the numerical solution gradually changes from the solution on one side to the other over these grid points, which only matches the sharp exact solution with first order accuracy. 

It is important to note that this result is not valid when the forcing is unknown a priori and instead has to be determined based on a constraint on the solution that involves interpolating the solution to the interface using Eq.~\eqref{eq:peskin_interp}. In this setting, the accuracy of the interpolation operation affects the accuracy of the solution. \citet{beyer_analysis_1992} show that the interpolation with discrete DDFs is dictated by the discrete moment conditions when the interpolation operation is associated with a smooth solution. However, for non-smooth solutions, as are typical when singular IB forcing creates distinct solutions that connect at an interface, the interpolation accuracy is reduced to first order. As the interpolated result is used to constrain the solution at the IB, the global accuracy of the solution is reduced to first order~\citep{beyer_analysis_1992}. Note that even a higher order interpolation scheme for non-smooth solutions would still produce a first-order method, because the near-interface solution values used in interpolation would be first-order accurate due to the force regularization discussed previously. The focus of this work is, therefore, to provide a method with greater than first-order accuracy by addressing these limitations.

Some works have addressed the accuracy limitations of the immersed boundary method while retaining the use of regularized DDFs. For example, \citet{griffith_order_2005} demonstrated that formally second-order accuracy can be achieved for sufficiently smooth solutions by giving the immersed boundary a finite thickness, which eliminates the jump discontinuities responsible for first-order errors in the classical formulation. Likewise, \citet{stein_immersed_2017} showed that higher-order accuracy can be obtained by constructing a smooth extension of the solution across the interface, such that the solution varies smoothly over the IB and the accuracy of Eq.~\eqref{eq:peskin_interp} is governed by the discrete moment conditions~\eqref{eq:moment_conditions}. However, both approaches impose structural constraints that reduce the generality of the IB method: the former requires assigning the interface a nonphysical thickness, and the latter effectively suppresses independent physical solution branches on either side of the interface.

The second focus of this work concerns the computation of the IB force. In Peskin’s original formulation \cite{peskin_flow_1972} the immersed boundary is a flexible, infinitesimally-thin membrane whose Lagrangian markers move with the surrounding incompressible viscous fluid. The force applied by the membrane is obtained from a constitutive relation between the marker displacements and stresses in massless elastic fibers connecting the markers. While this approach is natural for flexible structures, it becomes problematic when used for rigid bodies by imposing a large stiffness for the elastic fibers, which can lead to severe stability and stiffness issues in explicit or weakly-coupled time stepping.

Several alternative methods have been proposed that avoid these stability restrictions (see, e.g., the recent reviews of \citet{mittal_origin_2023} and \citet{verzicco_immersed_2023}). Among these, the distributed Lagrange multiplier (DLM) and the IB projection method are common for their ability to exactly enforce an arbitrary no-slip interface constraint without imposing specific constitutive relations or requiring ad-hoc parameters. They achieve this functionality by incorporating IB forces as Lagrange multipliers that enforce the no-penetration and no-slip conditions on the immersed boundary, analogous to the role of pressure as the Lagrange multiplier enforcing incompressibility. However, while these methods avoid stability restrictions in time-advancement, DLM and projection-based formulations become increasingly ill-conditioned as the ratio $\ds/\dx$ of the Lagrangian marker spacing $\ds$ to the Eulerian grid spacing $\dx$, decreases. A tradeoff must be found between numerical conditioning---which deteriorates rapidly for $\ds/\dx \lesssim 1$---and flow penetrating the immersed surface, which becomes significant for $\ds/\dx \gtrsim 1$, especially at high Reynolds numbers~\citep{beckers_high-fidelity_2026}. Values of $\ds/\dx$ used in prior studies vary: for example, \citet{taira_immersed_2007} recommend $\ds/\dx \approx 1$, whereas \citet{kallemov_immersed_2016} use $\ds/\dx \approx 2$.

\citet{goza_accurate_2016} showed that, for these methods that utilize the force as a Lagrange multiplier, ill-conditioning for small $\ds/\dx$ arises because of the nature of the linear equation that arises to compute the IB forcing strength. This equation, when cast using regularized DDFs, is an ill-posed first-kind integral equation; consequently, its discrete counterpart is inherently ill-conditioned. Intuitively, when $\ds/\dx$ is small, high-frequency variations in the surface-force distribution are smoothed by key ingredients in the IB method: regularization to the grid, solution of the Navier–Stokes equations with this forcing (including diffusion), and interpolation back to the boundary. As a result, the inverse problem is highly sensitive to small errors (e.g., from discretization or roundoff). That is, the correct IB forces to produce a physical flow are hard to separate from IB forces with the right nominal behavior but spurious deviations about it; either yields a similar velocity field. The numerical result of this sensitivity is that Lagrange multiplier-based IB methods produce stresses with incorrect high-frequency oscillations, computed from ill-conditioned linear systems that are more expensive to solve via iterative techniques than well-conditioned problems. \citet{goza_accurate_2016} proposed a postprocessing technique to smooth out this noise; however, it does not reduce the intrinsic cost of solving the system. Other approaches to mitigating the ill-posedness have also been applied: \citet{kallemov_immersed_2016} employ a preconditioner but still face constraints on the $\ds/\dx$ ratio, while \citet{yu_regularized_2024} use Tikhonov regularization, which requires selecting an ad-hoc regularization parameter.

In this work, we simultaneously address the accuracy challenge associated with continuous-forcing IB methods and the ill-conditioning of the projection-based approach for satisfying Dirichlet boundary conditions (such as the no-slip condition). The proposed method is framed within a projection-based formulation for the IB forces, to enable a robust treatment of flows past bodies moving with arbitrary prescribed kinematics. It retains a comparable computational cost per timestep to the first-order variants of that method, removes the ill conditioning of the linear system for the surface stress, and does not rely on any ad-hoc parameters. 
The approach formally distinguishes the solution on either side of the immersed interface using indicator (Heaviside) functions, as is common in immersed algorithms for multiphase flows~\citep{kataoka_local_1986,tryggvason_direct_2011}. We follow~\citet{tryggvason_front-tracking_2001} and, in particular,~\citet{eldredge_method_2022} to connect this methodology to the immersed boundary method. The key enabler of the proposed methodology is to use a Taylor series to expand the solution on either side of the immersed interface within the support of the regularized DDF in both the regularization and interpolation operations, while accounting for the smoothness of the indicator functions that arises from using that regularized DDF. This approach produces additional terms in the governing and constraint equations compared to the prototypical continuous-forcing IB formulation that simultaneously improve the accuracy of the operations and produce a second-kind (well posed) linear system for the unknown IB force. We present this approach first for a simple Dirichlet Poisson problem in Section~\ref{sec:poisson} and afterward for the incompressible Navier--Stokes equations with no-slip boundary conditions in Section~\ref{sec:navier-stokes}. Although we demonstrate the improved accuracy only using the projection-based approach to find the forces, our approach of modeling the solution behavior within the regularized DDF support to increase the accuracy of the formulation is applicable to any form of the IB with regularized DDFs.

\section{A higher-order accurate, well-conditioned projection IB method for a model Dirichlet Poisson problem}
\label{sec:poisson}

The surface $\surf$ partitions the domain $\domain$ into two subregions, $\domain^+$ and $\domain^-$, which we refer to as the \emph{exterior} and \emph{interior} regions, respectively. We will make use of the unit normal vector $\normal$ and the two tangent vectors $\tangentone$ and $\tangenttwo$ defined on $\surf$. By convention, the normal vector $\normal$ is oriented to point into $\Omega^+$. We also define the indicator functions $\heavi^+(\x)$ and $\heavi^-(\x)$ associated with the exterior and interior subdomains, respectively, where $\heavi^\pm(\x)=1$ for $\x\in\Omega^\pm$ and $\heavi^\pm(\x)=0$ otherwise.
In the distributional sense, the gradient of the indicator fields is a distribution of multidimensional Dirac delta functions supported on the interface multiplied by the interface normal, 
\begin{equation}
    \grad \heavi^\pm(\x) = \pm\int_\surf \normal(\surfc) \ddf \big(\x - \X(\surfc) \big) \dsurf .
    \label{eq:contheavi}
\end{equation}

We develop our method first for a Dirichlet Poisson problem on the global domain $\domain$ for a scalar field $\pot(x) \in \R$ that is defined as the composite (or piecewise) field,
\begin{equation}
    \pot \coloneqq \heavi^+ \pot^+ + \heavi^- \pot^-,
    \label{eq:composite_analytic_pot}
\end{equation}
where $\pot^{+}$ and $\pot^{-}$ are the exterior and interior solutions, respectively, each satisfying the Dirichlet Poisson problems
\begin{gather}
    \lap \pot^{\pm}(\x) = \source(\x), \quad \x \in \domain^\pm, 
    \label{eq:poisson_governeq} \\
    \pot^{\pm}(\x) = \pot^{\pm}_{\partial \domain}, \quad \x \in \partial \domain^\pm \setminus \surf,
    \label{eq:poisson_bc} \\
    \pot^{\pm}(\X(\surfc)) = \pot_{\surf}^{\pm}(\surfc), \quad \surfc \in \surf,
    \label{eq:poisson_ibc}
\end{gather}
where $\source$ is an integrable source function, $\pot^{\pm}_{\partial \domain}$ is a prescribed boundary condition at the boundary of the exterior and interior domains (excluding $\surf$), and $\pot^{\pm}_{\surf}$ is a prescribed function on $\surf$.
For concreteness, we are only treating the case where the interface condition on the IB is a Dirichlet condition that is the same for the exterior and interior solutions, i.e., $\pot^{\pm}_{\surf} = \pot_{\surf}$, resulting in a continuous solution, and we assume $\source$ varies smoothly across the interface. However, the method we develop can easily be adapted to account for discontinuities in the solution and source function or different types of interface conditions, such as Neumann conditions.

By applying the product rule for the gradient and divergence operators consecutively to the composite solution \eqref{eq:composite_analytic_pot}, we obtain a governing equation for $\pot$:
\begin{equation}
    \lap \pot = \source + \grad \heavi^+ \cdot \big( \grad \pot^{+} - \grad \pot^{-} \big) + \div \Big( \grad \heavi^+ \big( \pot^{+} - \pot^{-} \big) \Big), \quad \x \in \domain. 
    \label{eq:poisson_ilm_governeq}
\end{equation}

The immersed layers methodology of \citet{eldredge_method_2022} uses \eqref{eq:contheavi} together with the properties of the DDF to show that, in the case of interest where $\pot^{+}_{\surf} = \pot^{-}_{\surf}$, Eq.~\eqref{eq:poisson_ilm_governeq} reduces to an equation with a single-layer source term due to the jump in the solution gradient across the interface. A continuous-forcing IB method then discretizes this equation using regularized DDFs \citep{goza_accurate_2016,beckers_planar_2022} to solve \eqref{eq:poisson_governeq}–\eqref{eq:poisson_ibc}. In this formulation, the single-layer term is formally represented as an unknown source term that enforces the Dirichlet interface condition on $\surfc$, i.e., that $\pot$ evaluated on the interface---via the sifting property of the DDF---matches the prescribed surface function $\pot_{\surf}$. This treatment incurs a challenge related to accuracy and a challenge related to numerical stability. 

For the accuracy-related challenge, once the computed source term is regularized onto the governing equations, it produces a solution field $\pot$ that is not discretely consistent with its underlying masked identity \eqref{eq:composite_analytic_pot} near the interface for a staggered grid. This inconsistency, in turn, means that one cannot improve the global accuracy of the method by only increasing the accuracy of the interpolation scheme using knowledge of the discrete behavior of \eqref{eq:composite_analytic_pot} near the interface.

For the numerical challenge, the equation for the source term---
when cast using regularized DDFs to discretize the continuous governing equation and interface condition---is a discretization of a first-kind integral equation that is ill posed. The discrete equation is therefore ill conditioned and prone to producing inaccurate, unphysically oscillatory stresses \cite{goza_accurate_2016,kallemov_immersed_2016,yu_regularized_2024}.

To avoid both challenges, we work directly with a discretization of \eqref{eq:composite_analytic_pot} and induce its differential equations through application of discrete divergence and gradient operations. In deriving the formulation, we ensure that the discrete solution fields satisfy the desired Poisson problem on the exterior and interior regions of the flow, and formally keep any terms along the interface that arise. This derivation produces a discrete analog of \eqref{eq:poisson_ilm_governeq} that rigorously preserves the discrete representation of the composite solution field. 

We will demonstrate in this manuscript that our proposed approach addresses both challenges to continuous-forcing IB methods. First, the algorithm discretely preserves the composite solution representation. This property directly allows for the solution to be higher order accurate, since the interpolation step is now applied to a discretely consistent composite solution field. We will moreover show that a Taylor series expansion of the solution field near the interface allows for a solution algorithm that preserves the essential elements of the first-order prototypical continuous-forcing IB procedure, so that the increase in accuracy comes at a marginal increase in cost. 

Second, we retain the connection of the single-layer source term to the jump in the solution gradient directly in the formulation. This identity, when incorporated with appropriate Taylor series expansions utilized near the interface, leads to a \emph{second-kind} integral equation for these unknown jump terms. This equation is inherently well posed, irrespective of the ratio of point spacing for the IB versus the flow domain. The outcome is a well-conditioned linear system for computing the unknown jump terms, that is conducive to iterative methods and produces smooth, physical surface stresses.

\subsection{Discrete problem formulation for a composite solution}

\begin{figure}
    \centering
  \tikzsetnextfilename{figures/3D_grid_cell}%
  \input{figures/3D_grid_cell.tikz}%

    \caption{Three-dimensional grid cell. The labeled locations in the right panel each share the same index tuple $(i, j, k)$.}
    \label{fig:cell}
\end{figure}
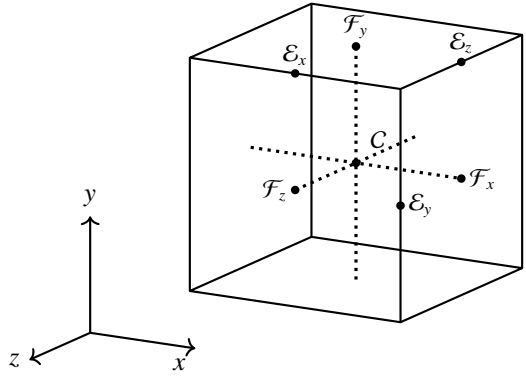

We discretize $\Omega$ using a three-dimensional, staggered Cartesian grid. 
The grid has constant cell dimensions $\dx, \dy, \dz \in \R$ and contains three sets of discrete locations (see Fig.~\ref{fig:cell}): cell centers $\centers$, cell faces $\faces = (\facescomponent_x, \facescomponent_y, \facescomponent_z)$, and cell edges $\edges = (\edgescomponent_x, \edgescomponent_y, \edgescomponent_z)$. We define the vector spaces $\R^{\centers}$, $\R^{\facescomponent_x}$, $\R^{\facescomponent_y}$, $\R^{\facescomponent_z}$, $\R^{\edgescomponent_x}$, $\R^{\edgescomponent_y}$, and $\R^{\edgescomponent_z}$ for the real-valued fields discretized at those cell locations, and will use the notation $\R^{\faces} \coloneqq \R^{\facescomponent_x} \times \R^{\facescomponent_y} \times \R^{\facescomponent_z}$ and $\R^{\edges} \coloneqq \R^{\edgescomponent_x} \times \R^{\edgescomponent_y} \times \R^{\edgescomponent_z}$ for discrete spaces containing all three components of vector-valued fields. For example, the data structure containing all the elements of a discretized field of scalar data (such as pressure) on the cell centers is in $\R^{\centers}$. Similarly, the data structure containing all elements of all three components of a discretized vector field (such as velocity) on the respective cell faces is in $\R^{\faces}$. The cell edges store the results of curl operations applied to vector fields defined on cell faces, and in Sec.~\ref{sec:navier-stokes} they are also used to store tensor-valued data. Additionally, we describe the positions of the cell centers by $\xgridvec[\centers] = (\xgrid[\centers], \ygrid[\centers], \zgrid[\centers]) \in \R^{\centers} \times \R^{\centers} \times \R^{\centers}$, and similarly for the faces and edges, for example, $\xgridvec[\facescomponent_x] = (\xgrid[\facescomponent_x], \ygrid[\facescomponent_x], \zgrid[\facescomponent_x]) \in \R^{\facescomponent_x} \times \R^{\facescomponent_x} \times \R^{\facescomponent_x}$. Throughout the rest of this paper we will use second-order finite differences to discretize our differential operators.

We denote the discrete approximations to the indicator fields as $\hgridcenters^\pm \in \R^{\centers}$ (their construction is detailed in the next section) and we define the discrete composite solution $\mask{\potgrid} \in \R^{\centers}$ on the cell centers,
\begin{equation}
    \mask{\potgrid} = \hgridcenters^+ \had \potgrid^+ + \hgridcenters^- \had \potgrid^-,
    \label{eq:composite_pot}
\end{equation}
where $\potgrid^\pm \in \R^{\centers}$ are the discrete interior and exterior solutions and $\had$ is the element-wise product between two fields defined on the same space.
We then aim to formulate a system of equations for this composite solution such that the discrete interior and exterior solutions satisfy Eq.~\eqref{eq:poisson_governeq}--\eqref{eq:poisson_bc} when discretized on our grid,
\begin{equation}
    \lapgrid \potgrid^\pm = \sourcegrid + \bcgrid^\pm,
    \label{eq:int_ext_governeq}
\end{equation}
and also satisfy the interface conditions at the points that discretize $\surf$. Here, $\lapgrid: \R^{\centers} \mapsto \R^{\centers}$ is the discrete Laplacian, $\bcgrid^\pm \in \R^{\centers}$ represents the boundary-condition term that arises from discretizing the Laplacian while enforcing Eq.~\eqref{eq:poisson_bc}, and $\sourcegrid \in \R^{\centers}$ is the discretized source function. We use the overbar notation in Eq.~\eqref{eq:composite_pot} to indicate the composite nature of the solution and the potential smoothing near the IB, which depends on whether the discrete indicator functions are themselves smoothed. This contrasts with the analytic solution~\eqref{eq:composite_analytic_pot}, which is expressed using sharp indicator functions.

When the indicator functions are smoothed, the solution from each side will bleed into the other side and we need an extension of the solution over that support. One of the key aspects of our method is that we enforce that $\pot^\pm$ vary smoothly across the IB. Away from the IB, and outside of their respective regions, their extension is irrelevant provided it is finite, well-defined, and will still make the solution satisfy the governing equation and interface condition. The proposed approach works directly with the composite solution that stitches together the exterior and interior solutions, so that the superfluous interior (exterior) part of the exterior (interior) solution is not needed, in contrast to the method of~\citet{stein_immersed_2017}.

We can multiply Eq.~\eqref{eq:int_ext_governeq} for the interior and exterior solutions by their respective indicator fields and add the two equations:
\begin{equation}
   \hgridcenters^+ \had \lapgrid \potgrid^+ + \hgridcenters^- \had \lapgrid \potgrid^- = \sourcegrid + \mask{\bcgrid}
\end{equation}
where $\mask{\bcgrid} \coloneqq \hgridcenters^+ \had \bcgrid^+ + \hgridcenters^- \had \bcgrid^-$ and we required that $\hgridcenters^+ + \hgridcenters^- = 1$. Then, using the discrete product rule~\eqref{eq:composite_scalar_lap} for the Laplacian, we can recast this into a governing equation for $\mask{\potgrid}$ that is the discrete equivalent of Eq.~\eqref{eq:poisson_ilm_governeq}, following the approach of~\citet{eldredge_method_2022}:
\begin{equation}
    \lapgrid \mask{\potgrid} = \sourcegrid + \mask{\bcgrid}
     + \spacetransform{\faces}{\centers} \left( \gradgrid \hgridcenters^+ \had \left( \gradgrid \potgrid^+ - \gradgrid \potgrid^- \right) \right)
    + \divgrid \left( \gradgrid \hgridcenters^+ \had \spacetransform{\centers}{\faces}\left( \potgrid^+ - \potgrid^- \right) \right),
    \label{eq:poisson_composite}
\end{equation}
where $\divgrid: \R^{\faces} \mapsto \R^{\centers}$ and $\gradgrid: \R^{\centers} \mapsto \R^{\faces}$ are the discrete divergence and gradient operators, respectively, defined in \ref{app:differential_operators}. The operation $\spacetransform{\centers}{\faces}: \R^{\faces} \mapsto \R^{\centers}$ transforms a vector on the cell faces to a scalar on the cell centers by interpolating each vector component to the cell center and adding them together, and $\spacetransform{\faces}{\centers}: \R^{\centers} \mapsto \R^{\faces}$ is the adjoint operation, which interpolates a scalar on the cell faces to each vector component on the cell faces (see~\ref{app:grid_transformations} for more details). 

Equation \eqref{eq:poisson_composite} is not yet in a convenient form for solving because it involves both the composite and the individual solutions. However, the gradient of the indicator function is zero everywhere except at the interface, where it becomes unbounded for sharp indicator functions. For smoothed indicator functions---where the smoothing is limited to a finite-width region around the interface---the gradient remains nonzero but finite in that region. We can take advantage of this compact support to rewrite the terms in Eq.~\eqref{eq:poisson_composite} that involve the interior and exterior solutions so that they depend only on their behavior at the interface. Moreover, because of this compact support, discretizing the divergence and gradient operators does not introduce any additional boundary-condition terms.

\subsection{Finding the discrete indicator fields and their gradients from the regularized DDF}
\label{sec:poisson_taylor_series}

 At this stage---and only at this stage---we apply the continuous-forcing discretization formula~\eqref{eq:peskin_reg} to find the discrete analogue of the relation~\eqref{eq:contheavi} between the gradient of the indicator fields and the DDF. To proceed, we first introduce a more detailed notation tailored to our grid.

We discretize $\surf$ into a collection of $\numpts$ surface patches whose centers we track as Lagrangian points. We discretize both scalar- and vector-valued surface quantities on these points and define the spaces of these data structures as $\R^{\spoints}$ and $\R^{\vpoints} \coloneqq \R^{\spoints} \times \R^{\spoints} \times \R^{\spoints}$, respectively. 
The positions of the points are described by $\xsurfvec = (\xsurf, \ysurf, \zsurf) \in \R^{\vpoints}$ and the local normal and tangential vectors by $\normalvec, \tanonevec, \tantwovec \in \R^{\vpoints}$. The areas of the surface patches are given by $\surfarea \in \R^{\spoints}$. We use the notation $\had$ to denote element-wise multiplication between these quantities. This operation applies component-by-component multiplication between arrays of matching type: for scalar-valued arrays, it corresponds to the standard Hadamard product; for vector-valued arrays, it multiplies corresponding vector components. 

To simplify our notation, we define the field $\dddf_{\centers,l} \in \R^{\centers}$ that results from evaluating the regularized DDF $\ddf_\dxvec$ centered at the $l$-th surface point at the grid points:
\begin{equation}
    \dddf_{\centers,l} \coloneqq \ddf_\dxvec \Big( \xgrid[\centers] - \xsurf[l], \ygrid[\centers] - \ysurf[l], \zgrid[\centers] - \zsurf[l] \Big),
    \label{eq:dddf}
\end{equation}
and we define $\dddf_{\facescomponent_x,l}$, $\dddf_{\facescomponent_y,l}$, and $\dddf_{\facescomponent_z,l}$ similarly. In this work we use the smoothed three-point delta function from \citet{yang_smoothing_2009} to construct $\ddf_\dxvec$ via~\eqref{eq:ddf}. Using these fields, we can define the regularization operators $\regds_{\centers}: \R^{\spoints} \mapsto \R^{\centers}$ and $\regds_{\faces}: \R^{\vpoints} \mapsto \R^{\faces}$ to regularize surface scalar data to the cell centers and surface vector data to the cell faces by summing over all the discrete surface points. For example, the vector field (with each component indexed by $i,j,k$) that one gets from regularizing the vector-valued surface data $\vectorpoint = (\vectorpointcomponent[x], \vectorpointcomponent[y], \vectorpointcomponent[z]) \in \R^{\vpoints}$ (indexed by $l$) to the grid can be formulated as
\begin{equation}
    \indicesbig{\regds_{\faces} \vectorpoint}{(i,j,k)}
    =
    \begin{pmatrix}
        \indicesbig{\regds_{\facescomponent_x} \vectorpointcomponent[x]}{(i,j,k)} \\
        \indicesbig{\regds_{\facescomponent_y} \vectorpointcomponent[y]}{(i,j,k)} \\
        \indicesbig{\regds_{\facescomponent_z} \vectorpointcomponent[z]}{(i,j,k)}
    \end{pmatrix}
    \coloneqq \sum_l
    \begin{pmatrix}
        \indicesbig{\dddf_{\facescomponent_x,l}}{(i,j,k)} \vectorpointcomponent[x,l] \\
        \indicesbig{\dddf_{\facescomponent_y,l}}{(i,j,k)} \vectorpointcomponent[y,l] \\
        \indicesbig{\dddf_{\facescomponent_z,l}}{(i,j,k)} \vectorpointcomponent[z,l]
    \end{pmatrix} \surfarea[l],
    \label{eq:disc_reg}
\end{equation}
which is the discrete vector-valued version of Eq.~\eqref{eq:cont_reg}, specialized for our staggered grid. Similarly, we can discretize Eq.~\eqref{eq:cont_interp} and its vector-valued variant for our grid using the operators $\interpbody_{\centers}: \R^{\centers} \mapsto \R^{\spoints}$ and $\interpbody_{\faces}: \R^{\faces} \mapsto \R^{\vpoints}$. For example, the scalar-valued surface data (indexed by $l$) that one gets from interpolating $\scalargrid \in \R^{\centers}$ (indexed by $i,j,k$) to the IB points can be formulated as
\begin{equation}
    \indicesbig{\interpbody_{\centers} \scalargrid}{l}
    \coloneqq
    \dx\dy\dz \sum_{i,j,k} \indicesbig{\dddf_{\centers,l}}{(i,j,k)} \scalargrid[(i,j,k)].
    \label{eq:disc_interp}
\end{equation}
We refer the reader to \ref{app:regularization} and \ref{app:interpolation} for a complete list of the regularization and interpolation operations used in this work.

Using the previous notation, we can discretize the relation~\eqref{eq:contheavi} between the indicator fields and DDF as in~\cite{unverdi_front-tracking_1992,tryggvason_direct_2011,eldredge_method_2022}:
\begin{equation}
    \gradgrid \hgridcenters^\pm = \pm \regds_{\faces} \normalvec.
    \label{eq:peskin_heavi}
\end{equation}
If we apply the discrete divergence operation to both sides and require that our discrete operators satisfy mimetic properties~\cite{liska_fast_2016}, we obtain the following discrete Poisson equation involving the discrete Laplacian (since $\divgrid \gradgrid = \lapgrid$)~\cite{unverdi_front-tracking_1992,tryggvason_direct_2011,eldredge_method_2022}:
\begin{equation}
    \lapgrid \hgridcenters^\pm = \pm \divgrid \regds_{\faces} \normalvec + \bcgrid_\hgrid.
    \label{eq:poisson_hgrid}
\end{equation}
By applying appropriate Dirichlet boundary conditions $\bcgrid_\hgrid$, we can solve for the indicator field whose smoothness is determined by the smoothness of the regularized DDF used in~\eqref{eq:dddf} and which satisfies $\hgridcenters^+ + \hgridcenters^- = 1$. However, it is important to note that if $\regds_{\faces} \normalvec$ has a component in the null space of $\divgrid$ (i.e., is not curl-free)---as is the case for most constructions of $\regds$ in combination with curved immersed boundaries---then the discrete indicator field obtained by solving~Eq.~\eqref{eq:poisson_hgrid} will generally not satisfy Eq.~\eqref{eq:peskin_heavi}. Instead, the gradient satisfies
\begin{equation}
    \gradgrid \hgridcenters^\pm = \pm \regds_{\faces} \normalvec \pm \curlgrid \lapgrid^{-1}_{\edges} \curlgrid^\top \regds_{\faces} \normalvec ,
    \label{eq:true_heavi}
\end{equation}
where $\lapgrid_{\edges}: \R^{\edges} \mapsto \R^{\edges}$ is the discrete Laplacian for variables on edges and $\curlgrid: \R^{\edges} \mapsto \R^{\faces}$ and $\curlgrid^\top: \R^{\faces} \mapsto \R^{\edges}$ are the discrete curl operation and its adjoint, respectively (see \ref{app:differential_operators}).
We will ignore this last term in \eqref{eq:true_heavi} for the rest of our derivation. This omission is a source of error in the global solution, which can be avoided by working with a curl-free $\regds_{\faces} \normalvec$. This topic is part of our ongoing investigation.


If we substitute Eq.~\eqref{eq:peskin_heavi} into Eq.~\eqref{eq:poisson_composite}, we obtain
\begin{equation}
    \lapgrid \mask{\potgrid} = \sourcegrid + \mask{\bcgrid}
     + \spacetransform{\faces}{\centers} \Big( \regds_{\faces} \normalvec \had ( \gradgrid \potgrid^+ - \gradgrid \potgrid^- ) \Big)
    + \divgrid \Big( \regds_{\faces} \normalvec \had \spacetransform{\centers}{\faces} ( \potgrid^+ - \potgrid^- ) \Big),
    \label{eq:poisson_composite_Rn}
\end{equation}
which still contains the interior and exterior solutions, but now along with the DDF, which is only nonzero is a small region at the interface that scales with the grid spacing. The last step to formulate our governing equation for $\mask{\potgrid}$ is to replace the interior and exterior solutions (and their gradients) by their behavior inside this region using a Taylor series expansion. This way, the governing equation can be expressed purely in terms of $\mask{\potgrid}$ and variables defined on the IB, as we have no access to the $\potgrid^+$ and $\potgrid^-$ fields directly. 

\subsection{Using Taylor series to find leading-order terms}

When the DDFs are centered at the discrete surface points and element-wise multiplied by a field variable, we can approximate the variation of this variable within the support of the discrete DDF by a truncated Taylor series expansion about the discrete surface point at which the DDF is centered. Using the truncated Taylor series~\eqref{eq:s_TS} and~\eqref{eq:Gs_TS}, and the normal-distance-weighted regularization operation $\regds_{\faces, 1 \normalbase}: \R^{\vpoints} \mapsto \R^{\faces}$, defined in~\eqref{eq:RF1n}, our final formulation for Eq.~\eqref{eq:poisson_composite_Rn} becomes
\begin{equation}
    \lapgrid \mask{\potgrid} = \sourcegrid + \mask{\bcgrid}
     + \spacetransform{\faces}{\centers} \regds_{\faces} \big( \normalvec \had \normalvec \had \spacetransform{\spoints}{\vpoints} \jump{\potpointnormal} \big)
    + \divgrid \Big( \regds_{\faces, 1 \normalbase} \big( \normalvec \had \spacetransform{\spoints}{\vpoints} \jump{\potpointnormal} \big) \Big),
    \label{eq:poisson_composite_TS}
\end{equation}
where we used the notation $\jump{\scalarpoint} \coloneqq \scalarpointsurfplus - \scalarpointsurfminus$ to represent the difference between the exterior and interior values (or the \emph{jump}) of a variable $\scalarpoint \in \R^{\spoints}$ at the IB and where $\spacetransform{\spoints}{\vpoints}$ represents the transformation~\eqref{eq:vIs} of scalar-valued surface data to vector-valued surface data. We use the superscripts $\normalbase$, $\tanonebase$, and $\tantwobase$ to denote the directional derivatives along $\normalvec$, $\tangentone$, and $\tangenttwo$, respectively; i.e., $s^\normalbase \coloneqq \normalvec \cdot \nabla s$ and $\scalarpointnormal \in \R^{\spoints}$ is its discrete counterpart, with analogous definitions for $\tanonebase$ and $\tantwobase$. To obtain Eq.~\eqref{eq:poisson_composite_TS}, we used the fact that $\jump{\potpoint} = 0$, $\jump{\potpointtanone} = 0$, and $\jump{\potpointtantwo} = 0$, since $\potpointsurfplus = \potpointsurfminus$. This choice is for clarity of presentation; jumps in the solution value across the interface can be incorporated directly.

We can contrast this formulation with the prototypical continuous-forcing IB formulation, which can be written generically here as
\begin{equation}
    \lapgrid \mask{\potgrid} = \sourcegrid + \mask{\bcgrid} + \regds_{\centers} \forcepoint.
    \label{eq:poisson_classic_IB}
\end{equation}
This is also the discrete formulation that results from the immersed layers approach in~\citet{eldredge_method_2022} if the double layer strength (i.e., the solution jump over the IB) is set to zero. The key differences with our proposed formulation are that Eq.~\eqref{eq:poisson_classic_IB} ignores the last term on the right-hand side of Eq.~\eqref{eq:poisson_composite_TS} and uses $\regds_{\centers}$ instead of $\spacetransform{\faces}{\centers}\regds_{\faces}$. 
However, the DDF evaluated at the cell centers is not exactly equal to the DDF evaluated at the cell faces and, subsequently, interpolated to the cell centers. We found empirically that using $\regds_{\centers}$ to regularize the single-layer forcing to our staggered grid limits the global accuracy of our approach to first order when the indicator fields are found using \eqref{eq:poisson_hgrid}.
Because it is known that, in this Poisson problem, the forcing represents the jump in the normal derivative of the solution \citep{wiegmann_explicit_1998,eldredge_method_2022}, we will also refer to $\forcepoint = \jump{\potpointnormal}$ as the (IB) forcing strength.

We can also establish the error scaling (with respect to the flow grid spacing) in satisfying the interface (interpolation) condition using our proposed approach. To do so, we perform a truncated Taylor series expansion about the IB points for the interior and exterior solutions using formula~\eqref{eq:disc_interp}:
\begin{align}
    &\indicesbig{\potpointsurf - \interpbody_{\centers} \mask{\potgrid}}{l} \nonumber
    \\
    &\quad= \potpointsurf[,l] - \dx\dy\dz \sum_{i,j,k} \indicesbig{\dddf_{\centers,l}}{(i,j,k)} \left( \indicesbig{\hgridcenters^+ \had \potgrid^+}{(i,j,k)} + \indicesbig{\hgridcenters^- \had \potgrid^-}{(i,j,k)} \right) \nonumber
    \\
    &\quad= \potpointsurf[,l] - \dx\dy\dz \sum_{i,j,k} \indicesbig{\dddf_{\centers,l}}{(i,j,k)} \bigg[ \hgridcenters[,(i,j,k)]^{+} \Big( \potpointsurfplus[,l] + \normalvec[l] \cdot \big( \xgridvec[\centers,(i,j,k)] - \xsurfvec[l] \big) \potpointsurfnormalplus[,l]
    \Big) \nonumber\\
    &\quad\quad + \hgridcenters[,(i,j,k)]^{-} \Big( \potpointsurfminus[,l] + \normalvec[l] \cdot \big( \xgridvec[\centers,(i,j,k)] - \xsurfvec[l] \big) \potpointsurfnormalminus[,l]
    \Big) + \mathcal{O}(\dx^2) + \mathcal{O}(\dy^2) + \mathcal{O}(\dz^2) \bigg].
    \label{eq:Eu_TS}
\end{align}

Using the first discrete moment condition for the DDF (i.e., Eq. \eqref{eq:moment_conditions} with $m=0$) and the facts that $\potpointsurfplusminus=\potpointsurf$ and $\hgridcenters^+ + \hgridcenters^- = 1$, $\potpointsurf$ cancels out with the first term of the interior and exterior solutions. Additionally, for DDFs that satisfy two discrete moment conditions, we have that
\begin{equation}
    \begin{split}
        &\sum_{i,j,k} \indices{\dddf_{\centers,l}}{(i,j,k)} \normalvec[l] \cdot \big( \xgridvec[\centers,(i,j,k)] - \xsurfvec[l] \big) = 
        \\
        &\quad = \sum_{i,j,k} 
        \ddf_\dx \big( \xgrid[\centers,(i)] - \xsurf[l] \big)
        \ddf_\dy \big( \ygrid[\centers,(j)] - \ysurf[l] \big)
        \ddf_\dz \big( \zgrid[\centers,(k)] - \zsurf[l] \big)
        \Big[
             \normalvecc{x,l} \big( \xgrid[\centers,(i)] - \xsurf[l] \big) 
            +\normalvecc{y,l} \big( \ygrid[\centers,(j)] - \ysurf[l] \big) 
            +\normalvecc{z,l} \big( \zgrid[\centers,(k)] - \zsurf[l] \big)         
        \Big]
        \\
        &\quad = 
        \normalvecc{x,l} \bigg(\sum_{i} \ddf_\dx \big( \xgrid[\centers,(i)] - \xsurf[l] \big) 
         \big( \xgrid[\centers,(i)] - \xsurf[l] \big) \bigg)
        \bigg(\sum_{j}\ddf_\dy \big( \ygrid[\centers,(j)] - \ysurf[l] \big) \bigg)
        \bigg(\sum_{k}\ddf_\dz \big( \zgrid[\centers,(k)] - \zsurf[l] \big) \bigg)
        \\
        &\quad \quad + 
        \normalvecc{y,l}\bigg(\sum_{i} \ddf_\dx \big( \xgrid[\centers,(i)] - \xsurf[l] \big) \bigg)
        \bigg(\sum_{j}\ddf_\dy \big( \ygrid[\centers,(j)] - \ysurf[l] \big)  \big( \ygrid[\centers,(j)] - \ysurf[l] \big) \bigg)
        \bigg(\sum_{k}\ddf_\dz \big( \zgrid[\centers,(k)] - \zsurf[l] \big) \bigg)
        \\
        &\quad \quad + 
        \normalvecc{z,l}\bigg(\sum_{i} \ddf_\dx \big( \xgrid[\centers,(i)] - \xsurf[l] \big) \bigg)
        \bigg(\sum_{j}\ddf_\dy \big( \ygrid[\centers,(j)] - \ysurf[l] \big) \bigg)
        \bigg(\sum_{k}\ddf_\dz \big( \zgrid[\centers,(k)] - \zsurf[l] \big)  \big( \zgrid[\centers,(k)] - \zsurf[l] \big) \bigg)
        \\
        &\quad=0,
    \label{eq:EC1n1}
    \end{split}
\end{equation}
for all IB points.
Then, using $\hgridcenters^-= 1 - \hgridcenters^+$, we can formulate the following second-order accurate interpolation for our composite solution,
\begin{equation}
    \interpbody_{\centers} \mask{\potgrid} - \jump{\potpointnormal} \had \interpbody_{\centers, 1n} \hgridcenters^+ \approx \potpointsurf,
    \label{eq:disc_pot_interp_TS}
\end{equation} 
where we introduced the normal-distance-weighted interpolation operation $\interpbody_{\centers, 1\normalbase} : \R^{\centers} \mapsto \R^{\spoints}$, defined in~\eqref{eq:EC1n}. Prototypical continuous-forcing methods do not account for the second term on the left-hand side, as it does not arise from directly discretizing the continuous constraint formulated using a singular DDF. This omission, as well as the use of 
Eq.~\eqref{eq:poisson_classic_IB} for the governing equation instead of Eq.~\eqref{eq:poisson_composite_TS} are the two causes for the first-order accuracy of these methods. Eq.~\eqref{eq:disc_pot_interp_TS} has to be combined with the new governing equation~\eqref{eq:poisson_composite_TS} to obtain a second-order accurate solution away from the interface, since otherwise the Taylor series expansions that underlie the derivation of the new constraint equation do not hold.






\subsection{Projection-based approach to the force computation}
\label{sec:projection}

Our proposed IB formulation, consisting of the governing equation~\eqref{eq:poisson_composite_TS} and the corrected interface condition~\eqref{eq:disc_pot_interp_TS}, forms a saddle-point system. To write out this system, we replace the element-wise products $\had$ with multiplications by diagonal matrices, so that the IB formulation becomes
\begin{align}
    \begin{bmatrix}
        \lapgrid & \spacetransform{\faces}{\centers}\regds_{\faces} \diag(\normalvec)^2 \spacetransform{\spoints}{\vpoints} + \divgrid \regds_{\faces,1\normalbase} \diag(\normalvec) \spacetransform{\spoints}{\vpoints} \\
        \interpbody_{\centers} & \diag \big( \interpbody_{\centers,1\normalbase} \hgridcenters^+ \big)
    \end{bmatrix}
    \begin{pmatrix}
        \mask{\potgrid} \\
        -\jump{\potpointnormal}
    \end{pmatrix} =
    \begin{pmatrix}
        \sourcegrid + \mask{\bcgrid}\\
        \potpointsurf
    \end{pmatrix} \,.
    \label{eq:poisson_new_saddle_system}
\end{align}
The solution for this system is computed by performing an analytical block-LU factorization of the matrix in \eqref{eq:poisson_new_saddle_system} and pre-multiplying both sides by the inverse of the lower triangle matrix in the factorization (see \ref{app:saddlesystems}). The outcome of this process yields a sequence of equations, whose solution exactly satisfies \eqref{eq:poisson_new_saddle_system}:
\begin{align}
    \lapgrid \mask{\potgrid}^* &= (\sourcegrid + \mask{\bcgrid}),
    \label{eq:poisson_new_intermediate_potgrid_equation} \\
    \schur \jump{\potpointnormal} &= -\big( \potpointsurf - \interpbody_{\centers} \mask{\potgrid}^* \big), \label{eq:poisson_new_force_equation} \\
    \mask{\potgrid} &= \mask{\potgrid}^* + \invlapgrid \big( \spacetransform{\faces}{\centers}\regds_{\faces,n} \diag(\normalvec)^2 + \divgrid \regds_{\faces,1n} \diag(\normalvec) \big) \jump{\potpointnormal}, \label{eq:poisson_new_potgrid_equation}
\end{align}
where the Schur complement of the system~\eqref{eq:poisson_new_saddle_system} is
\begin{equation}
    \schur = -\diag \big( \interpbody_{\centers,1n} \hgridcenters^+ \big) - \interpbody_{\centers} \invlapgrid \big( \spacetransform{\faces}{\centers}\regds_{\faces} \diag(\normalvec)^2 + \divgrid \regds_{\faces,1 \normalbase} \diag(\normalvec) \big).
    \label{eq:poisson_new_schur}
\end{equation}

This solution strategy is commonly used in immersed methods. For example, \citet{taira_immersed_2007} use it to solve the prototypical IB formulation of the Navier--Stokes equations for rigid bodies in velocity-pressure formulation and streamwise-vorticity formulation~\citep{colonius_fast_2008}, and they referred to their method as the \emph{immersed boundary projection method} (IBPM), after the projection method for solving the incompressible Navier--Stokes equations~\citep{perot1993analysis}. Analogous projection-based solution processes exist formulated for the Poisson's equation~\citep{goza_accurate_2016, beckers_planar_2022}. A similar strategy was also used earlier in several versions of the immersed interface method (IIM) for Poisson-type problems~\citep{li_fast_1998, wiegmann_explicit_1998} and Navier–Stokes problems in a streamfunction–vorticity formulation~\citep{li_immersed_2001}.
Similar to the IIM, the unknown jumps in our method (specifically the first derivative in the current Dirichlet Poisson problem) appear in the constraint equation, which fills in the lower right-hand block of the system~\eqref{eq:poisson_new_saddle_system}.

Even though both the IBPM and the proposed method use the Schur complement reduction method, their conditioning is very different. A typical continuous-forcing IB formulation of the governing equations, such as the IBPM, would produce a linear system of equations analogous to that of \eqref{eq:poisson_new_saddle_system}, but with a zero matrix in the $(2,2)$-block of the left-hand side matrix \cite{goza_accurate_2016}. \citet{goza_accurate_2016} show that the block-LU process produces a Schur complement $\tilde{\schur}$ that represents the discrete analogue of an ill-posed Fredholm integral equation of the first kind. On the other hand, the new Schur complement $\schur$ given by~\eqref{eq:poisson_new_schur} represents the discrete analogue of a well-posed Fredholm integral equation of the second kind, which \citet{li_fast_1998} and \citet{wiegmann_explicit_1998} also identified in their IIM formulation. The key reason behind the improved conditioning of the Schur complement in the proposed IB projection formulation is the nonzero bottom diagonal term. This new term arises by explicitly incorporating the connection between the unknown IB source term and the jump in solution gradients near the interface. This connection, encoded into the equations via a Taylor series about the immersed interface, produces the nonzero diagonal matrix block. This term shifts the rapidly decaying eigenvalues---associated with the smoothing of high-frequency components in the Schur complement involving the regularization, inverse Laplacian, and interpolation operations---to the nonzero eigenvalues of $\diag \big( \interpbody_{\centers,1n} \hgridcenters^+ \big)$, thereby rendering the inverse bounded and the system well-conditioned. This approach to computing more accurate surface stresses addresses the ill-conditioning through a formal and careful formulation of the governing equations, using Taylor series near the interface. It accordingly avoids an aesthetic post-processing procedure that leaves the underlying poorly conditioned problem unresolved \cite{goza_accurate_2016}, or introducing heuristic parameters in a Tikhonov-style regularization procedure \cite{yu_regularized_2024}. 

The proposed method only uses operations that preserve the elements of the prototypical IB methods. The new operations involve the distance-weighted regularization and interpolation operators, which involve only an additional multiplication by the distance between the body point and grid points inside the support of the discrete DDF compared to the IB operators in the prototypical formulation. This distance is also used in the prototypical formulation as the argument of the discrete DDF and is therefore available at no extra cost. The proposed method also requires the computation of the indicator fields using Eq.~\eqref{eq:poisson_hgrid}, which only involves one extra regularization operation and Poisson solve. The overall algorithmic complexity of the proposed method is therefore only slightly larger than that of the prototypical formulation when ignoring the solution cost of the most expensive linear system in the LU-process, associated with the Schur complement. When accounting for the improved conditioning, the proposed method is potentially faster than the prototypical formulation, especially if an iterative solver is used to solve the Schur complement system.
We will now apply both our proposed and a typical projection-type continuous-forcing IB method for this Dirichlet Poisson problem to a 1D and 2D example to demonstrate the differences in accuracy and conditioning.

\subsection{Results for a 1D example}

To build intuition for how the proposed approach provides higher-order accuracy, we first focus on a 1D example with a constant source term and homogeneous Dirichlet boundary and interface conditions. The interface is represented by a single IB point with its normal pointing in the positive $x$-direction. Accordingly, we consider the subdomains $[x_L,x_\surf]$ and $[x_\surf,x_R]$ to be the interior and exterior regions, respectively. In this setting, Poisson's equation reduces to a second-order ordinary differential equation:
\begin{gather}
    \der{\pot^+}{x}{2} = \source \quad x \in [x_L,x_\surf],
    \\
    \der{\pot^-}{x}{2} = \source \quad x \in [x_\surf,x_R],
    \\
    \pot^-(x_L) = 0, \quad \pot^+(x_R) = 0, \quad \pot^\pm(x_\surf) = 0.
\end{gather}
The composite analytical solution for this problem is
\begin{equation}
    \pot = 
    \begin{cases}
        \displaystyle (x - x_L) (x - x_\surf) \source /2, & x \in [x_L,x_\surf],\\
        \displaystyle (x - x_\surf) (x - x_R) \source /2, & x \in [x_\surf,x_R], \\
    \end{cases}
\end{equation}
which has a jump in its first derivative $[\inlder{\pot}{x}{}]_\surf = -(x_R - x_L) \source / 2$ at $x = x_\surf$. In our example we use the numerical values $x_L = 0$, $x_R = 2$, and $\source(x) = -4$, which results in a first-derivative jump $[\inlder{\pot}{x}{}]_\surf = 4$.

\begin{figure}
    \centering
    \tikzexternaldisable
  \tikzsetnextfilename{figures/1D_example_diagram}%
  \input{figures/1D_example_diagram.tikz}%

    \tikzexternalenable
    \caption{Comparison of three different IB solutions on 16 grid cells for the 1D Dirichlet Poisson example. (Left) Solution of the prototypical IB equation~\eqref{eq:poisson_classic_IB} with the analytical forcing strength (green, stars) and the forcing strength obtained using Eq.~\eqref{eq:poisson_old_force_equation} (red, diamonds). (Right) Solution (blue, triangles) of the proposed IB equation~\eqref{eq:poisson_composite_TS} with the forcing strength obtained using Eq.~\eqref{eq:poisson_new_force_equation}. The exact solution is shown in gray with filled circles. The symbol $\epsilon_i = \big( \potgrid_i - \pot(\xgrid[(i)]) \big) / || \pot(\xgrid) ||_\infty$ denotes the local relative error with respect to the analytical solution.}
    \label{fig:1D_accuracy_diagram}
\end{figure}
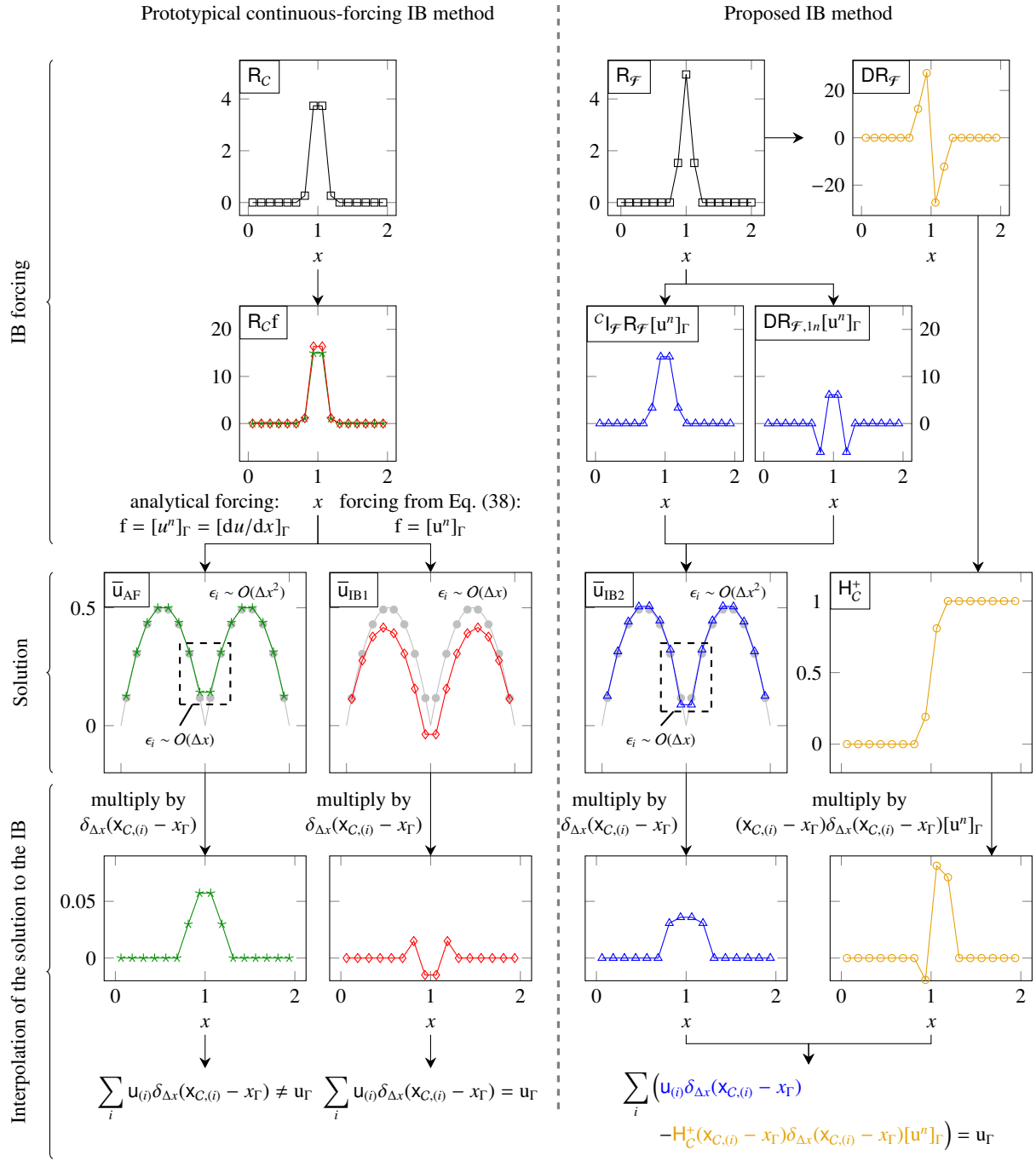

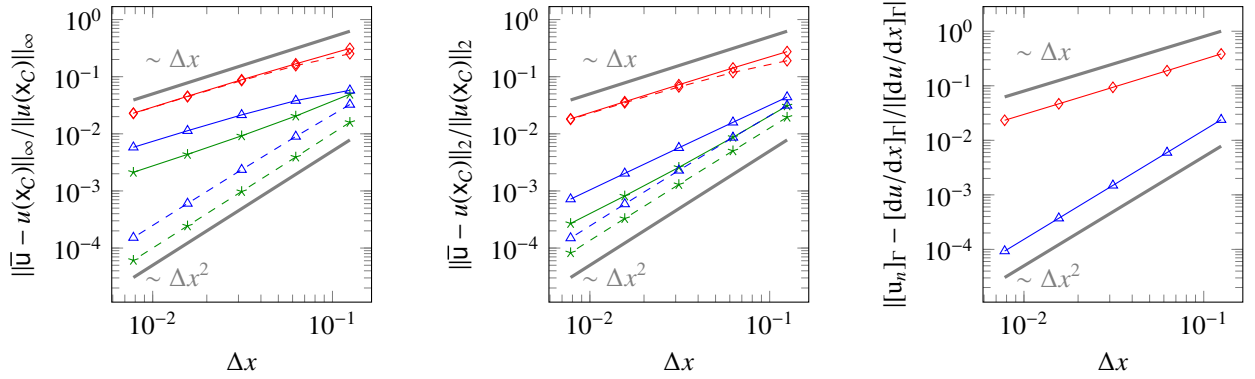
\begin{figure}
    \centering
  \tikzsetnextfilename{figures/1D_poisson_error_Linf_norm}%
  \input{figures/1D_poisson_error_Linf_norm.tikz}%

    \hfill
  \tikzsetnextfilename{figures/1D_poisson_error_L2_norm}%
  \input{figures/1D_poisson_error_L2_norm.tikz}%

    \hfill
  \tikzsetnextfilename{figures/1D_poisson_forcing_error}%
  \input{figures/1D_poisson_forcing_error.tikz}%

    \caption{Error of the numerical solution computed over the entire domain (solid lines) and over the domain excluding cells within the support of the DDF (dashed lines) using the infinity norm (left) and the 2-norm (center) and the error of the forcing strength (right) for the 1D Dirichlet Poisson example. The line colors and markers indicate the error using the prototypical continuous-forcing IB equation~\eqref{eq:poisson_classic_IB} with the analytical forcing (green, stars), the same equation~\eqref{eq:poisson_classic_IB} with the forcing strength obtained through Eq.~\eqref{eq:poisson_old_force_equation} (red, diamonds), and the proposed IB equation~\eqref{eq:poisson_composite_TS} with the forcing strength obtained through Eq.~\eqref{eq:poisson_new_force_equation} (blue, triangles). First- and second-order error scales are also shown (gray, solid).}
    \label{fig:1D_poisson_error}
\end{figure}

In Fig.~\ref{fig:1D_accuracy_diagram} we compare three different solutions and their forcing: the solution $\mask{\potgrid}_{\mathrm{AF}}$ of the prototypical continuous-forcing IB equation \eqref{eq:poisson_classic_IB} using the analytical first derivative jump as a prescribed IB forcing, the solution $\mask{\potgrid}_{\mathrm{IB1}}$ of the same equation, but where we solve for the first-derivative jump such that the solution satisfies the prototypical continuous-forcing IB constraint $\interpbody_\centers \potgrid = \potpointsurf$, and the solution $\mask{\potgrid}_{\mathrm{IB2}}$ of the proposed IB system~\eqref{eq:poisson_new_saddle_system}. In Fig.~\ref{fig:1D_poisson_error} we compare the relative error of each solution with respect to the analytical solution $\pot(\xgrid)$ evaluated at the same grid points, with and without including the error near the interface, as well as the error of the first-derivative jump at the IB (i.e., the IB forcing).

As expected from the analysis of \citet{tornberg_numerical_2004}, the solution $\mask{\potgrid}_{\mathrm{AF}}$, which uses the analytical first-derivative jump as its forcing, is second-order accurate away from the interface and only first-order accurate within the support of the DDF. However, applying the discrete interpolation~\eqref{eq:disc_interp} does not recover the prescribed interface value $\pot_{\surf}$. Because smoothed DDFs cause the IB method to overestimate the solution near a positive first-derivative jump (and underestimate it for a negative one), the DDF-weighted interpolation yields values that in this case are biased high, leading to an overestimation of the interpolated interface value. Moreover, in cases of interest to us the IB forcing will not generally be known---it must be computed to enforce the constraint---so this approach used for illustrative purposes is inaccessible.

The solution $\mask{\potgrid}_{\mathrm{IB1}}$ uses the first-derivative jump computed from
\begin{equation}
    \tilde{\schur} \jump{\potpointnormal} = -\big( \potpointsurf - \interpbody_{\centers} \invlapgrid (\sourcegrid + \mask{\bcgrid}) \big)
    \label{eq:poisson_old_force_equation} 
\end{equation}
with $\tilde{\schur} = \interpbody_{\centers} \invlapgrid \regds_{\centers}$. This solution comes from an analogous block-LU factorization process to that defined for the proposed solution above, using the prototypical continuous-forcing IB equation~\eqref{eq:poisson_classic_IB} coupled with $\interpbody_\centers \potgrid = \potpointsurf$. The resulting IB force (i.e., jump in solution gradient), is of larger magnitude near the interface than the analytical one. This outcome follows from the need to satisfy the discrete constraint: for a positive jump, the method must underestimate the solution within the DDF support, and for a negative jump it must overestimate it. As a result, the solution near the interface consists either of values equal to the imposed interface value (possible with a two-point DDF) or of mixed values that balance to yield the correct DDF-weighted average. In the present example, the positive jump forces an underestimation of the solution, producing a larger effective first-derivative jump. This inability of the interpolation to represent a first-derivative jump is what leads to the characteristic over- and underestimation and ultimately limits both the interpolation and the overall method to first-order accuracy~\cite{beyer_analysis_1992}.

Lastly, the third solution $\mask{\potgrid}_{\mathrm{IB2}}$, obtained with our proposed IB method, is still first-order accurate within the support of the (interpolated) DDF, but is second-order accurate away from the IB while simultaneously satisfying the discrete constraint equation~\eqref{eq:disc_pot_interp_TS}, which involves the indicator field obtained from the same DDF that is used in the governing equation. The second-order accuracy away from the IB could also be obtained by applying the analytical forcing using only the $\spacetransform{\faces}{\centers} \regds_{\faces} \forcepoint$ term, because the term $\divgrid \regds_{\faces, 1 \normalbase} \forcepoint$ has no effect on the solution away from the IB. Its only purpose is to adapt the solution near the IB such that the Taylor series used to obtain~\eqref{eq:disc_pot_interp_TS} are valid. 
While not shown here, this extra term ensures that $\mask{\potgrid}_{\mathrm{IB2}}$ converges with second-order accuracy to $\hgridcenters^+ \had \pot^+(\xgrid[\centers]) + \hgridcenters^- \had \pot^-(\xgrid[\centers])$ over the entire domain (including within the DDF support), which differs from the exact composite solution $\pot$, which is defined with sharp indicator functions in Eq.~\eqref{eq:composite_analytic_pot}. 

\subsection{Results for a 2D example}
\label{subsec:poisson_2d_example}

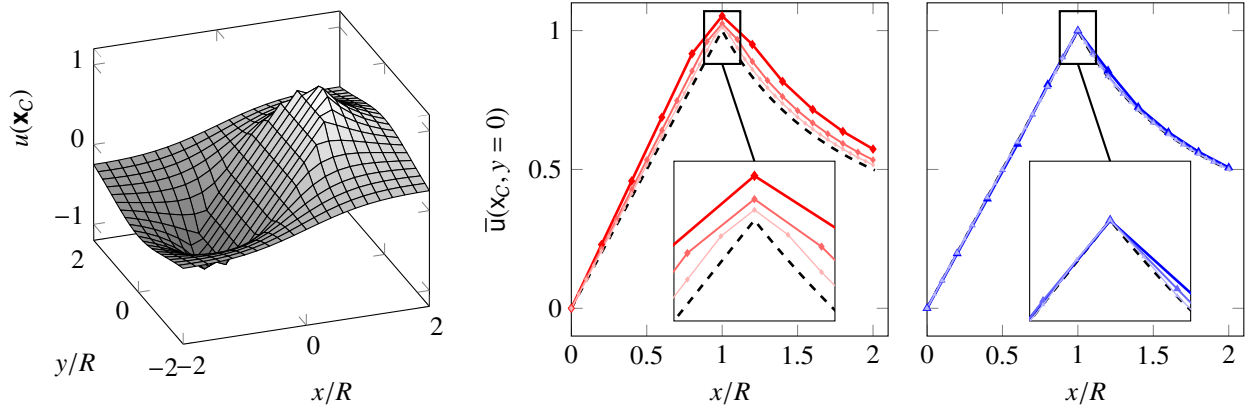
\begin{figure}
    \centering
  \tikzsetnextfilename{figures/2D_poisson_mesh}%
  \input{figures/2D_poisson_mesh.tikz}%

    \hfill
  \tikzsetnextfilename{figures/2D_poisson_solution_original}%
  \input{figures/2D_poisson_solution_original.tikz}%

    \hspace{0.0cm}
  \tikzsetnextfilename{figures/2D_poisson_solution_new}%
  \input{figures/2D_poisson_solution_new.tikz}%

    \caption{(Left) Analytical solution to the 2D Poisson equation on an unbounded domain subject to a sinusoidal Dirichlet interface condition on a circle centered at the origin with radius $R$, evaluated at the cell centers of a grid with $\dx / R = 0.2$. (Center and right) Comparison of the numerical solution from the proposed method (blue, dashed) with the solution from the prototypical continuous-forcing IB method (red, dotted) on the line along $y=0$ near the interface for $\ds / \dx \approx 1$ and three different grid cell sizes $\dx/R$: 0.2 (darkest, thick), 0.1, and 0.05 (lightest, thin). The exact solution is also shown (black, dashed).}
    \label{fig:poisson_solution}
\end{figure}

Now we apply the previous IB methods to the 2D example of a circle with radius $R$ and outward pointing normals centered at the origin in an unbounded domain, where the boundary condition~\eqref{eq:poisson_bc} is replaced by $\pot(\x) \to 0$ as $\x \to \infty$. In our discrete method, we can treat this problem using the lattice Green's function for unbounded domains on a bounded grid~\cite{liska_fast_2016}. The prescribed interface condition on $\surf$, described by a single coordinate $\xi$, is $\pot_{\surf}(\xi) = X(\xi)$. The analytical interior solution (inside the circle) is then $\pot^-(\x) = x$ and the exterior solution (outside the circle) in polar coordinates is $\pot(r,\theta) = R^2/r \cos(\theta)$. Fig.~\ref{fig:poisson_solution} shows the composite analytical solution. The analytical normal-derivative jump is $\jump{\pot^{\normalbase}}(\xi) = -2 \cos\big(\theta(\xi)\big)$.

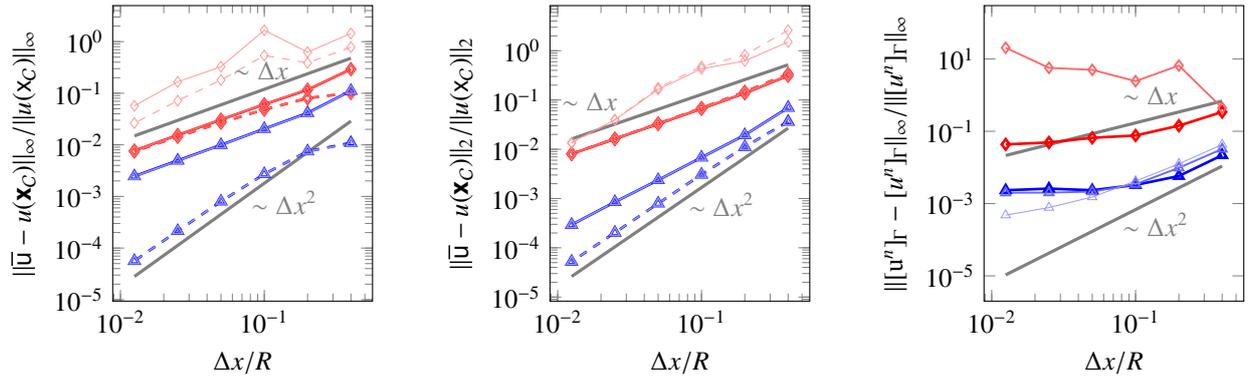
\begin{figure}
    \centering
  \tikzsetnextfilename{figures/2D_poisson_error_Linf_norm_dsdx_variation}%
  \input{figures/2D_poisson_error_Linf_norm_dsdx_variation.tikz}%

    \hfill
  \tikzsetnextfilename{figures/2D_poisson_error_L2_norm_dsdx_variation}%
  \input{figures/2D_poisson_error_L2_norm_dsdx_variation.tikz}%

    \hfill
  \tikzsetnextfilename{figures/2D_poisson_forcing_error_Linf_norm_dsdx_variation}%
  \input{figures/2D_poisson_forcing_error_Linf_norm_dsdx_variation.tikz}%

    \caption{Infinity norm (left) and 2-norm (center) of the numerical solution error computed over all cell centers (solid lines) and over all cell centers excluding those within the support of the DDF (dashed lines) and infinity norm of the forcing strength error (right) for the 2D Dirichlet Poisson example using the prototypical continuous-forcing IB method (red, diamonds) and the proposed method (blue, triangles) for three approximate values of $\ds/\dx$: 1.3 (darkest, thick), 0.7, and 0.1 (lightest, thin). First- and second-order error scales are also shown (gray, solid). The forcing error values of the prototypical formulation for $\ds/\dx \approx 0.1$ are not included to limit the upper range of the plot.}
    \label{fig:2D_poisson_error}
\end{figure}

We compute and evaluate the error of the numerical solution and the interface forcing term, for various discretizations of $\domain = \left\{ (x,y) \in \R^2 \;\middle|\; -2R \le x \le 2R,\; -2R \le y \le 2R \right\}$. These numerical solutions are computed both using the proposed IB system~\eqref{eq:poisson_new_saddle_system} and the prototypical continuous-forcing IB equation~\eqref{eq:poisson_classic_IB} coupled with $\interpbody_\centers \potgrid = \potpointsurf$, which results in the forcing equation~\eqref{eq:poisson_old_force_equation}. To quantify convergence, numerical solutions for both cases are computed using discretizations with several grid sizes. To assess the impact of problem conditioning on the solution, we also utilize a range of values for the body-to-grid spacing ratio $\ds/\dx$. In all cases, the grid spacing is uniform throughout the domain and identical in the $x$- and $y$-directions.

We show the solutions obtained with the prototypical method and with our proposed method on the middle and right panels of Fig.~\ref{fig:poisson_solution}, respectively. Figure~\ref{fig:2D_poisson_error} shows the associated relative error, with respect to the true solution, of the solution $\mask{\potgrid}$ and the forcing $\jump{\potpointnormal}$ versus the grid cell size for different values of $\ds/\dx$. 

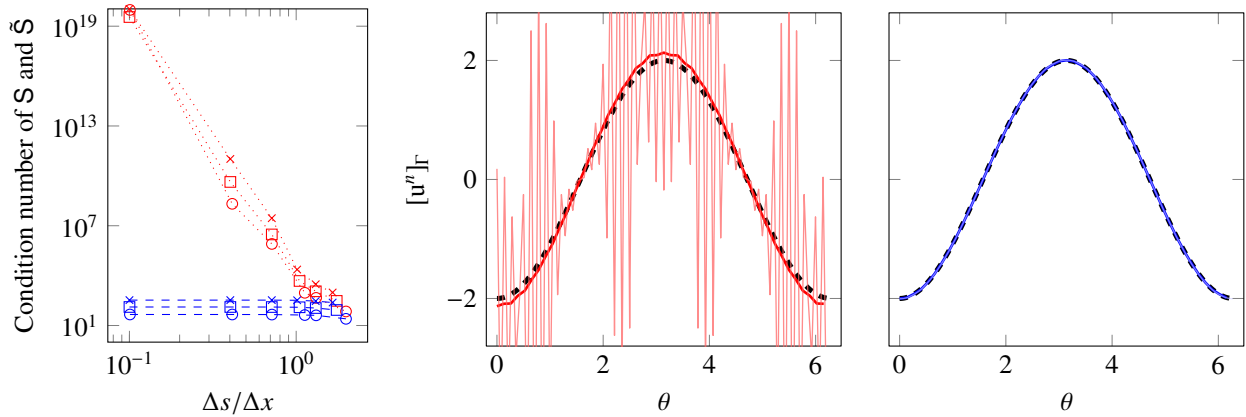
\begin{figure}
    \centering
  \tikzsetnextfilename{figures/2D_poisson_condition_number}%
  \input{figures/2D_poisson_condition_number.tikz}%

    \hfill
  \tikzsetnextfilename{figures/2D_poisson_forcing_original}%
  \input{figures/2D_poisson_forcing_original.tikz}%

    \hfill
  \tikzsetnextfilename{figures/2D_poisson_forcing_new}%
  \input{figures/2D_poisson_forcing_new.tikz}%

    \caption{(Left) Condition number of the Schur complement  versus the surface-to-grid spacing ratio for the proposed method ($\schur$, blue, dashed) and the prototypical continuous-forcing equation and constraint ($\tilde{\schur}$, red, dotted) using three different grid cell sizes $\dx/R$: 0.2 (circles), 0.1 (squares), and 0.05 (crosses). (Center and right) Forcing solution obtained with the prototypical continuous-forcing IB method (red, dotted) and the proposed method (blue, dashed) for $\dx/R= 0.2$ and three approximate values of $\ds/\dx$: 1.3 (darkest, thick), 0.7, and 0.1 (lightest, thin). The solution from the prototypical formulation for $\ds/\dx \approx 0.1$ cannot be meaningfully represented on the chosen $y$-axis scale and is therefore not included. The exact solution is also shown (black, dashed).}
    \label{fig:2D_poisson_forcing}
\end{figure}

The solution variable $\mask{\potgrid}$ provides a convergent approximation to the true solution for both IB methods, as expected. In the case when the prototypical continuous-forcing IB method is used, the error in this solution variable converges at first order with respect to both the infinity and $2$-norms. By contrast, the proposed approach provides a solution that converges at first order in the global infinity norm, second order in the infinity norm omitting spatial points within the support of the DDF, and between first and second order in the two norm. This improved accuracy is a result of using the Taylor series to obtain the new interpolation formula in~\eqref{eq:disc_pot_interp_TS} and in the governing equation~\eqref{eq:poisson_composite_TS}, so that the interpolation accuracy of~\eqref{eq:disc_pot_interp_TS} can be realized.

The effect of $\Delta s/\Delta x$ is distinct for the two methods. Using the prototypical interface constraint, decreasing $\Delta s/\Delta x$ to the smallest considered value of 0.1 leads to a poorer approximation in the solution variable $\mask{\potgrid}$ (though it still converges with $\dx$ at first order). This degraded solution accuracy is reflected in progressively poor approximations to the IB forcing strength. Figure~\ref{fig:2D_poisson_error} shows that this forcing does not converge for the IB method that uses a prototypical constraint condition. The left and middle columns of Figure~\ref{fig:2D_poisson_forcing} clarify this issue. For decreasing values of $\Delta s/\Delta x$, when an IB method is built using the prototypical IB constraint, the condition number of the Schur complement $\tilde{\schur}$ (left panel), defined as the ratio of its largest to smallest singular value, becomes sufficiently high to introduce significant numerical error in the solution. This poor conditioning and its corresponding effect to produce highly oscillatory IB forcing along the surface (middle panel) is a direct consequence of the ill-posed nature of the underlying integral equation that the Schur complement discretizes when the traditional IB constraint is used \cite{goza_accurate_2016,yu_regularized_2024}. Note that the IB forcing could not even be pictured for the smallest value of $\Delta s/\Delta x$, because of how large amplitude the oscillations are for that case.

By contrast, the proposed method, which modifies the IB constraint using Taylor series, produces a solution $\mask{\potgrid}$ that is essentially independent of the ratio $\Delta s/\Delta x$. Figure~\ref{fig:2D_poisson_forcing} demonstrates that this robust computation of the solution variable is associated with a more reliably well-conditioned Schur complement $\schur$ (left panel). The figure also shows IB forcing along the surface that is non-oscillatory and that effectively approximates the true forcing. In fact, the right panel of Figure~\ref{fig:2D_poisson_error} indicates that the IB forcing term might even be converging to the true forcing along the interface with fine enough ratios of $\Delta s/\Delta x$. We leave a full exploration of the proposed approach as a pathway to compute convergent surface stresses to future work. At such fine interface resolutions, the Schur complement becomes increasingly high-dimensional and high-performance linear solvers would need to be employed to fully explore the topic. For now, we observe that accurate, smooth IB forcing can be computed via well-conditioned linear systems that are strongly robust to the ratio $\Delta s/\Delta x$. Moreover, these IB forces are computed by fundamentally changing the nature of the integral equation from an ill-posed to well-posed problem, without introducing heuristic regularization parameters. This functionality may be especially useful in problems with moving and compliant bodies, where maintaining a perfectly constant grid spacing ratio is impractical.

\section{Treatment of the incompressible Navier--Stokes equations}
\label{sec:navier-stokes}

In this section we derive our new formulation for the non-dimensionalized, incompressible Navier--Stokes equations for the velocity $\vel^\pm(\x,t) \in \R^3$ and pressure $p^\pm(\x,t) \in \R$ in the interior and exterior regions with uniform properties,
\begin{gather}
    \frac{\partial \vel^\pm(\x,t)}{\partial t} + \div \big( \vel^\pm(\x,t) \vel^\pm(\x,t) \big)  = -\grad p^\pm(\x,t) + \frac{1}{\Re} \lap \vel^\pm(\x,t), \quad \x \in \domain^\pm
    \label{eq:ns_momentum}
    \\
    \div \vel^\pm(\x,t) = 0, \quad \x \in \domain^\pm 
    \label{eq:ns_continuity}
    \\
    \vel^{\pm}(\X(\surfc,t),t) = \vel_\surf(\surfc,t), \quad \surfc \in \surf(t),
    \label{eq:ns_ibc}
\end{gather}
where $\Re$ is the Reynolds number, $\vel_\surf(\surfc,t)$ is the prescribed velocity on the IB at the Lagrangian coordinate $\X(\surfc)$ and time $t$, and we assume that appropriate boundary conditions are provided on $\partial \domain^\pm \setminus \surf(t)$.

Following again the approach of~\citet{eldredge_method_2022}, we define the composite velocity and pressure
\begin{gather}
    \vel \coloneqq \heavi^+ \vel^+ + \heavi^- \vel^-,
    \label{eq:composite_analytic_vel}
    \\
    p \coloneqq \heavi^+ p^+ + \heavi^- p^-,
    \label{eq:composite_analytic_p}
\end{gather}
and we apply the relevant product rules and indicator function identities to obtain (see \ref{app:ns_ib_derivation}):
\begin{gather}
    \frac{\partial \vel}{\partial t} + \div ( \vel \vel ) 
    - \vecforce_{\inlder{\vel}{t}{}}
    - \vecforce_{\div ( \vel \vel )}
    =
    -\grad p + \frac{1}{\Re} \lap \vel
    + \vecforce_{\grad p}
    - \frac{1}{\Re} \vecforce_{\lap \vel}, \quad \x \in \domain,
    \label{eq:ns_momentum_composite}
    \\
    \div \vel 
    - \scaforce_{\div \vel}
    =
    0, \quad \x \in \domain,
    \label{eq:ns_continuity_composite}
    \\
    \vel(\X(\surfc,t),t) = \vel_\surf(\surfc,t), \quad \surfc \in \surf(t),
    \label{eq:ns_ibc_composite}
\end{gather}
where
\begin{align*}
    \vecforce_{\grad p}
    &= 
    \grad \heavi^+ ( p^+ - p^- ),
    \\
    \vecforce_{\lap \vel} 
    &= 
    \grad \heavi^+ \cdot ( \grad \vel^+ - \grad \vel^- ) + \div \big( \grad \heavi^+ (\vel^+ - \vel^- ) \big),
    \\
    \vecforce_{\div ( \vel \vel )}
    &= \grad \heavi^+ \cdot (\vel^+ \vel^+ - \vel^- \vel^-),
    \\
    \vecforce_{\inlder{\vel}{t}{}}
    &= \big( \grad \heavi^+ \cdot \Xdot \big) (\vel^+ - \vel^-),
    \\
    \scaforce_{\div \vel} 
    &=  \grad \heavi^+ \cdot ( \vel^+ - \vel^- ).
\end{align*}
In the current case, we prescribe no-slip boundary conditions on the IB for both the exterior and interior solutions. As a result, $\vel^+ \big( \X(\surfc),t \big) = \vel^{-} \big( \X(\surfc), t \big)$ for all $\surfc \in \surf$ and $t$. Discontinuous solutions across the interface can be accommodated directly into the formulation. 

In Section~\ref{subsec:ns_discrete_formulation} and Section~\ref{subsec:ns_regularization}, we will follow an analogous set of steps to the Poisson problem to arrive at the discretized version of \eqref{eq:ns_momentum_composite}--\eqref{eq:ns_ibc_composite}. To start, we combine the discrete Navier--Stokes equations for the interior and exterior fields and recast them into the governing equations for the discrete composite field, which contain forcing terms involving the gradient of the indicator function multiplying the difference of the exterior and interior solutions (or variants thereof), similar to the five forcing terms above. Then, in each of the forcing terms, we replace the indicator function gradient with the normals regularized using the smoothed DDFs (Eq.\eqref{eq:peskin_heavi}), which are only nonzero in a small region near the interface. 
This substitution, in turn, allows us to use Taylor series expressed using surface quantities to represent the difference terms, whose only parts that matter are near the interface where the DDFs are nonzero. Here, the goal is again to obtain an expression for the governing equations that is consistent with the formulation of the discrete composite solution, such that we can improve the interpolation accuracy based on the knowledge of the indicator function behavior near the interface.

The relevant surface quantities that appear in the leading-order terms of our discrete formulation are the jump in the normal derivative of the velocity and the jump in pressure over the interface. We treat these as the Lagrange multipliers that enforce the discrete analogue of \eqref{eq:ns_ibc_composite}, which requires one additional equation that links the pressure field with the pressure jump. These two equations are formulated in Section~\ref{subsec:ns_interpolations}, where we use a Taylor series expansion to increase the order of accuracy of the velocity interpolation, similar to our approach in the Poisson problem. The governing equations and the two constraint equations are then solved together in a projection-based approach, described in Section~\ref{subsec:ns_projection_method}.


\subsection{Discrete problem formulation for a composite solution}
\label{subsec:ns_discrete_formulation}

As in the Poisson problem, we derive the discrete formulation starting from the discrete composite velocity field $\mask{\velgrid} = (\mask{\velgridcomponent}_x, \mask{\velgridcomponent}_y, \mask{\velgridcomponent}_z) \in \R^{\faces}$ and pressure field $\mask{\pgrid} \in \R^{\centers}$, defined as
\begin{gather}
    \mask{\velgrid} \coloneqq \hgridfaces^+ \had \velgrid^+ + \hgridfaces^- \had \velgrid^-
    \\
    \mask{\pgrid} \coloneqq \hgridcenters^+ \had \pgrid^+ + \hgridcenters^- \had \pgrid^- ,
\end{gather}
where $\velgrid^\pm$ and $\pgrid^\pm$ are the discrete approximations to the exterior and interior velocity and pressure fields, respectively, and $\hgridfaces^\pm \coloneqq \spacetransform{\faces}{\centers} \hgridcenters^\pm \in \R^{\faces}$ is the interpolation of $\hgridcenters^\pm$ to each cell face.

We seek the solution for $\mask{\velgrid}$ and $\mask{\pgrid}$ such that the interior and exterior velocity and pressure fields satisfy the following spatially-discretized, non-dimensionalized, incompressible Navier--Stokes equations in their respective domains (assuming uniform properties in the entire domain),
\begin{gather}
    \der{\velgrid^\pm}{t}{} + \convecgrid \left( \velgrid^\pm \right) = - \gradgrid \pgrid^\pm + \frac{1}{\Re} \lapgrid_{\faces} \velgrid^\pm + \bcgridvec^\pm_1,
    \\
    \divgrid \velgrid^\pm = \bcgrid^\pm_2,
\end{gather}
where $\convecgrid$ is the discrete approximation to the convective term, $\lapgrid_{\faces} : \R^{\faces} \mapsto \R^{\faces}$ is the vector Laplacian for cell face data (see \ref{app:differential_operators} for details), and $\bcgridvec^\pm_1$ and $\bcgrid^\pm_2$ represent the boundary-condition terms that arise from discretizing the operators in the momentum and continuity equations, respectively, while enforcing any provided boundary conditions.

The corresponding equations for $\mask{\velgrid}$ and $\mask{\pgrid}$ can be found by first multiplying the Navier--Stokes equations for the interior and exterior fields by their respective indicator functions and summing them together, 
\begin{gather}
    \hgridfaces^+ \had \left(\der{\velgrid^+}{t}{} + \convecgrid \left( \velgrid^+ \right)\right) + \hgridfaces^- \had \left(\der{\velgrid^-}{t}{} + \convecgrid \left( \velgrid^- \right)\right)
    = 
    - \hgridfaces^+ \had \gradgrid \pgrid^+  - \hgridfaces^- \had  \gradgrid \pgrid^- + \frac{1}{\Re} \left( \hgridfaces^+ \had \lapgrid_{\faces} \velgrid^+ +\hgridfaces^- \had \lapgrid_{\faces} \velgrid^- \right) + \mask{\bcgridvec}_1 ,
    \\
    \hgridcenters^+ \had \divgrid \velgrid^+ + \hgridcenters^- \had \divgrid \velgrid^-
    = \mask{\bcgrid}_2,
\end{gather}
where $\mask{\bcgridvec}_1 \coloneqq \hgridfaces^+ \had \bcgridvec^+_1 + \hgridfaces^- \had \bcgridvec^-_1$ and $\mask{\bcgrid}_2 \coloneqq \hgridcenters^+ \had \bcgrid^+_1 + \hgridcenters^- \had \bcgrid^-_2$. Then, we replace the exterior and interior terms with matching operators by that operator acting on the composite field minus a forcing term that involves the gradient or time derivative of an indicator field from the relations~\eqref{eq:composite_scalar_grad}, \eqref{eq:composite_vector_div}, \eqref{eq:composite_vector_lap}, \eqref{eq:composite_vector_time_der}, and \eqref{eq:composite_vector_conv} in~\ref{app:differential_operators_on_composite_fields}, resulting in the discrete analogue of~\eqref{eq:ns_momentum_composite}--\eqref{eq:ns_continuity_composite}:
\begin{equation}
    \der{\mask{\velgrid}}{t}{} 
    + \convecgrid ( \mask{\velgrid} )
    - \forcegridvec_{\inlder{\mask{\velgrid}}{t}{}}
    - \forcegridvec_{\convecgrid ( \mask{\velgrid} )}
    = 
    - \gradgrid \mask{\pgrid}
    + \frac{1}{\Re} \lapgrid_{\faces} \mask{\velgrid}
    + \forcegridvec_{\gradgrid \mask{\pgrid}}
    - \frac{1}{\Re} \forcegridvec_{\lapgrid_{\faces} \mask{\velgrid}} + \mask{\bcgridvec}_1 ,
\end{equation}
\begin{equation}
    \divgrid \mask{\velgrid} - \forcegrid_{\divgrid \mask{\velgrid}} = \mask{\bcgrid}_2 ,
\end{equation}
where
\begin{align*}
    \forcegridvec_{\gradgrid \mask{\pgrid}} 
    &= 
    \gradgrid \hgridcenters^+ \circ \spacetransform{\centers}{\faces} ( \pgrid^+ - \pgrid^- ),
    \\
    \forcegridvec_{\lapgrid_{\faces} \mask{\velgrid}} 
    &= 
    \spacetransform{\tensorspace}{\faces} \Big( \big(\spacetransform{\faces}{\tensorspace} \gradgrid \hgridcenters^+ \big)^\top \had  \big( \gradgrid_{\faces} \velgrid^+ - \gradgrid_{\faces} \velgrid^- \big) \Big)
    +
    \divgrid_{\tensorspace}\Big( \big(\spacetransform{\faces}{\tensorspace}\gradgrid \hgridcenters^+\big)^\top \had \spacetransform{\faces}{\tensorspace} \big( \velgrid^+ - \velgrid^- \big) \Big),
    \\
    \forcegridvec_{\convecgrid ( \mask{\velgrid} )}
    &= \spacetransform{\tensorspace}{\faces} 
        \bigg( 
            \gradgrid_{\faces} \hgridfaces^+
            \had 
            \Big( \big( \spacetransform{\faces}{\tensorspace}\velgrid^+ \big)^\top \had \spacetransform{\faces}{\tensorspace}  \velgrid^+  - \big( \spacetransform{\faces}{\tensorspace}\velgrid^- \big)^\top \had \spacetransform{\faces}{\tensorspace}  \velgrid^-  \Big)
        \bigg)
        \\
        &\quad-
        \divgrid_{\tensorspace} 
        \Big( 
            \hgridtensorspace^+ \had \hgridtensorspace^- \had \big( \spacetransform{\faces}{\tensorspace} ( \velgrid^+ - \velgrid^- ) \big)^\top \had \spacetransform{\faces}{\tensorspace} ( \velgrid^+ - \velgrid^- )  
        \Big)
        + \mathcal{O}(\dxvec^2),
    \\
    \forcegridvec_{\inlder{\mask{\velgrid}}{t}{}}
    &= \der{\hgridfaces^+}{t}{} \had ( \velgrid^+ - \velgrid^- ),
    \\
    \forcegrid_{\divgrid \mask{\velgrid}} 
    &= \spacetransform{\faces}{\centers} \big( \gradgrid \hgridcenters^+ \had ( \velgrid^+ - \velgrid^- ) \big) .
\end{align*}
In the previous terms, we introduced several operators acting between cell face data and discrete second-order tensor fields, such as the velocity gradient, defined on the space $\R^{\tensorspace}$ (see \ref{app:discrete_operations_and_tools} for their details). These are the discrete gradient for cell face data, $\gradgrid_{\faces} : \R^{\faces} \mapsto \R^{\tensorspace}$, the discrete divergence for second-order tensor data, $\divgrid_{\tensorspace} : \R^{\tensorspace} \mapsto \R^{\faces}$, and the interpolation and expansion/contraction from vectors to second-order tensors, $\spacetransform{\faces}{\tensorspace} : \R^{\faces} \mapsto \R^{\tensorspace}$, and vice versa, $\spacetransform{\tensorspace}{\faces} : \R^{\tensorspace} \mapsto \R^{\faces}$. 

The next step is to replace $\gradgrid \hgridcenters^+$ with $\regds_{\faces} \normalvec$ in each of the preceding terms, except $\forcegridvec_{\inlder{\mask{\velgrid}}{t}{}}$. Then, in the next section, we apply Taylor series to express the difference terms involving the exterior and interior solution using surface quantities.

\subsection{Regularization operations}
\label{subsec:ns_regularization}

Similar to our treatment of the Poisson problem, we approximate the variation of the pressure and velocity within the support of a discrete DDF by a truncated Taylor series expansion about the discrete surface point at which the DDF is centered. For the pressure we use only the zeroth-order term of~\eqref{eq:s_TS} and for the velocity and its gradient we use~\eqref{eq:v_TS} and \eqref{eq:Gv_TS}. We will make use of the notation $\velpointsurfplusminus$, $\velpointsurfnormalplusminus$, $\velpointsurftanoneplusminus$, and $\velpointsurftantwoplusminus$ to denote the value and normal and tangential derivatives of the exterior and interior velocity fields at the IB points, and we will use $\ppointsurfplusminus$ to denote the exterior and interior pressure at the IB points.
Because $\velpointsurfplus = \velpointsurfminus$, we have that $\jump{\velpoint} = 0$ and also $\jump{\velpointtanone} = 0$ and $\jump{\velpointtantwo} = 0$. We can then  formulate the following approximations for substitution into $\forcegrid_{\divgrid \mask{\velgrid}}$ and $\forcegridvec_{\gradgrid \mask{\pgrid}}$:
\begin{gather}
    \regds_{\faces} \normalvec \had ( \velgrid^+ - \velgrid^- )
    \approx \regds_{\faces, 1 \normalbase} \big( \normalvec \had \jump{\velpointnormal} \big),
    \label{eq:fdv_TS}
    \\
    \regds_{\faces} \normalvec \had \spacetransform{\centers}{\faces} ( \pgrid^+ - \pgrid^- )
    \approx
    \regds_{\faces} \big( \normalvec \had \spacetransform{\spoints}{\vpoints} \jump{\ppoint} \big).
    \label{eq:fgp_TS}
\end{gather}

Next, we define the space $\R^{\tpoints} \coloneqq \big(\R^{\spoints}\big)^9$ for second-order tensor data at the IB points and the operation $\outerproductoperator_{\vectorpoint[1]}: \R^{\vpoints} \mapsto \R^{\tpoints}$, which acts on vector-valued point data and produces the outer product with $\vectorpoint[1]$ as $\indicesbig{\outerproductoperator_{\vectorpoint[1]}(\vectorpoint[2])}{l} \coloneqq \vectorpoint[1,l] \vectorpoint[2,l]^\top$. If we then use the regularization of tensor-valued data $\regds_{(\spacetransform{}{}\faces)^\top}: \R^{\tpoints} \mapsto \R^{\tensorspace}$ and its normal-distance-weighted variant $\regds_{(\spacetransform{}{}\faces)^\top, 1 \normalbase} \in \R^{\tpoints} \mapsto \R^{\tensorspace}$ defined in~\ref{app:regularization}, we can approximate the expressions in the outermost parentheses of $ \forcegridvec_{\lapgrid_{\faces} \mask{\velgrid}}$ as
\begin{gather}
    \big(\spacetransform{\faces}{\tensorspace} \regds_{\faces} \normalvec \big)^\top \had  \big( \gradgrid_{\faces} \velgrid^+ - \gradgrid_{\faces} \velgrid^- \big) = \regds_{(\spacetransform{}{}\faces)^\top} \outerproductoperator_{\normalvec \had \normalvec} \big( \jump{\velpointnormal} \big),
    \label{eq:flv_term_1_TS}
    \\
    \big(\spacetransform{\faces}{\tensorspace} \regds_{\faces} \normalvec \big)^\top \had \spacetransform{\faces}{\tensorspace} \big( \velgrid^+ - \velgrid^- \big)
    \approx
    \regds_{(\spacetransform{}{}\faces)^\top, 1 \normalbase} \outerproductoperator_{\normalvec} \big( \jump{\velpointnormal} \big).
    \label{eq:flv_term_2_TS}
\end{gather}

The term $\forcegridvec_{\convecgrid \left( \mask{\velgrid} \right)}$---which is only nonzero within the support of the discrete DDF---can be ignored since its magnitude inside this support scales as $\mathcal{O}(\dxvec)$ while the other terms scale as $\mathcal{O}(1/\dxvec)$. This can be seen by considering that the Riemann sum~\eqref{eq:disc_reg} with the discrete DDF~\eqref{eq:ddf} scales as $\mathcal{O}(1/\dxvec)$, while the Taylor series expansion of the velocity fields difference scales as $\mathcal{O}(\dxvec)$ within the DDF support. A divergence or gradient operation increases the scaling by another factor of $\dxvec$. Since the first term of $\forcegridvec_{\convecgrid \left( \mask{\velgrid} \right)}$ involves the difference of the product of two velocity fields, which scales as $\mathcal{O}(\dxvec^2)$, the entire first term involving $\gradgrid_{\faces} \hgridfaces^+$ scales as $\mathcal{O}(\dxvec)$. Similarly, one can show that the second term of $\forcegridvec_{\convecgrid \left( \mask{\velgrid} \right)}$ scales as $\mathcal{O}(\dxvec)$. 

For the time derivative of the discrete indicator function, which appears in the term $\forcegridvec_{\inlder{\mask{\velgrid}}{t}{}}$, we directly discretize $\grad \heavi^+ \cdot \Xdot = \int (\normal \cdot \Xdot) \ddf (\x - \X) \dsurf$~\citep{juric_computations_1998, eldredge_method_2022} as $\regds_{\faces} (\dot{\xsurfvec}_\normalbase)$ and obtain
\begin{equation}
    \der{\hgridfaces^+}{t}{} \had (\velgrid^+ - \velgrid^-) \approx - \regds_{\faces, 1n} (\dot{\xsurfvec}_\normalbase \had \jump{\velpointnormal}),
    \label{eq:fdvdt_TS}
\end{equation}
where $\dot{\xsurfvec}_{\normalbase,l} \coloneqq \normalvec[l] \cdot \dot{\xsurfvec}_{l} $ with $\dot{\xsurfvec}$ the velocity of the $l$-th IB point.

Using the approximations~\eqref{eq:fdv_TS}--\eqref{eq:fdvdt_TS}, we can formulate the final versions of the spatially-discretized momentum and continuity equations for the composite fields:
\begin{multline}
    \der{\mask{\velgrid}}{t}{} 
    + \convecgrid ( \mask{\velgrid} )
    + \regds_{\faces, 1 \normalbase} \big( \dot{\xsurfvec}_\normalbase  \had \jump{\velpointnormal} \big)
    = \\ 
    - \gradgrid \mask{\pgrid}
    + \frac{1}{\Re} \lapgrid_{\faces} \mask{\velgrid}
    + \regds_{\faces} \big( \normalvec \had \spacetransform{\spoints}{\vpoints} \jump{\ppoint} \big)
    - \frac{1}{\Re} \Big(
        \spacetransform{\tensorspace}{\faces}\regds_{(\spacetransform{}{}\faces)^\top} \outerproductoperator_{\normalvec \had \normalvec} \big( \jump{\velpointnormal} \big)
        + \divgrid_{\tensorspace} \regds_{(\spacetransform{}{}\faces)^\top, 1 \normalbase} \outerproductoperator_{\normalvec} \big( \jump{\velpointnormal} \big)
    \Big) + \mask{\bcgridvec}_1 ,
    \label{eq:momentum_composite_TS}
\end{multline}
\begin{equation}
    \divgrid \mask{\velgrid} = \spacetransform{\faces}{\centers}\regds_{\faces, 1 \normalbase} \big( \normalvec \had \jump{\velpointnormal} \big) + \mask{\bcgrid}_2 .
    \label{eq:continuity_composite_TS}
\end{equation}

\subsection{Interpolation of the velocity and pressure to the immersed boundary}
\label{subsec:ns_interpolations}

We also use the Taylor series construction~\eqref{eq:v_TS} for $\velgrid$ near the IB points to improve the accuracy of the interpolation operator~\eqref{eq:disc_interp} for vector-valued data. In our case that $\velpointsurfplusminus = \velpointsurf$, we follow the approach equivalent to~\eqref{eq:Eu_TS}, but now for vector-valued data, to obtain the following second-order accurate interpolation for the composite solution (provided we use a DDF that satisfies the first-order discrete moment condition):
\begin{equation}
    \interpbody_{\faces} \mask{\velgrid} - \jump{\velpointnormal} \had \interpbody_{\faces, 1n} \hgridfaces^+ \approx \velpointsurf,
    \label{eq:disc_vel_interp_TS}
\end{equation}
where we introduced the normal-distance-weighted interpolation operation $\interpbody_{\faces, 1\normalbase} : \R^{\vpoints} \mapsto \R^{\faces}$, defined in~\eqref{eq:EF1n}.

The discrete momentum equation~\eqref{eq:momentum_composite_TS} contains the jump of the pressure $\jump{\ppoint}$ across the interface in addition to the composite pressure field $\mask{\pgrid}$. We can find a relation between both by using the normal-distance-weighted interpolation operation~\eqref{eq:EC1n}, the relation \eqref{eq:EC1n1} for DDFs that satisfy the first-order discrete moment condition, and the zeroth-order term of the Taylor series expansion~\eqref{eq:s_TS} for $\pgrid^\pm$:
\begin{align}
    \indicesbig{\interpbody_{\centers,1n} \mask{\pgrid}}{l}
    &=
    \dx\dy\dz \sum_{i,j,k} \indicesbig{\dddf_{\centers,l}}{(i,j,k)} \Big(\normalvec[l] \cdot \big( \xgridvec[\centers,(i,j,k)] - \xsurfvec[l] \big)\Big) \Big( \indicesbig{\hgridcenters^+ \had \pgrid^+}{(i,j,k)} + \indicesbig{\hgridcenters^- \had \pgrid^-}{(i,j,k)} \Big) \nonumber    \\
    &= \dx\dy\dz \sum_{i,j,k} \indicesbig{\dddf_{\centers,l}}{(i,j,k)} \Big(\normalvec[l] \cdot \big( \xgridvec[\centers,(i,j,k)] - \xsurfvec[l] \big)\Big) \Big( \hgridcenters[,(i,j,k)]^{+} \ppointsurfplus[,l]+ \indicesbig{1 - \hgridcenters^{+}}{(i,j,k)}  \ppointsurfminus[,l]
    \Big) \nonumber \\
    &= \jump[\surf,l]{\ppoint} \dx\dy\dz \sum_{i,j,k} \indicesbig{\dddf_{\centers,l}}{(i,j,k)} \Big(\normalvec[l] \cdot \big( \xgridvec[\centers,(i,j,k)] - \xsurfvec[l] \big)\Big)  \hgridcenters[,(i,j,k)]^{+} \nonumber \\
    &= \indicesBig{\jump{\ppoint} \had \interpbody_{\centers,1n} \hgridcenters^+}{l} .
    \label{eq:disc_p_interp_TS}
\end{align}

We also define a modified normal-distance-weighted interpolation operation $\tilde{\interpbody}_{\centers,1n}$ that removes the average value from the final result:
\begin{equation}
    \tilde{\interpbody}_{\centers,1n} \mask{\pgrid} = \interpbody_{\centers,1n} \mask{\pgrid} - \sum_l \indices{\interpbody_{\centers,1n} \mask{\pgrid}}{l} / \numpts.
    \label{eq:disc_p_interp_TS_zero_average}
\end{equation}

\subsection{Projection-based approach to solving the constrained system}
\label{subsec:ns_projection_method}

For simplicity, we discretize the discrete momentum equation \eqref{eq:momentum_composite_TS} in time using a fully-explicit time-stepping scheme and we define the right-hand-side vector $\bm{\grid{r}} \in \R^{\faces}$,
\begin{equation}
    \bm{\grid{r}} =  \mask{\velgrid} + \timestepcoef\dt \left( -\convecgrid(\mask{\velgrid}) + \frac{1}{\Re} \lapgrid_{\faces} \mask{\velgrid} + \mask{\bcgridvec}_1 \right),
\end{equation}
where $\dt \in \R$ is the time-step size and $\timestepcoef \in \R$ is a coefficient that depends on the chosen time-stepping scheme.
Then, we can formulate the discrete Navier--Stokes system that we have to solve to advance the solution by one stage of the timestepping scheme. This system consists of the temporally-discretized momentum equation \eqref{eq:momentum_composite_TS}, constrained by the continuity equation \eqref{eq:continuity_composite_TS}, the IB boundary condition \eqref{eq:disc_vel_interp_TS}, and the relation between the pressure field and pressure jump \eqref{eq:disc_p_interp_TS}, which we can formulate as a saddle point system,
\begin{equation}
    \begin{bmatrix}
    \id
    &
    \dt \gradgrid
    &
    \dt \regds_{\faces, 1 \normalbase} \diag ( \dot{\xsurfvec}_\normalbase )
    + \dfrac{\dt}{\Re} \big( \spacetransform{\tensorspace}{\faces}\regds_{(\spacetransform{}{}\faces)^\top} \outerproductoperator_{\normalvec \had \normalvec} 
    + \divgrid_{\tensorspace} \regds_{(\spacetransform{}{}\faces)^\top, 1 \normalbase} \outerproductoperator_{\normalvec} \big) & -\dt \regds_{\faces} \diag (\normalvec) \spacetransform{\spoints}{\vpoints}
    \\
    \divgrid & 0 & -\spacetransform{\faces}{\centers}\regds_{\faces,1n} \diag (\normalvec) & 0
    \\
    \interpbody_{\faces} & 0 & -\diag \big( \interpbody_{\faces,1n} \hgridfaces^+ \big) & 0
    \\
    0 & \tilde{ \interpbody }_{\centers,1n} & 0 & -\diag \big( \interpbody_{\centers,1n} \hgridcenters^+ \big)
    \end{bmatrix}
    \begin{pmatrix}
    \mask{\velgrid}
    \\
    \mask{\pgrid}
    \\
    \jump{\velpointnormal}
    \\
    \jump{\ppoint}
    \end{pmatrix}
    =
    \begin{pmatrix}
    \bm{\grid{r}}
    \\
    \mask{\bcgrid}_2
    \\
    \velpointsurf
    \\
    \zerosspoint
    \end{pmatrix},
    \label{eq:ns_new_saddle_system}
\end{equation}
where $\bm{\grid{r}}$ is constructed using the solution values from the previous timestepping stage. Following a projection-based approach, the pressure field is computed to ensure that the velocity field at the new timestep satisfies the continuity equation. The variable $\jump{\ppoint}$, representing the jump in the pressure field, is accordingly also treated as unknown. Similar to the pressure field, $\jump{\velpointnormal}$ is computed to ensure that the velocity field satisfies the no-slip boundary condition at the new timestep and is treated as unknown, despite being derived from the explicitly handled viscous term.

For comparison, we also provide the formulation of a prototypical continuous-forcing IB method adapted for our timestepping scheme and incorporated into our notation as
\begin{equation}
    \begin{bmatrix}
    \id & \dt \gradgrid & \dt \regds_{\faces} 
    \\
    \divgrid & 0  & 0
    \\
    \interpbody_{\faces} & 0 & 0
    \end{bmatrix}
    \begin{pmatrix}
    \mask{\velgrid}
    \\
    \mask{\pgrid}
    \\
    \forcepointvec
    \end{pmatrix}
    =
    \begin{pmatrix}
    \bm{\grid{r}}
    \\
    \mask{\bcgrid}_2
    \\
    \velpointsurf
    \end{pmatrix},
    \label{eq:ns_old_saddle_system}
\end{equation}
where $\forcepointvec$ is the IB forcing strength. A variety of solutions to the prototypical continuous-forcing IB method exist, for which the most direct analog is the immersed boundary projection method of \citet{taira_immersed_2007}. However, the system as written in \eqref{eq:ns_old_saddle_system} applies to a number of continuous-forcing IB methods. It can be adapted to many others by allowing the IB force term to be computed explicitly in time (from motion of the immersed surface at prior time instances), so that it is applied directly on the right hand side as opposed to being incorporated as an unknown in the governing system of equations. 
One key difference between the proposed system and the prototypical system is that the former yields the diagonal matrices in the lower-right blocks. Similar to the Poisson problem, the presence of these matrices has the effect of removing the ill-conditioning associated with the prototypical treatment.

To solve the proposed system~\eqref{eq:ns_new_saddle_system}, we first define the block matrices,
\begin{align}
    B^\top_{1,1} &=
    \begin{bmatrix}
    \dt \regds_{\faces, 1 \normalbase} \diag (\dot{\xsurfvec}_\normalbase )
    + \dfrac{\dt}{\Re} \big( \spacetransform{\tensorspace}{\faces}\regds_{(\spacetransform{}{}\faces)^\top} \outerproductoperator_{\normalvec \had \normalvec} 
    + \divgrid_{\tensorspace} \regds_{(\spacetransform{}{}\faces)^\top, 1 \normalbase} \outerproductoperator_{\normalvec} \big) & -\dt \regds_{\faces} \diag (\normalvec) \spacetransform{\spoints}{\vpoints}
    \end{bmatrix},
    \\
    B^\top_{1,2} &=
    \begin{bmatrix}
    -\spacetransform{\faces}{\centers}\regds_{\faces,1n} \diag (\normalvec) & 0
    \end{bmatrix} ,
\end{align}
and then formulate the solution procedure based on Schur complement reduction method,
\begin{align}
    \dt \lapgrid \mask{\pgrid}^{*} &= \divgrid \bm{\grid{r}} - \mask{\bcgrid}_2
    \label{eq:ns_intermediate_pressure_equation}
    \\
    \mask{\velgrid}^{*} &= \bm{\grid{r}} - \dt \gradgrid \mask{\pgrid}^{*} 
    \label{eq:ns_intermediate_velocity_equation}
    \\
    \schur
    \begin{pmatrix} 
    \jump{\velpointnormal} \\ 
    \jump{\ppoint} 
    \end{pmatrix} 
    &= 
    \begin{pmatrix} 
    \velpointsurf \\ 
    \zerosspoint 
    \end{pmatrix} 
    - 
    \begin{bmatrix} 
    \interpbody_{\faces} & 0 \\ 
    0 & \tilde{ \interpbody }_{\centers,1n} 
    \end{bmatrix} 
    \begin{pmatrix} 
    \mask{\velgrid}^* \\ 
    \mask{\pgrid}^* 
    \end{pmatrix} 
    \label{eq:ns_new_force_equation}
    \\
    \mask{\pgrid} &= \mask{\pgrid}^{*} - \frac{1}{\dt}\invlapgrid \big( \divgrid B^\top_{1,1} - B^\top_{1,2} \big) \begin{pmatrix} \jump{\velpointnormal} \\ \jump{\ppoint} \end{pmatrix} 
    \\
    \mask{\velgrid} &= \bm{\grid{r}} - B^\top_{1,1} \begin{pmatrix} \jump{\velpointnormal} \\ \jump{\ppoint} \end{pmatrix} - \dt \gradgrid \mask{\pgrid}
    \label{eq:ns_final_velocity_equation}
\end{align}
where $\schur$ is the Schur complement of the system:
\begin{equation}
    \schur =
    \begin{bmatrix}
        -\interpbody_{\faces} B^\top_{1,1}  + \interpbody_{\faces} \gradgrid \invlapgrid \big( \divgrid B^\top_{1,1} - B^\top_{1,2} \big)
        \\
        -\tilde{ \interpbody }_{\centers,1n} \invlapgrid \big( \divgrid B^\top_{1,1} - B^\top_{1,2} \big)
    \end{bmatrix}
    +
    \begin{bmatrix}
        -\diag \big( \interpbody_{\faces,1n} \hgridfaces^+ \big) & 0 \\
        0 & - \diag \big( \interpbody_{\centers,1n} \hgridcenters^+ \big)
    \end{bmatrix}.
    \label{eq:ns_new_schur}
\end{equation}

The solution procedure \eqref{eq:ns_intermediate_pressure_equation}--\eqref{eq:ns_final_velocity_equation} involves several steps that require the solution of Poisson equations. In our implementation, these are solved using a fast Fourier transform (FFT)–based method with the zero Fourier mode set to zero. The procedure also requires the solution of the linear system \eqref{eq:ns_new_force_equation}. For moving bodies or problems with a large number of degrees of freedom, this system is most efficiently solved using an iterative method, such as a biconjugate gradient algorithm. Alternatively, for smaller problems with fewer points used to represent the immersed surface, the linear system involving $\schur$ may be solved directly using, for example, an LU factorization. 

In our system~\eqref{eq:ns_new_saddle_system}, we make use of the pressure interpolation~\eqref{eq:disc_p_interp_TS_zero_average} that produces a result with a zero spatial mean. This choice allows the iterative solver to converge to a unique solution, as the solution for $\jump{\ppoint}$ is only defined up to an additive constant. When directly inverting $\schur$, one also has to replace one row per closed IB curve in the final block row of~\eqref{eq:ns_new_saddle_system}---specifically, Eq.~\eqref{eq:disc_p_interp_TS} for one point on each closed IB curve---by an equation that sets the average of $\jump{\ppoint}$ to zero (per closed IB curve). For example, in case there is a single closed curve, we can use
\begin{equation}
    \sum_l  \jump[\surf,l]{\ppoint} / \numpts = 0.
\end{equation}

\subsection{Results for a 2D circular Couette flow problem}

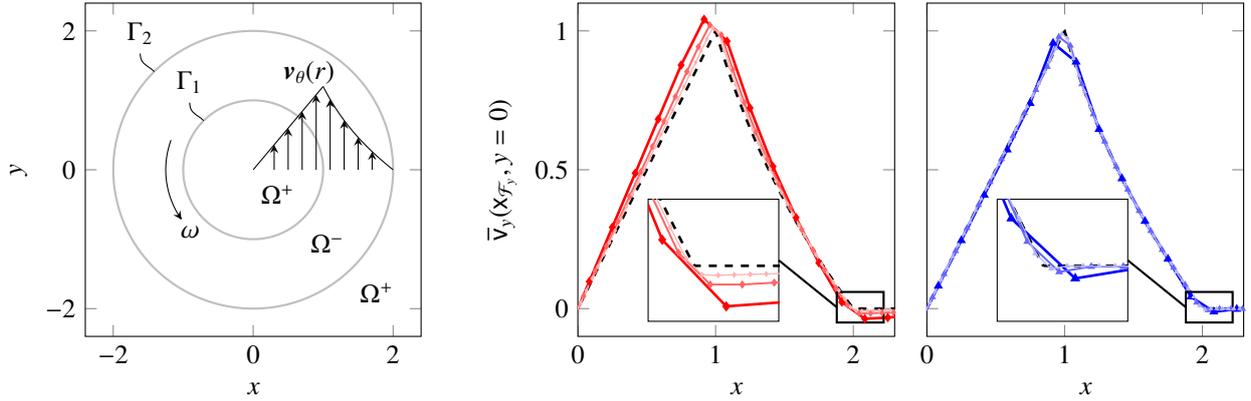
\begin{figure}
    \centering
  \tikzsetnextfilename{figures/2D_double_cylinder_diagram}%
  \input{figures/2D_double_cylinder_diagram.tikz}%

    \hfill
  \tikzsetnextfilename{figures/2D_double_cylinder_solution_original}%
  \input{figures/2D_double_cylinder_solution_original.tikz}%

  \tikzsetnextfilename{figures/2D_double_cylinder_solution_new}%
  \input{figures/2D_double_cylinder_solution_new.tikz}%

    \caption{(Left) Diagram of the circular Couette flow problem. (Center and right) Comparison of the vertical velocity solution on the line along $y=0$ near the inner cylinder for $\ds / \dx \approx 1$ obtained using the prototypical continuous-forcing IB method (red, diamonds) and the proposed method (blue, triangles) for $\ds / \dx \approx 1$ and three different grid cell sizes $\dx$: 0.167 (darkest, thick), 0.0833, and 0.0417 (lightest, thin). The exact solution is also shown (black, dashed).}
    \label{fig:double_cylinder_solution}
\end{figure}

To assess the accuracy and conditioning of the proposed method, we consider the two-dimensional Navier--Stokes problem of the circular Couette flow between two concentric cylinders with radii $R_1$ and $R_2$, shown in the left panel of Fig.~\ref{fig:double_cylinder_solution}. The inner cylinder $\surf_{1}$ rotates with a constant angular velocity $\omega$, while the outer cylinder $\surf_{2}$ remains stationary. 
In this study, we set the outer radius to be twice the inner radius, corresponding to a radius ratio $\kappa \coloneqq R_1 / R_2 = 0.5$.
The analytical solution to this problem, non-dimensionalized using the length scale $R_1$, velocity scale $\omega R_1$, and radius ratio $\kappa$ is
\begin{equation}
    \vel_\theta(r) = 
    \begin{cases}
        \displaystyle r, & \displaystyle r \le 1,\\
        \displaystyle \frac{\kappa^2}{1 - \kappa^2} \left( \frac{1}{\kappa^2 r} - r \right), & \displaystyle 1 < r \le \frac{1}{\kappa},\\
        0, & \displaystyle r > \frac{1}{\kappa},
    \end{cases}
\end{equation}
where $\vel_\theta$ is the azimuthal velocity. We treat the region between the inner and outer cylinder as the interior region $\domain^-$ and the region inside the inner cylinder and outside the outer cylinder as the exterior region $\domain^+$. In this case, the normals at the inner and outer cylinders point inward and outward, respectively. The analytical normal-derivative jumps over the inner and outer cylinders of the azimuthal velocity are then
\begin{equation}
    \jump[\surf_1]{\vel^\normalbase_\theta} = -\frac{2}{1 - \kappa^2}, \quad
    \jump[\surf_2]{\vel^\normalbase_\theta} = \frac{2 \kappa^2}{1 - \kappa^2}.
\end{equation}

We again compute and evaluate the error of the numerical solution and interface forcing term, for various discretizations of $\domain = \left\{ (x,y) \in \R^2 \;\middle|\; -2.67 \le x \le 2.67,\; -2.67 \le y \le 2.67 \right\}$. We use the lattice Green's function to treat the far-field boundary conditions, and we keep the grid spacing uniform throughout the domain and identical in the $x$- and $y$-directions in all cases. 
To assess problem conditioning and accuracy in the IB forcing, we compute results a range of values for the body-to-grid spacing ratio $\ds / \dx$. 
The solutions are computed by time-stepping the system from zero initial conditions until steady-state using the Schur complement reduction method applied to both the proposed IB system~\eqref{eq:ns_new_saddle_system} and the prototypical continuous-forcing IB system~\eqref{eq:ns_old_saddle_system}. For the former, this resulted in the solution procedure \eqref{eq:ns_intermediate_pressure_equation}--\eqref{eq:ns_final_velocity_equation}. For the latter, we adopt the projection solution procedure of \citet{taira_immersed_2007}. That is,~\eqref{eq:ns_old_saddle_system} is solved using the same block-LU factorization process we apply to~\eqref{eq:ns_new_force_equation}. The result is
identical to that presented in~\eqref{eq:ns_intermediate_pressure_equation}--\eqref{eq:ns_final_velocity_equation}, but with the last three steps replaced by
\begin{align}
    \tilde{\schur} \forcepointvec &= \velpointsurf - \interpbody_{\faces} \mask{\velgrid}^*,
    \label{eq:ns_old_force_equation} 
    \\
    \mask{\pgrid} &= \mask{\pgrid}^{*} - \invlapgrid \regds_{\faces} \forcepointvec 
    \\
    \mask{\velgrid} &= \bm{\grid{r}} - \dt \regds_{\faces} \forcepointvec - \dt \gradgrid \mask{\pgrid},
    \label{eq:ns_final_old_velocity_equation}
\end{align}
where $\tilde{\schur}$ is the Schur complement of the prototypical formulation~\eqref{eq:ns_old_saddle_system}; i.e., $\tilde{\schur} = - \dt \interpbody_{\faces} \regds_{\faces} + \dt \interpbody_{\faces} \gradgrid \invlapgrid \divgrid \regds_{\faces}$. Note that the proposed method has two forcing variables ($\jump{\velpointnormal}$ and $\jump{\ppoint}$) while the prototypical formulation has one ($\forcepointvec$). However, because the pressure jump in this problem is zero, we directly compare $\jump{\velpointnormal}$ and $\forcepointvec$. We utilize the projection-based continuous-IB method as the point of comparison because it is the most direct analog to the proposed methodology. Both approaches compute the IB strength via projection (enacted by analytical block-LU decomposition) to exactly enforce the interface constraint (to within machine precision). This enables a direct comparison of problem conditioning and accuracy of both the velocity and IB forcing variables under the same projection formalism. That said, we emphasize that the proposed approach---modifying the constraint equation via Taylor series to remove first-order errors in interpolating the velocity to the interface---can be used to bring higher than first-order accuracy across a broad range of continuous-forcing IB methods. Moreover, we will demonstrate that the projection-based solution procedure enables a well-conditioned, implicit treatment that provides accurate surface stresses. This functionality may be especially useful for moving and deforming body problems where body points move close to one another and the surface motion is coupled nonlinearly to the IB forcing at the interface.

\begin{figure}
    \centering  
  \tikzsetnextfilename{figures/2D_double_cylinder_error_Linf_norm}%
  \input{figures/2D_double_cylinder_error_Linf_norm.tikz}%

    \hspace{2cm}
  \tikzsetnextfilename{figures/2D_double_cylinder_error_L2_norm}%
  \input{figures/2D_double_cylinder_error_L2_norm.tikz}%

    \caption{(Infinity norm (left) and 2-norm (right) of the velocity solution error with respect to their analytical values for the circular Couette flow example for the prototypical continuous-forcing IB method (red, diamonds) and the proposed method (blue, triangles). The velocity error is computed over all cell faces (solid lines) and over all cell faces excluding those within the support of $\spacetransform{\centers}{\faces} \spacetransform{\faces}{\centers} \regds_{\faces}$ (dashed lines). 
    First- and second-order error scales are also shown (gray, solid).
    }
    \label{fig:double_cylinder_error}
\end{figure}
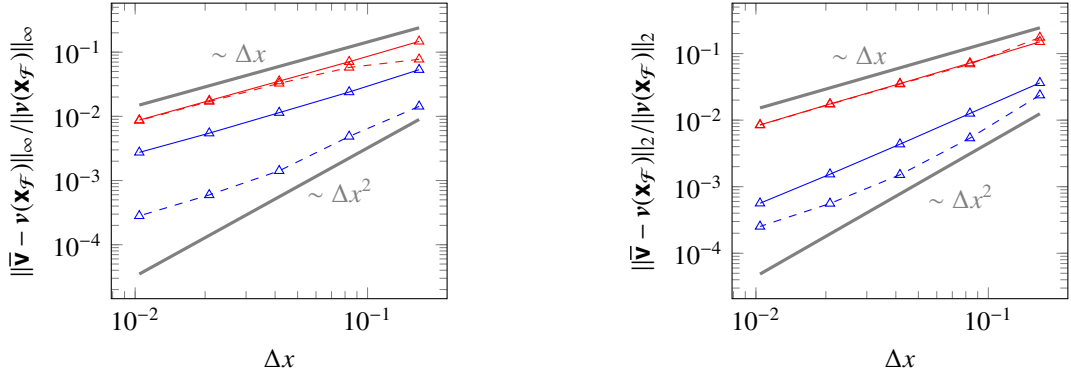

We show the solutions obtained with the prototypical method and with our proposed method on the middle and right panels of Fig.~\ref{fig:double_cylinder_solution}, respectively. Figure~\ref{fig:double_cylinder_error} shows the associated relative error, with respect to the true solution, of the solution $\mask{\velgrid}$ versus the grid cell size. Note that the solution obtained with the proposed method (right panel of Fig.~\ref{fig:double_cylinder_solution}) is visually close to the exact solution except near the interfaces, where it converges to $\hgridfaces^+ \had \vel^+(\xgrid[\faces]) + \hgridfaces^- \had \vel^-(\xgrid[\faces])$, which differs from the exact composite solution $\vel$, which is defined with sharp indicator functions in Eq.~\eqref{eq:composite_analytic_vel}.

We make similar observations to those for the Poisson problem test case in Section~\ref{subsec:poisson_2d_example}. Firstly, Fig.~\ref{fig:double_cylinder_error} provides the error convergence of the computed solutions. The solution variable $\mask{\velgrid}$ provides a convergent approximation to the true solution for both IB methods, with the error converging at first order with respect to both the infinity and $2$-norms for the prototypical continuous-forcing formulation. The solution error for the proposed approach converges at first order in the global infinity norm and between first and second order in the two norm. However, the error convergence when omitting the spatial points within the support of the DDF is close to second order for coarser grids, and slightly better than first order for the finer grids. This result is distinct from the Poisson case, where the error converged at second-order accuracy when excluding the support of the DDF.
We believe that this reduction in accuracy is due to the neglect of the second term on the right-hand side of \eqref{eq:true_heavi}, which has a significant impact on the continuity forcing term but can be neglected in the momentum equation and in the Poisson problem for the grid spacings considered here. A detailed analysis of this error source and strategies for its mitigation are the subject of ongoing work.

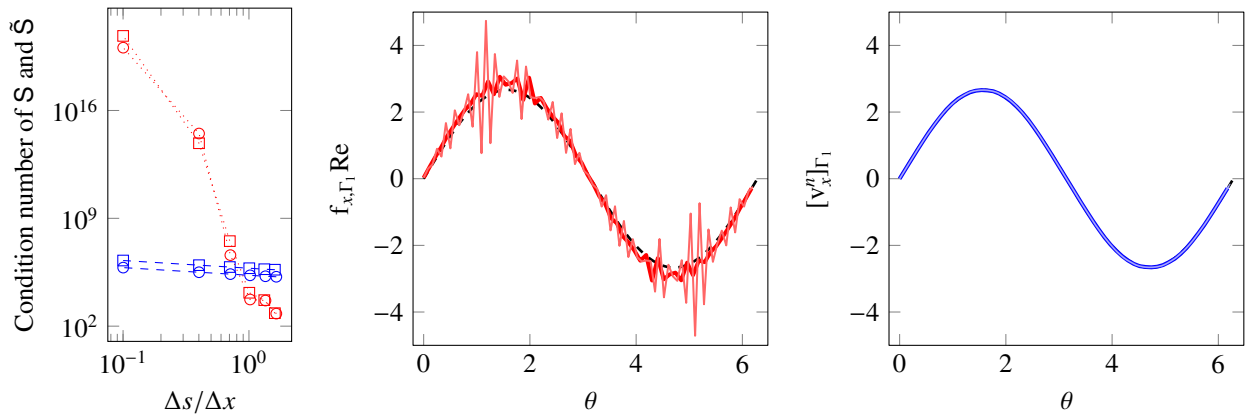
\begin{figure}
    \centering
  \tikzsetnextfilename{figures/2D_double_cylinder_condition_number}%
  \input{figures/2D_double_cylinder_condition_number.tikz}%

    \hfill
  \tikzsetnextfilename{figures/2D_double_cylinder_forcing_original}%
  \input{figures/2D_double_cylinder_forcing_original.tikz}%

    \hfill
  \tikzsetnextfilename{figures/2D_double_cylinder_forcing_new}%
  \input{figures/2D_double_cylinder_forcing_new.tikz}%

    \caption{(Left) Condition number of $\schur$ and $\tilde{\schur}$ versus the surface-to-grid spacing ratio for the proposed method (blue, dashed) and the prototypical continuous-forcing IB method (red, dotted) using two different grid cell sizes $\dx$: 0.167 (circles) and 0.0833 (squares). (Center and right) $x$-component of the forcing solution on the interior cylinder obtained with the prototypical formulation (red, dotted) and the proposed method (blue, dashed) for $\dx = 0.0833$ and three approximate values of $\ds / \dx$: 1.3 (darkest, thick), 1.0, and 0.7 (lightest, thin). The exact solution is also shown (black, dashed). Note that the forcing strength obtained in the prototypical formulation is scaled by the Reynolds number $\Re = \omega R_1^2 / \nu$ to facilitate direct comparison with the normal-derivative jump of the velocity obtained with the proposed method, which is possible because the pressure jumps over the interfaces are zero. The solution from the prototypical formulation for $\ds/\dx \approx 0.7$ cannot be meaningfully represented on the chosen $y$-axis scale and is therefore not included.}
    \label{fig:double_cylinder_forcing}
\end{figure}

Fig.~\ref{fig:double_cylinder_forcing} shows the effect of $\Delta s/\Delta x$ on the problem conditioning and IB force accuracy for the prototypical continuous-forcing and proposed methods. Decreasing $\Delta s/\Delta x$ leads to progressively poor approximations to the IB forcing strength when using the prototypical interface constraint. Indeed, for $\ds/\dx \approx 0.7$ the IB forcing could not be computed with sufficient accuracy to be conveyed on the vertical axis used for the IB force plot. 
As in the Poisson problem, this behavior is again a direct consequence of the ill-posed nature of the underlying integral equation that is discretized by Eq.~\eqref{eq:ns_old_force_equation} with $\tilde{\schur}$, the Schur complement for the prototypical formulation. As a result, for decreasing $\ds/\dx$ values, the condition number of $\tilde{\schur}$ (left panel of Fig.~\ref{fig:double_cylinder_forcing}) rapidly increases and the IB forcing along the surface (middle panel of Fig.~\ref{fig:double_cylinder_forcing}) shows increasingly large oscillations. By contrast, the proposed method, which modifies the IB constraint using Taylor series, has a well-conditioned Schur complement $\schur$ (left panel Fig.~\ref{fig:double_cylinder_forcing}) and provides a smooth forcing solution (right panel of Fig.~\ref{fig:double_cylinder_forcing}) that is visually clearly accurate and essentially independent of the ratio $\Delta s/\Delta x$. The forcing error for the proposed method initially converges but levels off at finer $\dx$ values (not shown here). As in the Poisson problem, exploring $\ds / \dx$ ratios could allow one to explicitly explore convergence of the forcing error for the proposed method.

\section{Conclusion}

We present an immersed boundary method that formally retains a DDF in a continuous-forcing IB formulation, and therefore avoids stencil modifications. Our proposed method simultaneously extends the order of accuracy of these formulations to greater than one, and removes the well-known ill-conditioning of the linear system that arises from applying a projection-based solution procedure to the formulation. 

Our approach leverages two key ingredients. The first is to formally cast the solution as a composite solution that connects two distinct solutions across the interface. This composition is performed in both a continuous and discrete setting using the indicator function representation proposed in \citet{eldredge_method_2022}. This representation allows us to modify the IB interpolation operation in the constraint equation to account for the non-smoothness of typical IB solution fields, which otherwise limits the global order of accuracy to one. This is achieved by explicitly accounting for the jump in the solution gradient and for the smoothing of the discrete indicator function that occurs due to its connection with the smoothed DDF. This connection is explicitly retained within the governing equation for the composite solution to ensure the validity of the composite solution definition. The second key ingredient is to use Taylor series approximations to recast the governing and constraint equations in terms of only the composite solution and surface distributions, which are treated as IB forcing strengths. This enables us to formulate a system similar to that of typical continuous-forcing IB methods but with higher-order accuracy, which we can easily solve using a projection-based solution procedure.

This approach---leveraging Taylor series to modify the interpolation operation in the constraint equation---has the added benefit of regularizing the computation of the IB force term. Typical continuous-forcing IB forcing methods produce surface stresses with unphysical oscillations, which can become numerically unbounded as the ratio of spacing between the IB surface and the flow domain, $ \ds/ \dx$, becomes small. We demonstrate that a projection-based solution procedure can be used for our higher-order formulation to produce accurate IB forces that are robust to $\ds / \dx$ variations. This outcome arises because the Taylor series treatment alters the equation for the IB force term from an ill-posed first-kind integral equation to a well-posed second-kind integral equation. This regularization is a natural by-product of the higher-order treatment, and provides robustly well-posed equations for the IB force without use of heuristic parameters or post-processing procedures. Moreover, the higher-order accurate, well-conditioned approach introduces low-cost computations (involving Taylor series and operations restricted to the support of the DDF), so that the proposed formulation can be implemented at a comparable cost to its first-order continuous-forcing analog. In fact, the problem conditioning in solving for the IB forces can be a bottleneck with projection-based continuous-forcing IB methods, particularly for moving-body problems where the IB force equation is solved via linear iterative techniques. In this setting, the poor conditioning can significantly increase the number of iterations needed to compute the solution. For these problems, the proposed robustly well-conditioned procedure may be faster than its first-order accurate counterparts.

We demonstrate higher-order accuracy and smooth, accurate IB forces with robust problem conditioning across a range of $\ds/ \dx$ values, many of which are inaccessible to some continuous-IB forcing methods. This demonstration is performed for 1D and 2D Poisson problems, and for the problem of concentric rotating cylinders generating a laminar, incompressible, viscous flow governed by the Navier--Stokes equations. In the Poisson problems, the solution converges at second-order accuracy omitting error contributions within the support of the DDF. When error within the support of the DDF is included, the convergence rate is between first and second order. This scaling is the same as if the exact IB force strength were provided analytically~\citep{tornberg_numerical_2004}, and is taken here as the upper bound for achievable accuracy without changing the problem formulation to enforce solution smoothness across the interface~\citep{griffith_order_2005,stein_immersed_2017}. In the Navier--Stokes setting, a similar scaling is obtained except that for fine values of $\dx$ the error outside of the support of the DDF drops below second-order accuracy. We hypothesize this result is because of inherent properties of the discrete DDF distributions considered, which have a component in the null space of the discrete divergence for curved surfaces. For both the Poisson and Navier--Stokes problems, the equation to compute the IB force is well conditioned and leads to accurate IB forces, including for $\ds/\dx$ ratios that yield unusable solutions using a first-order continuous-forcing projection IB approach. 

The results of our proposed method suggest that continuous-forcing IB methods are not confined to first-order accuracy and indicate the potential to achieve exact second-order or higher accuracy for a wide variety of PDE problems, provided that the gradient of the discrete indicator functions can be approximated more accurately by the regularized delta functions. Future work could enable higher-order formulations, including a treatment of the component of the DDF that lies in the null space of the divergence operator, and could more systematically explore the convergence properties of IB forcing. These efforts could help further clarify the scope and limitations of the proposed approach.

\section*{Acknowledgments}

This work was supported by the United States Air Force Office of Scientific Research (AFOSR MURI FA9550-23-1-0299).

\appendix

\section{Discrete operations and tools}
\label{app:discrete_operations_and_tools}

In this section, we formulate the discrete operations, identities, and shorthands used in this work. We will make use of a discrete scalar field $\scalargrid$, vector field $\vectorgrid$ on cell faces, vector field $\edgevectorgrid$ on cell edges , and tensor field $\tensorgrid$, where
\begin{gather}
    \scalargrid \in \R^{\centers}, 
    \qquad \vectorgrid = \begin{pmatrix}
       \vectorgridcomponent{x} \\
       \vectorgridcomponent{y} \\
       \vectorgridcomponent{z}
    \end{pmatrix} \in \R^{\faces}, 
    \qquad \edgevectorgrid = \begin{pmatrix}
       \edgevectorgridcomponent{x} \\
       \edgevectorgridcomponent{y} \\
       \edgevectorgridcomponent{z}
    \end{pmatrix} \in \R^{\edges},
    \qquad \tensorgrid = \begin{pmatrix}
        \tensorgridcomponent{11} & \tensorgridcomponent{12} & \tensorgridcomponent{13} \\
        \tensorgridcomponent{21} & \tensorgridcomponent{22} & \tensorgridcomponent{23} \\
        \tensorgridcomponent{31} & \tensorgridcomponent{32} & \tensorgridcomponent{33} \\
    \end{pmatrix} \in \R^{\tensorspace},
\end{gather}
and their composite counterparts (except for the vector on cell edges), $\mask{\scalargrid} = \hgridcenters^+ \had \scalargrid^+ + \hgridcenters^- \had \scalargrid^-$, $\mask{\vectorgrid} = \hgridfaces^+ \had \vectorgrid^+ + \hgridfaces^- \had \vectorgrid^-$, and $\mask{\tensorgrid} = \hgridtensorspace^+ \had \tensorgrid^+ + \hgridtensorspace^- \had \tensorgrid^-$. 

In the notation above, $\R^{\centers}$ represents the space of discrete scalar data, $\R^{\faces} \coloneqq \R^{\facescomponent_x} \times \R^{\facescomponent_y} \times \R^{\facescomponent_z}$ and $\R^{\edges} \coloneqq \R^{\edgescomponent_x} \times \R^{\edgescomponent_y} \times \R^{\edgescomponent_z}$ are the spaces of discrete vector values data on cell faces and edges, respectively, and $\R^{\tensorspace} \coloneqq \Pi_{m,n=1}^{3} \R^{\gspace{G}_{m,n}}$ is the space of discrete, second-order tensor fields whose components are defined on the following combination of cell locations:
\begin{equation}
    \left(\gspace{G}_{m,n} \right) 
    = \begin{pmatrix}
    \centers & \edgescomponent_z & \edgescomponent_y \\
    \edgescomponent_z & \centers & \edgescomponent_x \\
    \edgescomponent_y & \edgescomponent_x & \centers
    \end{pmatrix}.
\end{equation}

We will also use discrete scalar data $\scalarpoint$, vector data $\vectorpoint$, and second-order tensor data $\tensorpoint$ on the IB points, where
\begin{gather}
    \scalarpoint \in \R^{\spoints}, 
    \qquad \vectorpoint = \begin{pmatrix}
       \vectorpointcomponent[x] \\
       \vectorpointcomponent[y] \\
       \vectorpointcomponent[z]
    \end{pmatrix} \in \R^{\vpoints},
    \qquad \tensorgrid = \begin{pmatrix}
        \tensorpointcomponent[11] & \tensorpointcomponent[12] & \tensorpointcomponent[13] \\
        \tensorpointcomponent[21] & \tensorpointcomponent[22] & \tensorpointcomponent[23] \\
        \tensorpointcomponent[31] & \tensorpointcomponent[32] & \tensorpointcomponent[33] \\
    \end{pmatrix} \in \R^{\tpoints},
\end{gather}
Here, $\R^{\spoints}$, $\R^{\vpoints} \coloneqq \big(\R^{\spoints}\big)^3$, and $\R^{\tpoints} \coloneqq \big(\R^{\spoints}\big)^9$ denote the spaces of scalar-valued, vector-valued, and second-order tensor-valued data on the IB points, respectively.

Throughout, we will use second-order central differences and linear interpolations, and the following grid spacing vector $\dxgridvec \in \R^{\faces}$ and tensor $\dxgridtensor \in \R^{\tensorspace}$:
\begin{equation}
    \dxgridvec
    = 
    \begin{pmatrix}
        \dx \\ \dy \\ \dz
    \end{pmatrix},
    \qquad
    \dxgridtensor
    = 
    \begin{pmatrix}
        \dx & \dy & \dz \\
        \dx & \dy & \dz \\
        \dx & \dy & \dz
    \end{pmatrix}
\end{equation}

\subsection{Space transformations}
\label{app:grid_transformations}
If we use the notations $\spacetransform{x}{}$, $\spacetransform{y}{}$, $\spacetransform{z}{}$ to denote the interpolations of a field in the three Cartesian directions by taking the two-point average and placing the result midway, we can formulate the following transformations between the different grid and surface spaces defined in this work:
\begin{itemize}
\item Interpolation and expansion of a scalar on cell centers to a vector on cell faces
\begin{equation}
    \spacetransform{\centers}{\faces} \scalargrid =
    \begin{pmatrix}
        \spacetransform{x}{} \scalargrid \\
        \spacetransform{y}{} \scalargrid \\
        \spacetransform{z}{} \scalargrid
    \end{pmatrix} \in \R^{\faces}
\end{equation}
\item Interpolation and contraction of a vector on cell faces to a scalar on cell centers
\begin{equation}
    \spacetransform{\faces}{\centers} \vectorgrid = \spacetransform{x}{} \vectorgridcomponent{x} + \spacetransform{y}{} \vectorgridcomponent{y} + \spacetransform{z}{} \vectorgridcomponent{z} \in \R^{\centers}
\end{equation}
\item Interpolation and expansion of a vector on cell faces to a tensor in the cell tensor space
\begin{equation}
    \spacetransform{\faces}{\tensorspace} \vectorgrid =
    \begin{pmatrix}
        \spacetransform{x}{}\vectorpointcomponent[x] & \spacetransform{y}{}\vectorpointcomponent[x] & \spacetransform{z}{}\vectorpointcomponent[x] \\
        \spacetransform{x}{}\vectorpointcomponent[y] & \spacetransform{y}{}\vectorpointcomponent[y] & \spacetransform{z}{}\vectorpointcomponent[y] \\
        \spacetransform{x}{}\vectorpointcomponent[z] & \spacetransform{y}{}\vectorpointcomponent[z] & \spacetransform{z}{}\vectorpointcomponent[z] \\
    \end{pmatrix} \in \R^{\tensorspace}
\end{equation}
\item Interpolation and row-wise contraction of a tensor in the cell tensor space to a vector on cell faces
\begin{equation}
    \spacetransform{\tensorspace}{\faces} \tensorgrid =
    \begin{pmatrix}  
        \spacetransform{x}{} \tensorgridcomponent{11} +
        \spacetransform{y}{} \tensorgridcomponent{12} +
        \spacetransform{z}{} \tensorgridcomponent{13} 
        \\
        \spacetransform{x}{} \tensorgridcomponent{21} +
        \spacetransform{y}{} \tensorgridcomponent{22} + 
        \spacetransform{z}{} \tensorgridcomponent{23} 
        \\
        \spacetransform{x}{} \tensorgridcomponent{31} +
        \spacetransform{y}{} \tensorgridcomponent{32} + 
        \spacetransform{z}{} \tensorgridcomponent{33} 
    \end{pmatrix} \in \R^{\faces}
\end{equation}
\item Expansion of scalar-valued surface data to vector-valued surface data
\begin{equation}
    \spacetransform{\spoints}{\vpoints} \scalarpoint
    = \begin{pmatrix}
       \scalarpoint \\
       \scalarpoint \\
       \scalarpoint
    \end{pmatrix} \in \R^{\vpoints}
    \label{eq:vIs}
\end{equation}
\end{itemize}

\subsection{Differential operators}
\label{app:differential_operators}
We use the notations $\dergrid{x}, \dergrid{y}, \dergrid{z}$ to denote the finite-difference approximation to the first derivative of a field in the three Cartesian directions. In this work we use central differences and compute the first derivative by taking the difference between values on two adjacent grid points, divided by the grid spacing in that direction, and placing the result midway. We can then formulate the following discrete differential operators:
\begin{itemize}
\item Gradient of a scalar on cell centers
\begin{equation}
    \gradgrid \scalargrid 
    = 
    \begin{pmatrix}
        \dergrid{x} \scalargrid \\
        \dergrid{y} \scalargrid \\
        \dergrid{z} \scalargrid 
    \end{pmatrix} 
    \in \R^{\faces}
\end{equation}
\item Divergence of a vector on cell faces
\begin{equation}
    \divgrid \vectorgrid = \dergrid{x}\vectorgridcomponent{x} + \dergrid{y} \vectorgridcomponent{y} + \dergrid{z}\vectorgridcomponent{z} \in \R^{\centers}
\end{equation}
\item Gradient of a vector on cell faces
\begin{equation}
    \gradgrid_{\faces} \vectorgrid =
    \begin{pmatrix}
        \dergrid{x}\vectorgridcomponent{x} & \dergrid{y}\vectorgridcomponent{x} & \dergrid{z}\vectorgridcomponent{x} \\
        \dergrid{x}\vectorgridcomponent{y} & \dergrid{y}\vectorgridcomponent{y} & \dergrid{z}\vectorgridcomponent{y} \\
        \dergrid{x}\vectorgridcomponent{z} & \dergrid{y}\vectorgridcomponent{z} & \dergrid{z}\vectorgridcomponent{z} \\
    \end{pmatrix} \in \R^{\tensorspace}
\end{equation}
\item Divergence of a tensor in the cell tensor space
\begin{equation}
    \divgrid_{\tensorspace} \tensorgrid =
    \begin{pmatrix}  
        \dergrid{x} \tensorgridcomponent{11} +
        \dergrid{y} \tensorgridcomponent{12} +
        \dergrid{z} \tensorgridcomponent{13} 
        \\
        \dergrid{x} \tensorgridcomponent{21} +
        \dergrid{y} \tensorgridcomponent{22} + 
        \dergrid{z} \tensorgridcomponent{23} 
        \\
        \dergrid{x} \tensorgridcomponent{31} +
        \dergrid{y} \tensorgridcomponent{32} + 
        \dergrid{z} \tensorgridcomponent{33} 
    \end{pmatrix} \in \R^{\faces}
\end{equation}
\item Curl of a vector on cell edges
\begin{equation}
    \curlgrid \edgevectorgrid 
    = 
    \begin{pmatrix}
        \dergrid{y} \edgevectorgridcomponent{z} - \dergrid{z} \edgevectorgridcomponent{y} \\
        \dergrid{z} \edgevectorgridcomponent{x} - \dergrid{x} \edgevectorgridcomponent{z} \\
        \dergrid{x} \edgevectorgridcomponent{y} - \dergrid{y} \edgevectorgridcomponent{x}
    \end{pmatrix} 
    \in \R^{\faces}
\end{equation}
\item Curl of a vector on cell faces
\begin{equation}
    \curlgrid^\top \vectorgrid 
    = 
    \begin{pmatrix}
        \dergrid{y} \vectorgridcomponent{z} - \dergrid{z} \vectorgridcomponent{y} \\
        \dergrid{z} \vectorgridcomponent{x} - \dergrid{x} \vectorgridcomponent{z} \\
        \dergrid{x} \vectorgridcomponent{y} - \dergrid{y} \vectorgridcomponent{x}
    \end{pmatrix} 
    \in \R^{\edges}
\end{equation}
\end{itemize}
Similarly, if we use $\dergrid{x}^2, \dergrid{y}^2, \dergrid{z}^2$, to denote the finite-difference approximation to the second derivative in the three Cartesian directions and use the standard three-point stencil, second-order central difference scheme where the result is placed at the center, we can formulate the following Laplacian operators:
\begin{itemize}
\item Laplacian of a scalar on cell centers
\begin{equation}
    \lapgrid \scalargrid = \dergrid{x}^2 \scalarpoint + \dergrid{y}^2 \scalarpoint + \dergrid{z}^2 \scalarpoint  \in \R^{\centers}
\end{equation}
\item Laplacian of a vector on cell faces
\begin{equation}
    \lapgrid_{\faces} \vectorgrid = 
    \begin{pmatrix}
        \dergrid{x}^2 \vectorgridcomponent{x} + \dergrid{y}^2 \vectorgridcomponent{x} + \dergrid{z}^2 \vectorgridcomponent{x}
        \\
        \dergrid{x}^2 \vectorgridcomponent{y} + \dergrid{y}^2 \vectorgridcomponent{y} + \dergrid{z}^2 \vectorgridcomponent{y}
        \\
        \dergrid{x}^2 \vectorgridcomponent{z} + \dergrid{y}^2 \vectorgridcomponent{z} + \dergrid{z}^2 \vectorgridcomponent{z}
    \end{pmatrix} \in \R^{\faces}
\end{equation}
\item Laplacian of a vector on cell edges
\begin{equation}
    \lapgrid_{\edges} \edgevectorgrid = 
    \begin{pmatrix}
        \dergrid{x}^2 \edgevectorgridcomponent{x} + \dergrid{y}^2 \edgevectorgridcomponent{x} + \dergrid{z}^2 \edgevectorgridcomponent{x}
        \\
        \dergrid{x}^2 \edgevectorgridcomponent{y} + \dergrid{y}^2 \edgevectorgridcomponent{y} + \dergrid{z}^2 \edgevectorgridcomponent{y}
        \\
        \dergrid{x}^2 \edgevectorgridcomponent{z} + \dergrid{y}^2 \edgevectorgridcomponent{z} + \dergrid{z}^2 \edgevectorgridcomponent{z}
    \end{pmatrix} \in \R^{\edges}
\end{equation}
\end{itemize}

\subsection{Interpolations of products}
\label{app:interpolation_of_products}

\begin{itemize}
    \item Interpolation of an element-wise scalar product in the $x$-direction (similar for $y$- and $z$-directions). Note that the position of $\indices{\spacetransform{x}{} \scalargrid}{(i,j,k)} \in \R^{\edgescomponent_x}$ lies halfway between $\scalargrid[(i,j,k)]$ and $\scalargrid[(i+1,j,k)]$.
        \begin{align}
            \indicesbig{\spacetransform{x}{} \left( \scalargrid[1] \circ \scalargrid[2]\right)}{(i,j,k)} \nonumber
            &=
            \frac{\scalargrid[1,(i+1,j,k)] \scalargrid[2,(i+1,j,k)] + \scalargrid[1,(i,j,k)] \scalargrid[2,(i,j,k)]}{2} \nonumber
            \\
            &=
            \frac{\scalargrid[1,(i+1,j,k)] + \scalargrid[1,(i,j,k)]}{2} \frac{\scalargrid[2,(i+1,j,k)] + \scalargrid[2,(i,j,k)]}{2} \nonumber
            \\
            &\quad+ \frac{\scalargrid[1,(i+1,j,k)] - \scalargrid[1,(i,j,k)]}{2} \frac{\scalargrid[2,(i+1,j,k)] - \scalargrid[2,(i,j,k)]}{2} \nonumber
            \\
            &=
            \indicesbig{\spacetransform{x}{} \scalargrid[1] \circ \spacetransform{x}{} \scalargrid[2]}{(i,j,k)} + \frac{\dx^2}{4} \indicesbig{\dergrid{x} \scalargrid[1] \circ \dergrid{x} \scalargrid[2]}{(i,j,k)}
        \end{align}
    \item Interpolation and expansion of an element-wise scalar product on cell centers to a vector on cell faces
        \begin{align}
            \spacetransform{\centers}{\faces} \left( \scalargrid[1] \circ \scalargrid[2] \right) 
            = 
            \spacetransform{\centers}{\faces} \scalargrid[1] \circ \spacetransform{\centers}{\faces} \scalargrid[2] + \frac{\dxgridvec^2}{4} \had \gradgrid \scalargrid[1] \had \gradgrid \scalargrid[2] 
        \end{align}
    \item Interpolation and expansion of an element-wise vector product on cell faces to a tensor in the cell tensor space
        \begin{align}
            \spacetransform{\faces}{\tensorspace} \left( \vectorgrid[1] \circ \vectorgrid[2] \right) 
            = 
            \spacetransform{\faces}{\tensorspace} \vectorgrid[1] \circ \spacetransform{\faces}{\tensorspace} \vectorgrid[2] + \frac{\dxgridtensor^2}{4} \had \gradgrid_{\faces} \vectorgrid[1] \had \gradgrid_{\faces} \vectorgrid[2]
            \label{eq:interpolation_of_vector_product}
        \end{align}
\end{itemize}

\subsection{Derivatives of products}

\begin{itemize}
    \item The finite-difference derivative of the element-wise product of two scalar fields in the $x$-direction (similar for $y$- and $z$-directions)
    \begin{align}
        \indicesbig{\dergrid{x} \left( \scalargrid[1] \had \scalargrid[2] \right)}{(i,j,k)} \nonumber
        &=
        \frac{\scalargrid[1,(i+1,j,k)] \scalargrid[2,(i+1,j,k)] - \scalargrid[1,(i,j,k)] \scalargrid[2,(i,j,k)]}{\dx} \nonumber
        \\
        &=
        \frac{\scalargrid[1,(i+1,j,k)] + \scalargrid[1,(i,j,k)]}{2}\frac{\scalargrid[2,(i+1,j,k)] - \scalargrid[2,(i,j,k)]}{\dx} \nonumber
        \\
        &\quad+ \frac{\scalargrid[1,(i+1,j,k)] - \scalargrid[1,(i,j,k)]}{\dx}\frac{\scalargrid[2,(i+1,j,k)] + \scalargrid[2,(i,j,k)]}{2} \nonumber
        \\
        &= \indicesbig{\spacetransform{x}{} \scalargrid[1] \had \dergrid{x} \scalargrid[2] + \spacetransform{x}{} \scalargrid[2] \had \dergrid{x} \scalargrid[1]}{(i,j,k)}
    \end{align}
    \item Using the previous relation, one can show that the gradient of the product of two scalars fields, analogous to the continuous product rule $\grad \left( s_1 s_2 \right) = s_1 \grad s_2 + s_2 \grad s_1$, is
    \begin{align}
        \gradgrid \left( \scalargrid[1] \had \scalargrid[2] \right)
        &=
        \spacetransform{\centers}{\faces} \scalargrid[1] \had \gradgrid \scalargrid[2] + \spacetransform{\centers}{\faces} \scalargrid[2] \had \gradgrid \scalargrid[1]
    \end{align}
    \item Similarly, the gradient of the element-wise product of a scalar field (interpolated to cell faces) and vector field, analogous to the continuous product rule $\grad \left( s \bm{v} \right) = s \grad \bm{v} = \bm{v} \grad s$, is
    \begin{align}
        \gradgrid_{\faces} \left( \spacetransform{\centers}{\faces} \scalargrid \had \vectorgrid \right)
        &= \spacetransform{\faces}{\tensorspace}\spacetransform{\centers}{\faces} \scalargrid \had \gradgrid_{\faces} \vectorgrid + \spacetransform{\faces}{\tensorspace}\vectorgrid \had \gradgrid_{\faces} \spacetransform{\centers}{\faces} \scalargrid 
    \end{align}
    \item The divergence of the element-wise product of a scalar field (interpolated to the cell faces) and a vector field, analogous to the continuous product rule $\div \left( s \bm{v} \right) = \bm{v} \cdot \grad s + s \div \bm{v}$, is
    \begin{align}
        \begin{split}
            \indicesBig{\divgrid \left( \spacetransform{\centers}{\faces} \scalargrid \had \vectorgrid \right)}{(i,j,k)} 
            &      = \frac{1}{\dx} \left( \frac{\scalargrid[(i+1,j,k)] + \scalargrid[(i,j,k)]}{2} \vectorgridcomponent{x,(i,j,k)}
                                        - \frac{\scalargrid[(i,j,k)] + \scalargrid[(i-1,j,k)]}{2} \vectorgridcomponent[x,(i-1,j,k)] \right) \\
            &\quad + \frac{1}{\dy} \left( \frac{\scalargrid[(i,j+1,k)] + \scalargrid[(i,j,k)]}{2} \vectorgridcomponent{y,(i,j,k)}
                                        - \frac{\scalargrid[(i,j,k)] + \scalargrid[(i,j-1,k)]}{2} \vectorgridcomponent[y,(i,j-1,k)] \right) \\
            &\quad + \frac{1}{\dz} \left( \frac{\scalargrid[(i,j,k+1)] + \scalargrid[(i,j,k)]}{2} \vectorgridcomponent{z,(i,j,k)}
                                        - \frac{\scalargrid[(i,j,k)] + \scalargrid[(i,j,k-1)]}{2} \vectorgridcomponent[z,(i,j,k-1)] \right)
        \end{split} \nonumber
        \\
        \begin{split}
            &      = \frac{1}{2} \left( \frac{\scalargrid[(i+1,j,k)] - \scalargrid[(i,j,k)]}{\dx} \vectorgridcomponent{x,(i,j,k)}
                   +                    \frac{\scalargrid[(i,j,k)] - \scalargrid[(i-1,j,k)]}{\dx} \vectorgridcomponent[x,(i-1,j,k)] \right) \\
            &\quad + \frac{1}{2} \left( \frac{\scalargrid[(i,j+1,k)] - \scalargrid[(i,j,k)]}{\dy} \vectorgridcomponent{y,(i,j,k)} 
                   +                    \frac{\scalargrid[(i,j,k)] - \scalargrid[(i,j-1,k)]}{\dy} \vectorgridcomponent[y,(i,j-1,k)] \right)\\
            &\quad + \frac{1}{2} \left( \frac{\scalargrid[(i,j+1,k)] - \scalargrid[(i,j,k)]}{\dz} \vectorgridcomponent{z,(i,j,k)} 
                   +                    \frac{\scalargrid[(i,j,k)] - \scalargrid[(i,j-1,k)]}{\dz} \vectorgridcomponent[z,(i,j,k-1)] \right)\\
            &\quad + \scalargrid[(i,j,k)]
            \left( \frac{ \vectorgridcomponent{x,(i,j,k)} - \vectorgridcomponent[x,(i-1,j,k)]}{\dx}
                 + \frac{ \vectorgridcomponent{y,(i,j,k)} - \vectorgridcomponent[y,(i,j-1,k)]}{\dy} 
                 + \frac{ \vectorgridcomponent{z,(i,j,k)} - \vectorgridcomponent[z,(i,j,k-1)]}{\dz}\right)
        \end{split} \nonumber
        \\
        &= \indicesBig{\spacetransform{\faces}{\centers} \left(\vectorgrid \had \gradgrid \scalargrid \right)}{(i,j,k)} + \indicesbig{\scalargrid \had \divgrid \vectorgrid}{(i,j,k)} 
    \end{align}
    
    \item Similarly, one can show that the divergence of the product of a scalar and a tensor field, analogous to the continuous product rule $\div \left( s \boldsymbol{T} \right) = \boldsymbol{T} \cdot \grad s + s \div \boldsymbol{T}$, is
    \begin{align}
        \divgrid_{\tensorspace} \big( \spacetransform{\faces}{\tensorspace}\spacetransform{\centers}{\faces} \scalargrid \had \tensorgrid \big) 
        =
        \spacetransform{\tensorspace}{\faces} \Big( \tensorgrid \had \gradgrid_{\faces} \big( \spacetransform{\centers}{\faces} \scalargrid \big) \Big) + \spacetransform{\centers}{\faces} \scalargrid \had \divgrid_{\tensorspace} \tensorgrid
        \label{eq:scalar_tensor_product_rule}
    \end{align}
\end{itemize}

\subsection{Differential operators applied to composite fields}
\label{app:differential_operators_on_composite_fields}
The indicator fields on the different cell spaces are related as follows:
\begin{align}
    \hgridfaces^\pm &\coloneqq \spacetransform{\centers}{\faces} \hgridcenters^\pm \\
    \hgridtensorspace^\pm &\coloneqq \spacetransform{\faces}{\tensorspace} \hgridfaces^\pm
\end{align}
Because the one-dimensional interpolation operations commute on uniform grids (for example, $\spacetransform{x}{} \spacetransform{y}{} \scalargrid = \spacetransform{y}{} \spacetransform{x}{} \scalargrid $), it follows that $\hgridtensorspace^{\pm}= \left(\hgridtensorspace^{\pm} \right)^\top$. Similarly, because the one-dimensional interpolation and differentiation operations commute for any combination of grid directions on uniform grids, one can show that $\spacetransform{\faces}{\centers} \gradgrid \scalargrid = \divgrid \spacetransform{\centers}{\faces} \scalargrid$ and $\gradgrid_{\faces} \spacetransform{\centers}{\faces} \scalargrid = \left(\spacetransform{\faces}{\tensorspace} \gradgrid \scalargrid \right)^\top$
We can then work out the following operations:
\begin{itemize}
\item Gradient of a composite scalar field
\begin{equation}
    \gradgrid \mask{\scalargrid} = \hgridfaces^+ \had \gradgrid \scalargrid^+ + \hgridfaces^- \had \gradgrid \scalargrid^- + \gradgrid \hgridcenters^+ \had \spacetransform{\centers}{\faces}\left( \scalargrid^+ - \scalargrid^- \right)
    \label{eq:composite_scalar_grad}
\end{equation}
\item Gradient of a composite vector field
\begin{align}
    \gradgrid_{\faces} \mask{\vectorgrid} 
    &= 
    \hgridtensorspace^+ \had \gradgrid_{\faces} \vectorgrid^+ + \hgridtensorspace^- \had \gradgrid_{\faces} \vectorgrid^- + \gradgrid_{\faces} \hgridfaces^+ \had \spacetransform{\faces}{\tensorspace}\left( \vectorgrid^+ - \vectorgrid^- \right)
    \\
    &=
    \hgridtensorspace^+ \had \gradgrid_{\faces} \vectorgrid^+ + \hgridtensorspace^- \had \gradgrid_{\faces} \vectorgrid^- + \left(\spacetransform{\faces}{\tensorspace} \gradgrid \hgridcenters^+ \right)^\top \had \spacetransform{\faces}{\tensorspace}\left( \vectorgrid^+ - \vectorgrid^- \right)
    \label{eq:composite_vector_grad}
\end{align}

\item Divergence of a composite vector field
\begin{align}
    \divgrid \mask{\vectorgrid} &= \hgridcenters^+ \had \divgrid \vectorgrid^+ + \hgridcenters^- \had \divgrid \vectorgrid^- + \spacetransform{\faces}{\centers} \left( \gradgrid \hgridcenters^+ \had \left( \vectorgrid^+ - \vectorgrid^- \right) \right) 
    \label{eq:composite_vector_div}
\end{align}

\item Divergence of a composite tensor field
\begin{align}
    \divgrid_{\tensorspace} \mask{\tensorgrid} &= \hgridfaces^+ \had \divgrid_{\tensorspace} \tensorgrid^+ + \hgridfaces^- \had \divgrid_{\tensorspace} \tensorgrid^- + \spacetransform{\tensorspace}{\faces} \Big( \gradgrid_{\faces} \hgridfaces^+ \had \left( \tensorgrid^+ - \tensorgrid^- \right) \Big)
    \\
    &= \hgridfaces^+ \had \divgrid_{\tensorspace} \tensorgrid^+ + \hgridfaces^- \had \divgrid_{\tensorspace} \tensorgrid^- + \spacetransform{\tensorspace}{\faces}\left(\left(\spacetransform{\faces}{\tensorspace}\gradgrid \hgridcenters^+\right)^\top \had \left( \tensorgrid^+ - \tensorgrid^- \right) \right)
    \label{eq:composite_tensor_div}
\end{align}

\item Divergence of the gradient, or the Laplacian, of a composite scalar field (combining \eqref{eq:composite_scalar_grad} and \eqref{eq:composite_vector_div})
\begin{align}
    \lapgrid \mask{\scalargrid} &= \divgrid \gradgrid \mask{\scalargrid} \\
    &= \hgridcenters^+ \had \lapgrid \scalargrid^+ + \hgridcenters^- \had \lapgrid \scalargrid^-
     + \spacetransform{\faces}{\centers} \left( \gradgrid \hgridcenters^+ \had \left( \gradgrid \scalargrid^+ - \gradgrid \scalargrid^- \right) \right)
    + \divgrid \left( \gradgrid \hgridcenters^+ \had \spacetransform{\centers}{\faces}\left( \scalargrid^+ - \scalargrid^- \right) \right)
    \label{eq:composite_scalar_lap}
\end{align}

\item Laplacian of a composite vector field (combining \eqref{eq:composite_vector_grad} and \eqref{eq:composite_tensor_div})
\begin{align}
    \lapgrid_{\faces} \mask{\vectorgrid} &= \divgrid_{\tensorspace} \gradgrid_{\faces} \mask{\vectorgrid} \\
    &= \hgridfaces^+ \had \lapgrid_{\faces} \vectorgrid^+ + \hgridfaces^- \had \lapgrid_{\faces} \vectorgrid^-
    + \spacetransform{\tensorspace}{\faces}\left(\left(\spacetransform{\faces}{\tensorspace} \gradgrid \hgridcenters^+ \right)^\top \had  \left( \gradgrid_{\faces} \vectorgrid^+ - \gradgrid_{\faces} \vectorgrid^- \right)\right) \nonumber
    \\
    &\quad+ \divgrid_{\tensorspace}\left(\left(\spacetransform{\faces}{\tensorspace}\gradgrid \hgridcenters^+\right)^\top \had \spacetransform{\faces}{\tensorspace} \left( \vectorgrid^+ - \vectorgrid^- \right) \right)
    \label{eq:composite_vector_lap}
\end{align}

\item Time derivative of a composite velocity field
\begin{equation}
    \der{\mask{\vectorgrid}}{t}{} = \hgridfaces^+ \had \der{\vectorgrid^+}{t}{} + \hgridfaces^- \had \der{\vectorgrid^-}{t}{} + \der{\hgridfaces^+}{t}{} \had \left( \vectorgrid^+ - \vectorgrid^- \right) 
    \label{eq:composite_vector_time_der}
\end{equation}



\item Convective derivative of the composite velocity field in divergence form using Eq.~\eqref{eq:interpolation_of_vector_product}, Eq.~\eqref{eq:scalar_tensor_product_rule}, $(\hgridtensorspace^{\pm})^\top = \hgridtensorspace^{\pm}$ and $(\hgridtensorspace^{\pm})^2 = \hgridtensorspace^{\pm} (1 - \hgridtensorspace^{\pm})$. Here, we also use the notation $\otimes$ to indicate an outer product.
\begin{align}
    \begin{split}
        &\divgrid_{\tensorspace} \left(  \mask{\vectorgrid} \otimes \mask{\vectorgrid} \right) 
        \\
        &= 
        \divgrid_{\tensorspace} \left( \left( \spacetransform{\faces}{\tensorspace} \mask{\vectorgrid} \right)^\top \had \spacetransform{\faces}{\tensorspace} \mask{\vectorgrid} \right) 
        \\
        &= 
        \divgrid_{\tensorspace} \left( \hgridtensorspace^+ \had  \left( \spacetransform{\faces}{\tensorspace} \vectorgrid^+ \right)^\top \had \hgridtensorspace^+ \had \spacetransform{\faces}{\tensorspace} \vectorgrid^+  \right)
        +
        \divgrid_{\tensorspace} \left( \hgridtensorspace^+ \had  \left( \spacetransform{\faces}{\tensorspace} \vectorgrid^+ \right)^\top \had \hgridtensorspace^- \had \spacetransform{\faces}{\tensorspace} \vectorgrid^-  \right)
        \\
        &\quad+
        \divgrid_{\tensorspace} \left( \hgridtensorspace^- \had  \left( \spacetransform{\faces}{\tensorspace} \vectorgrid^- \right)^\top \had \hgridtensorspace^+ \had \spacetransform{\faces}{\tensorspace} \vectorgrid^+  \right)
        +
        \divgrid_{\tensorspace} \left( \hgridtensorspace^- \had  \left( \spacetransform{\faces}{\tensorspace} \vectorgrid^- \right)^\top \had \hgridtensorspace^- \had \spacetransform{\faces}{\tensorspace} \vectorgrid^-  \right)
        \\
        &= \divgrid_{\tensorspace} \left( \hgridtensorspace^+ \had \left( \spacetransform{\faces}{\tensorspace} \vectorgrid^+ \right)^\top \had \spacetransform{\faces}{\tensorspace} \vectorgrid^+ + \hgridtensorspace^- \had \left( \spacetransform{\faces}{\tensorspace} \vectorgrid^- \right)^\top \had \spacetransform{\faces}{\tensorspace} \vectorgrid^-  \right)
        \\
        &\quad-
        \divgrid_{\tensorspace} \left( \hgridtensorspace^+ \had \hgridtensorspace^- \had \left( \spacetransform{\faces}{\tensorspace} \left( \vectorgrid^+ - \vectorgrid^- \right) \right)^\top \had \spacetransform{\faces}{\tensorspace} \left( \vectorgrid^+ - \vectorgrid^- \right)  \right)
        \\
        &\quad+
        \divgrid_{\tensorspace} 
        \left( 
            \spacetransform{\faces}{\tensorspace}\mask{\vectorgrid}^\top \had \left(\frac{\dxgridtensor^2}{4} \had  \gradgrid_{\faces} \hgridfaces^+ \had \left( \gradgrid_{\faces} \vectorgrid^+ - \gradgrid_{\faces} \vectorgrid^- \right) \right)
        \right)
        \\
        &\quad+
        \divgrid_{\tensorspace} 
        \left( 
            \spacetransform{\faces}{\tensorspace}\mask{\vectorgrid} \had \left(\frac{\dxgridtensor^2}{4} \had  \gradgrid_{\faces} \hgridfaces^+ \had \left( \gradgrid_{\faces} \vectorgrid^+ - \gradgrid_{\faces} \vectorgrid^- \right) \right)^\top
        \right)
        \\
        &\quad-
        \divgrid_{\tensorspace} \left( 
            \left(\frac{\dxgridtensor^2}{4} \had  \gradgrid_{\faces} \hgridfaces^+ \had \left( \gradgrid_{\faces} \vectorgrid^+ - \gradgrid_{\faces} \vectorgrid^- \right) \right)^\top \had
            \left(\frac{\dxgridtensor^2}{4} \had  \gradgrid_{\faces} \hgridfaces^+ \had \left( \gradgrid_{\faces} \vectorgrid^+ - \gradgrid_{\faces} \vectorgrid^- \right) \right)
        \right)
                \\
        &= \hgridfaces^+ \had \divgrid_{\tensorspace} \left( \vectorgrid^+ \otimes \vectorgrid^+\right) + \hgridfaces^- \had \divgrid_{\tensorspace} \left( \vectorgrid^-  \otimes \vectorgrid^-  \right)
        \\
        &\quad+ \spacetransform{\tensorspace}{\faces} 
        \Bigg( 
            \gradgrid_{\faces} \hgridfaces^+
            \had 
            \Big( \left( \vectorgrid^+  \otimes \vectorgrid^+  \right) - \left( \vectorgrid^-  \otimes \vectorgrid^-  \right) \Big)
        \Bigg)
        \\
        &\quad-
        \divgrid_{\tensorspace} 
        \Bigg( 
            \hgridtensorspace^+ \had \hgridtensorspace^- \had \left( \vectorgrid^+ - \vectorgrid^- \right) \otimes \left( \vectorgrid^+ - \vectorgrid^- \right)  
        \Bigg)
        \\
        &\quad+
        \divgrid_{\tensorspace} 
        \left( 
            \spacetransform{\faces}{\tensorspace} \mask{\vectorgrid}^\top \had \left(\frac{\dxgridtensor^2}{4} \had  \gradgrid_{\faces} \hgridfaces^+ \had \left( \gradgrid_{\faces} \vectorgrid^+ - \gradgrid_{\faces} \vectorgrid^- \right) \right)
        \right)
        \\
        &\quad+
        \divgrid_{\tensorspace} 
        \left( 
            \spacetransform{\faces}{\tensorspace} \mask{\vectorgrid} \had \left(\frac{\dxgridtensor^2}{4} \had  \gradgrid_{\faces} \hgridfaces^+ \had \left( \gradgrid_{\faces} \vectorgrid^+ - \gradgrid_{\faces} \vectorgrid^- \right) \right)^\top
        \right)
        \\
        &\quad-
        \divgrid_{\tensorspace} \left( 
            \left(\frac{\dxgridtensor^2}{4} \had  \gradgrid_{\faces} \hgridfaces^+ \had \left( \gradgrid_{\faces} \vectorgrid^+ - \gradgrid_{\faces} \vectorgrid^- \right) \right)^\top \had
            \left(\frac{\dxgridtensor^2}{4} \had  \gradgrid_{\faces} \hgridfaces^+ \had \left( \gradgrid_{\faces} \vectorgrid^+ - \gradgrid_{\faces} \vectorgrid^- \right) \right)
        \right)
    \end{split}
    \label{eq:composite_vector_conv}
\end{align}
\end{itemize}


\subsection{Taylor series construction of discrete fields and their gradients}
\label{app:taylor_series}
Let $\scalarpoint, \scalarpointnormal, \scalarpointtanone, \scalarpointtantwo \in \R^{\spoints}$ denote the discrete approximations to $s$ and its normal and tangential derivatives on IB points, respectively. Similarly, let $\velpoint, \velpointnormal, \velpointtanone, \velpointtantwo \in \R^{\vpoints}$ denote the discrete approximations to $v$ and its normal and tangential derivatives on IB points, respectively. Then we can formulate the following Taylor series:
\begin{itemize}
\item Taylor series construction of the scalar field $\scalargrid$ about the point $\xsurfvec[l]$
\begin{equation}
    \indicesbig{(\scalargrid^\pm)^{\mathrm{TS}}_l}{(i,j,k)}
    \approx 
    \scalarpointsurfplusminus[,l]
    + \normalvec[l] \cdot \big( \xgridvec[\centers,(i,j,k)] - \xsurfvec[l] \big)\scalarpointsurfnormalplusminus[,l]
    + \tanonevec[l] \cdot \big( \xgridvec[\centers,(i,j,k)] - \xsurfvec[l] \big)\scalarpointsurftanoneplusminus[,l]
    + \tantwovec[l] \cdot \big( \xgridvec[\centers,(i,j,k)] - \xsurfvec[l] \big)\scalarpointsurftantwoplusminus[,l]
    ,
    \label{eq:s_TS}
\end{equation}
\item Gradient of the Taylor series construction of the scalar field $\scalargrid$ about the point $\xsurfvec[l]$
\begin{equation}
    \indicesbig{\gradgrid (\scalargrid^\pm)^{\mathrm{TS}}_l}{(i,j,k)}
    \approx 
    \begin{pmatrix}
    \scalarpointsurfnormalplusminus[,l] \normalvecc{x,l} + \scalarpointsurftanoneplusminus[,l] \tanonevecc{x,l} + \scalarpointsurftantwoplusminus[,l] \tantwovecc{x,l}
    \\
    \scalarpointsurfnormalplusminus[,l] \normalvecc{y,l} + \scalarpointsurftanoneplusminus[,l] \tanonevecc{y,l} + \scalarpointsurftantwoplusminus[,l] \tantwovecc{y,l}
    \\
    \scalarpointsurfnormalplusminus[,l] \normalvecc{z,l} + \scalarpointsurftanoneplusminus[,l] \tanonevecc{z,l} + \scalarpointsurftantwoplusminus[,l] \tantwovecc{z,l}
    \end{pmatrix},
    \label{eq:Gs_TS}
\end{equation}
\item Taylor series construction of the vector field $\velgrid$ about the point $\xsurfvec[l]$
\begin{equation}
    \indicesbig{(\velgrid^\pm)^{\mathrm{TS}}_l}{(i,j,k)}
    \approx 
    \left(\begin{smallmatrix}
    \indices{\velcomponentpointsurfplusminus{x}}{l}
    + \normalvec[l] \cdot \big( \xgridvec[\facescomponent_x, (i,j,k)] - \xsurfvec[l] \big)\indices{\velcomponentpointsurfnormalplusminus{x}}{l}
    + \tanonevec[l] \cdot \big( \xgridvec[\facescomponent_x, (i,j,k)] - \xsurfvec[l] \big)\indices{\velcomponentpointsurftanoneplusminus{x}}{l}
    + \tantwovec[l] \cdot \big( \xgridvec[\facescomponent_x, (i,j,k)] - \xsurfvec[l] \big)\indices{\velcomponentpointsurftantwoplusminus{x}}{l}
    \\
    \indices{\velcomponentpointsurfplusminus{y}}{l}
    + \normalvec[l] \cdot \big( \xgridvec[\facescomponent_y, (i,j,k)] - \xsurfvec[l] \big)\indices{\velcomponentpointsurfnormalplusminus{y}}{l}
    + \tanonevec[l] \cdot \big( \xgridvec[\facescomponent_y, (i,j,k)] - \xsurfvec[l] \big)\indices{\velcomponentpointsurftanoneplusminus{y}}{l}
    + \tantwovec[l] \cdot \big( \xgridvec[\facescomponent_y, (i,j,k)] - \xsurfvec[l] \big)\indices{\velcomponentpointsurftantwoplusminus{y}}{l}
    \\
    \indices{\velcomponentpointsurfplusminus{z}}{l}
    + \normalvec[l] \cdot \big( \xgridvec[\facescomponent_z, (i,j,k)] - \xsurfvec[l] \big)\indices{\velcomponentpointsurfnormalplusminus{z}}{l}
    + \tanonevec[l] \cdot \big( \xgridvec[\facescomponent_z, (i,j,k)] - \xsurfvec[l] \big)\indices{\velcomponentpointsurftanoneplusminus{z}}{l}
    + \tantwovec[l] \cdot \big( \xgridvec[\facescomponent_z, (i,j,k)] - \xsurfvec[l] \big)\indices{\velcomponentpointsurftantwoplusminus{z}}{l}
    \end{smallmatrix}\right),
    \label{eq:v_TS}
\end{equation}
\item Gradient of the Taylor series construction of the vector field $\velgrid$ about the point $\xsurfvec[l]$
\begin{equation}
    \indicesbig{\gradgrid_{\faces}(\velgrid^\pm)^{\mathrm{TS}}_l}{(i,j,k)} \approx 
    \left(\begin{smallmatrix}
    \indices{\velcomponentpointsurfnormalplusminus{x}}{l} \normalvecc{x,l} +
    \indices{\velcomponentpointsurftanoneplusminus{x}}{l} \tanonevecc{x,l} +
    \indices{\velcomponentpointsurftantwoplusminus{x}}{l} \tantwovecc{x,l}
    &
    \indices{\velcomponentpointsurfnormalplusminus{x}}{l} \normalvecc{y,l} +
    \indices{\velcomponentpointsurftanoneplusminus{x}}{l} \tanonevecc{y,l} +
    \indices{\velcomponentpointsurftantwoplusminus{x}}{l} \tantwovecc{y,l}
    &
    \indices{\velcomponentpointsurfnormalplusminus{x}}{l} \normalvecc{z,l} +
    \indices{\velcomponentpointsurftanoneplusminus{x}}{l} \tanonevecc{z,l} +
    \indices{\velcomponentpointsurftantwoplusminus{x}}{l} \tantwovecc{z,l}
    \\
    \indices{\velcomponentpointsurfnormalplusminus{y}}{l} \normalvecc{x,l} +
    \indices{\velcomponentpointsurftanoneplusminus{y}}{l} \tanonevecc{x,l} +
    \indices{\velcomponentpointsurftantwoplusminus{y}}{l} \tantwovecc{x,l}
    &
    \indices{\velcomponentpointsurfnormalplusminus{y}}{l} \normalvecc{y,l} +
    \indices{\velcomponentpointsurftanoneplusminus{y}}{l} \tanonevecc{y,l} +
    \indices{\velcomponentpointsurftantwoplusminus{y}}{l} \tantwovecc{y,l}
    &
    \indices{\velcomponentpointsurfnormalplusminus{y}}{l} \normalvecc{z,l} +
    \indices{\velcomponentpointsurftanoneplusminus{y}}{l} \tanonevecc{z,l} +
    \indices{\velcomponentpointsurftantwoplusminus{y}}{l} \tantwovecc{z,l}
    \\
    \indices{\velcomponentpointsurfnormalplusminus{z}}{l} \normalvecc{x,l} +
    \indices{\velcomponentpointsurftanoneplusminus{z}}{l} \tanonevecc{x,l} +
    \indices{\velcomponentpointsurftantwoplusminus{z}}{l} \tantwovecc{x,l}
    &
    \indices{\velcomponentpointsurfnormalplusminus{y}}{l} \normalvecc{y,l} +
    \indices{\velcomponentpointsurftanoneplusminus{y}}{l} \tanonevecc{y,l} +
    \indices{\velcomponentpointsurftantwoplusminus{y}}{l} \tantwovecc{y,l}
    &
    \indices{\velcomponentpointsurfnormalplusminus{z}}{l} \normalvecc{z,l} +
    \indices{\velcomponentpointsurftanoneplusminus{z}}{l} \tanonevecc{z,l} +
    \indices{\velcomponentpointsurftantwoplusminus{z}}{l} \tantwovecc{z,l}
    \end{smallmatrix} \right),
    \label{eq:Gv_TS}
\end{equation}
\end{itemize}

\subsection{Regularization operations}
\label{app:regularization}

\begin{itemize}
\item Regularization of scalar-valued surface data to cell centers
\begin{equation}
    \indicesbig{\regds_{\centers} \scalarpoint}{(i,j,k)} \coloneqq \sum_l \indicesbig{\dddf_{\centers,l}}{(i,j,k)} \scalarpoint[l] \surfarea[l]
    \label{eq:RC}
\end{equation}
\item Regularization of vector-valued surface data to cell faces
\begin{equation}
    \indicesbig{\regds_{\faces} \vectorpoint}{(i,j,k)} \coloneqq \sum_l
    \begin{pmatrix}
        \indicesbig{\dddf_{\facescomponent_x,l}}{(i,j,k)} \vectorpointcomponent[x,l] \\
        \indicesbig{\dddf_{\facescomponent_y,l}}{(i,j,k)} \vectorpointcomponent[y,l] \\
        \indicesbig{\dddf_{\facescomponent_z,l}}{(i,j,k)} \vectorpointcomponent[z,l]
    \end{pmatrix} \surfarea[l]
    \label{eq:RF}
\end{equation}
\item Normal-distance-weighted regularization of vector-valued surface data to cell faces
\begin{equation}
    \indicesbig{\regds_{\faces, 1\normalbase} (\vectorpoint)}{(i,j,k)} \coloneqq \sum_l
    \begin{pmatrix}
        \indicesbig{\dddf_{\facescomponent_x,l}}{(i,j,k)} \normalvec[l] \cdot \big( \xgridvec[\facescomponent_x] - \xsurfvec[l] \big) \vectorpointcomponent[x,l] \\
        \indicesbig{\dddf_{\facescomponent_y,l}}{(i,j,k)} \normalvec[l] \cdot \big( \xgridvec[\facescomponent_y] - \xsurfvec[l] \big) \vectorpointcomponent[y,l] \\
        \indicesbig{\dddf_{\facescomponent_z,l}}{(i,j,k)} \normalvec[l] \cdot \big( \xgridvec[\facescomponent_z] - \xsurfvec[l] \big) \vectorpointcomponent[z,l]
    \end{pmatrix} \surfarea[l]
    \label{eq:RF1n}
\end{equation}
\item Regularization of tensor-valued surface data to the tensorspace using the transpose of the DDFs evaluated at the faces and interpolated to the tensorspace
\begin{equation}
    \indicesbig{\regds_{(\spacetransform{}{}\faces)^\top} \tpoint{T}}{(i,j,k)} \coloneqq \sum_l
    \begin{pmatrix}
    \indicesbig{\spacetransform{x}{} \dddf_{\facescomponent_x,l}}{(i,j,k)} \tensorpointcomponent[11,l]
    &
    \indicesbig{\spacetransform{x}{} \dddf_{\facescomponent_y,l}}{(i,j,k)} \tensorpointcomponent[12,l]
    &
    \indicesbig{\spacetransform{x}{} \dddf_{\facescomponent_z,l}}{(i,j,k)} \tensorpointcomponent[13,l]
    \\
    \indicesbig{\spacetransform{y}{} \dddf_{\facescomponent_x,l}}{(i,j,k)} \tensorpointcomponent[21,l]
    &
    \indicesbig{\spacetransform{y}{} \dddf_{\facescomponent_y,l}}{(i,j,k)} \tensorpointcomponent[22,l]
    &
    \indicesbig{\spacetransform{y}{} \dddf_{\facescomponent_z,l}}{(i,j,k)} \tensorpointcomponent[23,l]
    \\
    \indicesbig{\spacetransform{z}{} \dddf_{\facescomponent_x,l}}{(i,j,k)} \tensorpointcomponent[31,l]
    &
    \indicesbig{\spacetransform{z}{} \dddf_{\facescomponent_y,l}}{(i,j,k)} \tensorpointcomponent[32,l]
    &
    \indicesbig{\spacetransform{z}{} \dddf_{\facescomponent_z,l}}{(i,j,k)} \tensorpointcomponent[33,l]
    \end{pmatrix} \surfarea[l]
    \label{eq:RIFT}
\end{equation}
\item Normal-distance-weighted regularization of tensor-valued surface data to tensorspace using the transpose of the DDFs evaluated at the faces and interpolated to the tensorspace
\begin{equation}
    \indicesbig{\regds_{(\spacetransform{}{}\faces)^\top, 1 \normalbase} \tpoint{T}}{(i,j,k)}\coloneqq
    \sum_l
    \left(\begin{smallmatrix}
    \indicesbig{\spacetransform{x}{} \dddf_{\facescomponent_x,l}}{(i,j,k)} \normalvec[l] \cdot \big( \xgridvec[\centers,(i,j,k)] - \xsurfvec[l] \big) \tensorpointcomponent[11,l]
    &
    \indicesbig{\spacetransform{x}{} \dddf_{\facescomponent_y,l}}{(i,j,k)} \normalvec[l] \cdot \big( \xgridvec[\edgescomponent_z,(i,j,k)] - \xsurfvec[l] \big) \tensorpointcomponent[12,l]
    &
    \indicesbig{\spacetransform{x}{} \dddf_{\facescomponent_z,l}}{(i,j,k)} \normalvec[l] \cdot \big( \xgridvec[\edgescomponent_y,(i,j,k)] - \xsurfvec[l] \big) \tensorpointcomponent[13,l]
    \\
    \indicesbig{\spacetransform{y}{} \dddf_{\facescomponent_x,l}}{(i,j,k)} \normalvec[l] \cdot \big( \xgridvec[\edgescomponent_z,(i,j,k)] - \xsurfvec[l] \big) \tensorpointcomponent[21,l]
    &
    \indicesbig{\spacetransform{y}{} \dddf_{\facescomponent_y,l}}{(i,j,k)} \normalvec[l] \cdot \big( \xgridvec[\centers,(i,j,k)] - \xsurfvec[l] \big) \tensorpointcomponent[22,l]
    &
    \indicesbig{\spacetransform{y}{} \dddf_{\facescomponent_z,l}}{(i,j,k)} \normalvec[l] \cdot \big( \xgridvec[\edgescomponent_x,(i,j,k)] - \xsurfvec[l] \big) \tensorpointcomponent[23,l]
    \\
    \indicesbig{\spacetransform{z}{} \dddf_{\facescomponent_x,l}}{(i,j,k)} \normalvec[l] \cdot \big( \xgridvec[\edgescomponent_y,(i,j,k)] - \xsurfvec[l] \big) \tensorpointcomponent[31,l]
    &
    \indicesbig{\spacetransform{z}{} \dddf_{\facescomponent_y,l}}{(i,j,k)} \normalvec[l] \cdot \big( \xgridvec[\edgescomponent_x,(i,j,k)] - \xsurfvec[l] \big) \tensorpointcomponent[32,l]
    &
    \indicesbig{\spacetransform{z}{} \dddf_{\facescomponent_z,l}}{(i,j,k)} \normalvec[l] \cdot \big( \xgridvec[\centers,(i,j,k)] - \xsurfvec[l] \big) \tensorpointcomponent[33,l]
    \end{smallmatrix}\right) \surfarea[l]
    \label{eq:RIFT1n}
\end{equation}
\end{itemize}

\subsection{Interpolation operations}
\label{app:interpolation}

\begin{itemize}
\item Interpolation of scalar-valued grid data on cell centers to scalar-valued surface data
\begin{equation}
    \indicesbig{\interpbody_{\centers} \scalargrid}{l}
    \coloneqq
    \dx\dy\dz \sum_{i,j,k} \indicesbig{\dddf_{\centers,l}}{(i,j,k)} \scalargrid[(i,j,k)]
    \label{eq:EC}
\end{equation}
\item Interpolation of vector-valued grid data on cell faces to vector-valued surface data
\begin{equation}
    \indicesbig{\interpbody_{\faces} \vectorgrid}{l}
    \coloneqq
    \dx\dy\dz \sum_{i,j,k}
    \begin{pmatrix}
        \indices{\dddf_{\facescomponent_x,l}}{(i,j,k)} \vectorgridcomponent{x,(i,j,k)} \\
        \indices{\dddf_{\facescomponent_y,l}}{(i,j,k)} \vectorgridcomponent{y,(i,j,k)} \\
        \indices{\dddf_{\facescomponent_z,l}}{(i,j,k)} \vectorgridcomponent{z,(i,j,k)}
    \end{pmatrix}
    \label{eq:EF}
\end{equation}
\item Normal-distance-weighted interpolation of scalar-valued grid data on cell centers to scalar-valued surface data 
\begin{equation}
    \indicesbig{\interpbody_{\centers, 1\normalbase}  \scalargrid}{l} \coloneqq \dx\dy\dz \sum_{i,j,k} \indicesbig{\dddf_{\centers,l}}{(i,j,k)} \normalvec[l] \cdot \big( \xgridvec[\centers,(i,j,k)] - \xsurfvec[l] \big) \scalargrid[(i,j,k)]
    \label{eq:EC1n}
\end{equation}
\item Normal-distance-weighted interpolation of vector-valued grid data on cell faces to vector-valued surface data 
\begin{equation}
    \indicesbig{\interpbody_{\faces, 1\normalbase}  \vectorgrid}{l} \coloneqq \dx\dy\dz \sum_{i,j,k} 
    \begin{pmatrix}
        \indicesbig{\dddf_{\facescomponent_x,l}}{(i,j,k)} \normalvec[l] \cdot \big( \xgridvec[\facescomponent_x, (i,j,k)] - \xsurfvec[l] \big) \vectorgridcomponent{x,(i,j,k)} \\
        \indicesbig{\dddf_{\facescomponent_y,l}}{(i,j,k)} \normalvec[l] \cdot \big( \xgridvec[\facescomponent_y, (i,j,k)] - \xsurfvec[l] \big) \vectorgridcomponent{y,(i,j,k)} \\
        \indicesbig{\dddf_{\facescomponent_z,l}}{(i,j,k)} \normalvec[l] \cdot \big( \xgridvec[\facescomponent_z, (i,j,k)] - \xsurfvec[l] \big) \vectorgridcomponent{z,(i,j,k)}
    \end{pmatrix}
    \label{eq:EF1n}
\end{equation}

\end{itemize}

\section{Solution of general saddle-point systems}
\label{app:saddlesystems}

A general block system (with non-singular matrix $\basicA$) can be factorized as follows:
\begin{equation}
\begin{bmatrix}
\basicA & \basicBoneT \\ \basicBtwo & -\basicC
\end{bmatrix} = \begin{bmatrix}
\basicA & 0\\ \basicBtwo & \basicS
\end{bmatrix}
\begin{bmatrix}
\basicI & \basicA^{-1}\basicBoneT\\ 0 & \basicI
\end{bmatrix},
\end{equation}
where
\begin{equation}
\label{eq:schur}
\basicS \equiv -\basicC - \basicBtwo\basicA^{-1}\basicBoneT
\end{equation}
is the {\em Schur complement} of the matrix system and $\basicI$ is the identity. By this decomposition, one can develop the \emph{Schur complement reduction} algorithm~\citep{benzi_numerical_2005} for the solution of the block system
\begin{equation}
\label{eq:saddlept}
\begin{bmatrix}
\basicA & \basicBoneT \\ \basicBtwo & -\basicC
\end{bmatrix} \begin{pmatrix} \basicx \\ \basicy \end{pmatrix} = \begin{pmatrix} \basicrone \\ \basicrtwo \end{pmatrix}.
\end{equation}
We will refer to $\basicx$ as the solution vector and $\basicy$ as the constraint force. We define the intermediate solution vector $(\basicx^{*},\, \basicy^{*})^\top$ as the solution of the lower-triangular system
\begin{equation}
\begin{bmatrix}
\basicA & 0\\ \basicBtwo & \basicS
\end{bmatrix} \begin{pmatrix} \basicx^{*} \\ \basicy^{*} \end{pmatrix} = \begin{pmatrix} \basicrone \\ \basicrtwo \end{pmatrix} 
\end{equation}
and then the solution we seek can be found by back substitution of
\begin{equation}
\begin{bmatrix}
\basicI & \basicA^{-1}\basicBoneT\\ 0 & \basicI
\end{bmatrix} \begin{pmatrix} \basicx \\ \basicy \end{pmatrix} = \begin{pmatrix} \basicx^{*} \\ \basicy^{*} \end{pmatrix} 
\end{equation}
The algorithm we derive from this is
\begin{align}
\basicA \basicx^{*} &= \basicrone, \nonumber\\
\basicS \basicy^{*} &=  \basicrtwo - \basicBtwo \basicx^{*}, \label{eq:algor}\\
\basicy &= \basicy^{*}, \nonumber\\
\basicx &= \basicx^{*} - \basicA^{-1}\basicBoneT \basicy \nonumber.
\end{align}
It is also useful to have an inverse representation of the block matrix system:
\begin{equation}
\label{eq:saddleinverse}
\begin{pmatrix} \basicx \\ \basicy \end{pmatrix} = \begin{bmatrix} \basicA^{-1} + \basicA^{-1} \basicBoneT\basicS^{-1} \basicBtwo \basicA^{-1} & -\basicA^{-1}\basicBoneT \basicS^{-1} \\ -\basicS^{-1} \basicBtwo \basicA^{-1} & \basicS^{-1}\end{bmatrix} \begin{pmatrix} \basicrone \\ \basicrtwo \end{pmatrix}. 
\end{equation}

For a more detailed discussion on the numerical solution of saddle-point systems, the reader is referred to~\citet{benzi_numerical_2005} and other related works cited therein.



\section{Derivation of the incompressible Navier--Stokes equations for composite solutions}
\label{app:ns_ib_derivation}

We start the derivation of the incompressible Navier--Stokes equations for the composite velocity field $\vel \coloneqq \heavi^+ \vel^+ + \heavi^- \vel^-$ and pressure field $p \coloneqq \heavi^+ p^+ + \heavi^- p^-$ by adding together the Navier--Stokes equations for the exterior and interior solutions, multiplied by their indicator fields:
\begin{gather}
    \heavi^+ \left(\frac{\partial \vel^+}{\partial t} + \div ( \vel^+ \vel^+ )\right)
    +
    \heavi^- \left(\frac{\partial \vel^-}{\partial t} + \div ( \vel^- \vel^- )\right)
    =
    \heavi^+ \left( -\grad p^+ + \frac{1}{\Re} \lap \vel^+ \right)
    +
    \heavi^- \left( -\grad p^- + \frac{1}{\Re} \lap \vel^- \right)
    \\
    \heavi^+ \left( \div \vel^+ \right)
    +
    \heavi^- \left( \div \vel^- \right)
    =
    0.
\end{gather}

Then we apply the product rules for the composite velocity, noting that $\grad\heavi^- = - \grad \heavi^+$, $\heavi^\pm \heavi^\pm = \heavi^\pm$, and $\heavi^\pm \heavi^\mp = 0$,
\begin{align}
\grad ( \heavi^+ p^+ + \heavi^- p^- ) 
&= 
\heavi^+ \grad p^+ + \heavi^- \grad p^- + (p^+ - p^-) \grad \heavi^+,
\\
\grad ( \heavi^+ \vel^+ + \heavi^- \vel^- ) 
&= 
\heavi^+ \grad \vel^+ + \heavi^- \grad \vel^- + (\vel^+ - \vel^-) \grad \heavi^+,
\\
\lap ( \heavi^+ \vel^+ + \heavi^- \vel^- ) 
&= \div \grad ( \heavi^+ \vel^+ + \heavi^- \vel^- ) \nonumber \\
&= 
\heavi^+ \lap \vel^+ + \heavi^- \lap \vel^- + \grad \heavi^+ \cdot (\grad \vel^+ - \grad \vel^-) + \div \left( \grad \heavi^+  (\vel^+ - \vel^-)\right),
\\
\div \left( ( \heavi^+ \vel^+ + \heavi^- \vel^- ) ( \heavi^+ \vel^+ + \heavi^- \vel^- ) \right)
&= 
\div (  
\heavi^+\heavi^+ \vel^+\vel^+ + 
\heavi^+\heavi^- \vel^+\vel^- + 
\heavi^-\heavi^+ \vel^-\vel^+ + 
\heavi^-\heavi^+ \vel^-\vel^-
) \nonumber\\
&=
\div ( \heavi^+ \vel^+\vel^+ ) +
\div ( \heavi^- \vel^-\vel^- ) \nonumber\\
&=
\heavi^+ \div (\vel^+\vel^+) + 
\heavi^- \div (\vel^-\vel^-) + 
\grad \heavi^+ \cdot ( \vel^+\vel^+ - \vel^-\vel^- )
\\
\der{}{t}{}( \heavi^+ \vel^+ + \heavi^- \vel^- ) 
&=
\heavi^+ \der{\vel^+}{t}{} +
\heavi^- \der{\vel^-}{t}{} + 
(\vel^+ - \vel^-) \der{\heavi^+}{t}{}
\\
\div ( \heavi^+ \vel^+ + \heavi^- \vel^- ) 
&= 
\heavi^+ \div \vel^+ + \heavi^- \div \vel^- + (\vel^+ - \vel^-) \grad \heavi^+,
\end{align}
and we substitute the expression for the time derivative of the indicator field~\cite{juric_computations_1998, eldredge_method_2022},
\begin{equation}
    \der{\heavi^+}{t}{} = \grad \heavi^+ \cdot \Xdot
\end{equation}
to obtain
\begin{gather}
    \frac{\partial \vel}{\partial t} + \div ( \vel \vel ) 
    - \vecforce_{\inlder{\vel}{t}{}}
    - \vecforce_{\div ( \vel \vel )}
    =
    -\grad p + \frac{1}{\Re} \lap \vel
    + \vecforce_{\grad p}
    - \frac{1}{\Re} \vecforce_{\lap \vel}, \quad \x \in \domain,
    \\
    \div \vel 
    - \scaforce_{\div ( \vel )}
    =
    0,
\end{gather}
where
\begin{align*}
    \vecforce_{\grad p}
    &= 
    \grad \heavi^+ ( p^+ - p^- ),
    \\
    \vecforce_{\lap \vel} 
    &= 
    \grad \heavi^+ \cdot ( \grad \vel^+ - \grad \vel^- ) + \div \big( \grad \heavi^+ (\vel^+ - \vel^- ) \big),
    \\
    \vecforce_{\div ( \vel \vel )}
    &= \grad \heavi^+ \cdot (\vel^+ \vel^+ - \vel^- \vel^-),
    \\
    \vecforce_{\inlder{\vel}{t}{}}
    &= \big( \grad \heavi^+ \cdot \Xdot \big) (\vel^+ - \vel^-),
    \\
    \scaforce_{\div \vel} 
    &=  \grad \heavi^+ \cdot ( \vel^+ - \vel^- ).
\end{align*}

\bibliographystyle{unsrtnat}
\bibliography{references}

\end{document}

%% file: macros.tex
\DeclareMathAlphabet{\mathsfit}{T1}{\sfdefault}{\mddefault}{\sldefault}
\SetMathAlphabet{\mathsfit}{bold}{T1}{\sfdefault}{\bfdefault}{\sldefault}

\newcommand{\domain}{\Omega}
\newcommand{\surf}{\Gamma}

\newcommand{\mask}[1]{\overline{#1}}
\newcommand{\jump}[2][\surf]{[{#2}]_{#1}}

\newcommand{\normal}{\bm{n}}
\newcommand{\tangentone}{\bm{t}}
\newcommand{\tangenttwo}{\bm{b}}

\renewcommand{\Re}{\mathrm{Re}}

\newcommand{\x}{\boldsymbol{x}}
\newcommand{\y}{\boldsymbol{y}}
\newcommand{\X}{\boldsymbol{X}}
\newcommand{\Xdot}{\dot{\boldsymbol{X}}}
\newcommand{\vel}{\boldsymbol{v}}
\newcommand{\pot}{u}

\newcommand{\source}{q}

\newcommand{\scaforcedist}{f}

\newcommand{\scaforce}{F}
\newcommand{\vecforce}{\boldsymbol{F}}
\newcommand{\heavi}{H}

\newcommand{\R}{\mathbb{R}}
\newcommand{\der}[3]{\frac{\mathrm{d}^{#3}#1}{\mathrm{d}#2^{#3}}}
\newcommand{\inlder}[3]{\mathrm{d}^{#3}#1/\mathrm{d}#2^{#3}}

\newcommand{\grad}{\nabla}
\renewcommand{\div}{\nabla\cdot}

\newcommand{\lap}{\nabla^{2}}
\newcommand{\dvol}{\mathrm{d}\x}

\newcommand{\surfc}[1][]{\boldsymbol{\xi}_{#1}}
\newcommand{\dsurf}{\mathrm{d} S(\surfc)}

\newcommand{\grid}[1]{\mathsf{#1}}
\newcommand{\pointwise}[1]{\mathrm{#1}}
\newcommand{\spoint}[1]{\pointwise{#1}}
\newcommand{\vpoint}[1]{\bm{\pointwise{#1}}}
\newcommand{\tpoint}[1]{\bm{\pointwise{#1}}}

\newcommand{\indices}[2]{( #1 )_{#2}}
\newcommand{\indicesbig}[2]{\big( #1 \big)_{#2}}
\newcommand{\indicesBig}[2]{\Big( #1 \Big)_{#2}}

\newcommand{\dt}{\Delta t}
\newcommand{\timestepcoef}{K}

\newcommand{\numpts}{N}

\newcommand{\xgrid}[1][]{\grid{x}_{#1}}
\newcommand{\ygrid}[1][]{\grid{y}_{#1}}
\newcommand{\zgrid}[1][]{\grid{z}_{#1}}
\newcommand{\xgridvec}[1][]{\bm{\grid{x}}_{#1}}

\newcommand{\xsurf}[1][]{\spoint{X}_{#1}}
\newcommand{\ysurf}[1][]{\spoint{Y}_{#1}}
\newcommand{\zsurf}[1][]{\spoint{Z}_{#1}}
\newcommand{\xsurfvec}[1][]{\vpoint{X}_{#1}}

\newcommand{\dx}{{\Delta x}}
\newcommand{\dy}{{\Delta y}}
\newcommand{\dz}{{\Delta z}}
\newcommand{\dxvec}{{\Delta \bm{x}}}

\newcommand{\ds}{\Delta s}

\newcommand{\surfarea}[1][]{\spoint{A}_{#1}}
\newcommand{\dxgridvec}{\Delta \bm{\grid{x}}}
\newcommand{\dxgridtensor}{\Delta \bm{\grid{X}}}

\newcommand{\normalbase}{n}
\newcommand{\normalvecc}[1]{\spoint{\normalbase}_{#1}}
\newcommand{\normalvec}[1][]{\vpoint{\normalbase}_{#1}}

\newcommand{\tanonebase}{t}
\newcommand{\tanonevecc}[1]{\spoint{\tanonebase}_{#1}}
\newcommand{\tanonevec}[1][]{\vpoint{\tanonebase}_{#1}}

\newcommand{\tantwobase}{b}
\newcommand{\tantwovecc}[1]{\spoint{\tantwobase}_{#1}}
\newcommand{\tantwovec}[1][]{\vpoint{\tantwobase}_{#1}}

\newcommand{\gspace}[1]{\mathcal{#1}}
\newcommand{\gspacevec}[1]{\boldsymbol{\mathcal{#1}}}
\newcommand{\faces}{\gspacevec{F}}
\newcommand{\facescomponent}{\gspace{F}}
\newcommand{\centers}{\gspace{C}}
\newcommand{\edges}{\gspacevec{E}}
\newcommand{\edgescomponent}{\gspace{E}}

\newcommand{\tensorspace}{\gspacevec{D}}
\newcommand{\spoints}{\gspace{S}}
\newcommand{\vpoints}{\gspacevec{V}}
\newcommand{\tpoints}{\gspacevec{T}}

\newcommand{\ddf}{\delta}

\newcommand{\dddf}{\grid{d}}

\newcommand{\hgrid}{\grid{H}}
\newcommand{\hgridcenters}[1][]{\hgrid_{\centers{#1}}}
\newcommand{\hgridfaces}[1][]{\bm{\hgrid}_{\faces{#1}}}
\newcommand{\hgridtensorspace}[1][]{\bm{\hgrid}_{\tensorspace{#1}}}

\newcommand{\had}{\circ}
\newcommand{\interpbody}{\grid{E}}
\newcommand{\spacetransform}[2]{{}^{#2}\grid{I}_{#1}}

\newcommand{\regds}{\grid{R}}
\newcommand{\dergrid}[1]{\grid{D}_{#1}}
\newcommand{\gradgrid}{\grid{G}}
\newcommand{\divgrid}{\grid{D}}
\newcommand{\curlgrid}{\grid{C}}
\newcommand{\lapgrid}{\grid{L}}
\newcommand{\convecgrid}{\grid{N}}
\newcommand{\invlapgrid}{\grid{L}^{-1}}
\newcommand{\outerproductoperator}{\grid{O}}
\newcommand{\schur}{\grid{S}}

\newcommand{\bcgrid}{\grid{b}}
\newcommand{\bcgridvec}{\bm{\grid{b}}}
\newcommand{\scalargrid}[1][]{\grid{s}_{#1}}
\newcommand{\vectorgrid}[1][]{\bm{\grid{v}}_{#1}}
\newcommand{\vectorgridcomponent}[1]{\grid{v}_{#1}}
\newcommand{\edgevectorgrid}[1][]{\bm{\grid{w}}_{#1}}
\newcommand{\edgevectorgridcomponent}[1]{\grid{w}_{#1}}
\newcommand{\tensorgrid}[1][]{\bm{\grid{T}}_{#1}}
\newcommand{\tensorgridcomponent}[1]{\grid{T}_{#1}}

\newcommand{\potgrid}{\grid{u}}
\newcommand{\sourcegrid}{\grid{q}}

\newcommand{\velgrid}{\bm{\grid{v}}}
\newcommand{\velgridcomponent}{\grid{v}}
\newcommand{\pgrid}{\grid{p}}

\newcommand{\forcegridvec}{\bm{\grid{F}}}
\newcommand{\forcegrid}{\grid{F}}

\newcommand{\zerosspoint}{\spoint{0}}

\newcommand{\forcepoint}{\spoint{f}}
\newcommand{\forcepointvec}{\vpoint{f}}
\newcommand{\vectorpoint}[1][]{\vpoint{v}_{#1}}
\newcommand{\vectorpointcomponent}[1][]{\spoint{v}_{#1}}
\newcommand{\tensorpoint}[1][]{\tpoint{T}_{#1}}
\newcommand{\tensorpointcomponent}[1][]{\spoint{T}_{#1}}

\newcommand{\scalarpoint}[1][]{\spoint{s}_{#1}}

\newcommand{\scalarpointsurfplus}[1][]{\spoint{s}^{+}_{\surf{#1}}}
\newcommand{\scalarpointsurfminus}[1][]{\spoint{s}^{-}_{\surf{#1}}}
\newcommand{\scalarpointsurfplusminus}[1][]{\spoint{s}^{\pm}_{\surf#1}}

\newcommand{\scalarpointnormal}[1][]{\spoint{s}^{\normalbase}_{#1}}
\newcommand{\scalarpointtanone}[1][]{\spoint{s}^{\tanonebase}_{#1}}
\newcommand{\scalarpointtantwo}[1][]{\spoint{s}^{\tantwobase}_{#1}}

\newcommand{\scalarpointsurfnormalplusminus}[1][]{\spoint{s}^{\normalbase \pm}_{\surf{#1}}}

\newcommand{\scalarpointsurftanoneplusminus}[1][]{\spoint{s}^{\tanonebase \pm}_{\surf{#1}}}

\newcommand{\scalarpointsurftantwoplusminus}[1][]{\spoint{s}^{\tantwobase \pm}_{\surf{#1}}}

\newcommand{\potpoint}[1][]{\spoint{u}_{#1}}
\newcommand{\potpointsurf}[1][]{\spoint{u}_{\surf{#1}}}
\newcommand{\potpointsurfplus}[1][]{\spoint{u}^{+}_{\surf{#1}}}
\newcommand{\potpointsurfminus}[1][]{\spoint{u}^{-}_{\surf{#1}}}
\newcommand{\potpointsurfplusminus}[1][]{\spoint{u}^{\pm}_{\surf{#1}}}

\newcommand{\potpointnormal}[1][]{\spoint{u}^{\normalbase}}
\newcommand{\potpointtanone}[1][]{\spoint{u}^{\tanonebase}}
\newcommand{\potpointtantwo}[1][]{\spoint{u}^{\tantwobase}}

\newcommand{\potpointsurfnormalplus}[1][]{ \spoint{u}^{\normalbase +}_{\surf{#1}}}
\newcommand{\potpointsurfnormalminus}[1][]{\spoint{u}^{\normalbase -}_{\surf{#1}}}

\newcommand{\velpoint}[1][]{\vpoint{v}_{#1}}
\newcommand{\velpointsurf}[1][]{\vpoint{v}_{\surf{#1}}}
\newcommand{\velpointsurfplus}[1][]{\vpoint{v}^{+}_{\surf{#1}}}
\newcommand{\velpointsurfminus}[1][]{\vpoint{v}^{-}_{\surf{#1}}}
\newcommand{\velpointsurfplusminus}[1][]{\vpoint{v}^{\pm}_{\surf{#1}}}

\newcommand{\velcomponentpointsurfplusminus}[1]{\spoint{v}^{\pm}_{#1,\surf}}

\newcommand{\velpointnormal}{\vpoint{v}^{\normalbase}}
\newcommand{\velpointtanone}{\vpoint{v}^{\tanonebase}}
\newcommand{\velpointtantwo}{\vpoint{v}^{\tantwobase}}

\newcommand{\velcomponentpointnormal}[1]{\spoint{v}^{\normalbase}_{#1}}

\newcommand{\velpointsurfnormalplusminus}{\vpoint{v}^{\normalbase \pm}_{\surf}}

\newcommand{\velpointsurftanoneplusminus}{\vpoint{v}^{\tanonebase \pm}_{\surf}}

\newcommand{\velpointsurftantwoplusminus}{\vpoint{v}^{\tantwobase \pm}_{\surf}}

\newcommand{\velcomponentpointsurfnormalplusminus}[1]{\spoint{v}^{\normalbase \pm}_{#1,\surf}}

\newcommand{\velcomponentpointsurftanoneplusminus}[1]{\spoint{v}^{\tanonebase \pm}_{#1,\surf}}

\newcommand{\velcomponentpointsurftantwoplusminus}[1]{\spoint{v}^{\tantwobase \pm}_{#1,\surf}}

\newcommand{\ppoint}[1][]{\spoint{p}_{#1}}

\newcommand{\ppointsurfplus}[1][]{\spoint{p}^{+}_{\surf{#1}}}
\newcommand{\ppointsurfminus}[1][]{\spoint{p}^{-}_{\surf{#1}}}
\newcommand{\ppointsurfplusminus}[1][]{\spoint{p}^{\pm}_{\surf{#1}}}

\newcommand{\Z}{\mathbb{Z}}

\newcommand{\diag}{\mathrm{diag}}

\newcommand{\basicI}{I}
\newcommand{\basicA}{A}
\newcommand{\basicBtwo}{B_{2}}
\newcommand{\basicBoneT}{B_{1}^\top}
\newcommand{\basicC}{C}
\newcommand{\basicS}{S}
\newcommand{\basicrone}{r_{1}}
\newcommand{\basicrtwo}{r_{2}}
\newcommand{\basicx}{x}
\newcommand{\basicy}{y}

\newcommand{\id}{\grid{I}}

%% file: colors.tex
\usepackage{xcolor}
\definecolor{myorange}{RGB}{230,159,0}
\definecolor{mylightblue}{RGB}{86,180,233}
\definecolor{mygreen}{RGB}{0,143,0}
\definecolor{myyellow}{RGB}{240,228,66}
\definecolor{mydarkblue}{RGB}{0,114,178}
\definecolor{myred}{RGB}{255,0,0}
\definecolor{mypurple}{RGB}{148,33,146}

%% file: figures/2D_surface_discretization.tikz
\begin{tikzpicture}

\def\xib{0.24052144192644503}
\def\yib{0.8214786128580098}

\begin{axis}[
    axis equal image,
    minor tick num=1,
    xtick=\empty,
    ytick=\empty,
    width=8cm,
    height=6cm,
    xmin=-0.5,
    xmax=1.7,
    ymin=-0.2,
    ymax=1.4,
    clip=false,
    after end axis/.code={%
        \begin{scope}[on background layer]
            \foreach \x in {-0.4,-0.2,...,1.6}{
                \foreach \y in {-0.1,0.1,...,1.5}{
                    \fill[lightgray] (axis cs:\x,\y) circle[radius=0.75pt];
                }
            }
        \end{scope}
    }]
    ]

    \addplot[
    thick, lightgray,
    ]
    table[x index=0, y index=1] {figures/data/2D_arbitrary_body/arbitrary_body_fine.txt};
    
    \addplot[
    only marks,
    mark=square*,
    mark size=1.2pt,
    color=red,
    ]
    table [x index=0, y index=1] {figures/data/2D_arbitrary_body/arbitrary_body_discretized.txt};
    
    
    \addplot[
    quiver={
      u=\thisrowno{2},
      v=\thisrowno{3},
      scale arrows=0.15, 
    },
    -{Latex[length=3pt]},
    color=black,
    restrict expr to domain={\coordindex}{1:3}
    ]
    table [x index=0, y index=1] {figures/data/2D_arbitrary_body/arbitrary_body_discretized.txt};

    \node[anchor=center,font=\small] at (axis cs:-0.3,1.2) {$\Omega$};
    \node[anchor=center,font=\small] at (axis cs:1.15,0.94) {$\normalvec[l+1]$};
    \node[anchor=center,font=\small] at (axis cs:1.35,0.67) {$\normalvec[l  ]$};
    \node[anchor=center,font=\small] at (axis cs:1.3,0.22)  {$\normalvec[l-1]$};
    \node[anchor=east,font=\small] at (axis cs:1.12,0.75) {$l+1$};
    \node[anchor=east,font=\small] at (axis cs:1.20,0.55) {$l$};
    \node[anchor=east,font=\small] at (axis cs:1.15,0.35) {$l-1$};

    
    
    \draw[->] (axis cs:-0.5,-0.3) -- (axis cs:1.8,-0.3);
    
    \draw[->] (axis cs:-0.6,-0.2) -- (axis cs:-0.6,1.5);
  
    \pgfplotsinvokeforeach{-0.4,-0.2,...,1.6}{%
        \draw[] (axis cs:#1,-0.32) -- (axis cs:#1,-0.28);
    }

    \pgfplotsinvokeforeach{-0.1,0.1,...,1.3}{%
        \draw[] (axis cs:-0.62,#1) -- (axis cs:-0.58,#1);
    }

    \node[anchor=north,font=\small] at (axis cs:0.0,-0.32) {$i-1$};    
    \node[anchor=north,font=\small] at (axis cs:0.2,-0.32) {$i$};    
    \node[anchor=north,font=\small] at (axis cs:0.4,-0.32) {$i+1$};    

    \node[anchor=east,font=\small] at (axis cs:-0.62,0.3) {$j-1$};    
    \node[anchor=east,font=\small] at (axis cs:-0.62,0.5) {$j$};    
    \node[anchor=east,font=\small] at (axis cs:-0.62,0.7) {$j+1$};    

    \addplot[only marks,mark=star,mark options={color=blue, scale=1}] coordinates {(0.2, 0.5)};

    \coordinate (ANCHOR1) at (axis cs: \xib, \yib);
    \coordinate (CTRL11)  at (axis cs: \xib - 0.05, \yib + 0.05);
    \coordinate (CTRL12)  at (axis cs: \xib - 0.1, \yib + 0.05);
    \coordinate (LABEL1)  at (axis cs: \xib - 0.1, \yib + 0.13);
    \draw [] (ANCHOR1) .. controls (CTRL11)  and (CTRL12) .. (LABEL1) node[above] {$\surf$};
\end{axis}

\end{tikzpicture}

%% file: figures/3D_grid_cell.tikz
\tdplotsetmaincoords{75}{30}
\begin{tikzpicture}[scale=1.6,tdplot_main_coords,
    cubeedge/.style={thick, line join=round},
    rib/.style={dotted,very thick},
    faceLabel/.style={font=\small,rounded corners=2pt,inner sep=1pt},
    edgeLabel/.style={font=\small,rounded corners=2pt,inner sep=1pt}
  ]

\def\a{2}

\coordinate (O)     at (0,0,0);
\coordinate (X)     at (\a,0,0);
\coordinate (Y)     at (0,\a,0);
\coordinate (Z)     at (0,0,\a);
\coordinate (XY)    at (\a,\a,0);
\coordinate (XZ)    at (\a,0,\a);
\coordinate (YZ)    at (0,\a,\a);
\coordinate (XYZ)   at (\a,\a,\a);

\draw[cubeedge] (O) -- (X) -- (XY) -- (Y) -- cycle;   
\draw[cubeedge] (Z) -- (XZ) -- (XYZ) -- (YZ) -- cycle; 
\draw[cubeedge] (O) -- (Z);
\draw[cubeedge] (X) -- (XZ);
\draw[cubeedge] (Y) -- (YZ);
\draw[cubeedge] (XY) -- (XYZ);

\coordinate (Ffront) at ({\a/2},0,{\a/2});   
\coordinate (Fback)  at ({\a/2},\a,{\a/2});  
\coordinate (FsideL) at (0,{\a/2},{\a/2});   
\coordinate (FsideR) at (\a,{\a/2},{\a/2});  
\coordinate (Ftop)   at ({\a/2},{\a/2},\a);  
\coordinate (Fbottom)at ({\a/2},{\a/2},0);   

\coordinate (Ex) at ({\a/2},0,\a);   
\coordinate (Ey)  at (\a,0,{\a/2});  
\coordinate (Ez) at (\a,{\a/2},\a);   

\coordinate (C) at ({\a/2},{\a/2},{\a/2});
\fill (C) circle (1pt) node[above right=2pt] {\small $\centers$};

\draw[rib] (Ffront) -- (Fback);
\draw[rib] (FsideL) -- (FsideR);
\draw[rib] (Fbottom) -- (Ftop);

\fill (Ffront)  circle (1pt) node[left=2pt,faceLabel] {\small  $\facescomponent_z$};
\fill (FsideR)  circle (1pt) node[right=2pt,faceLabel] {\small  $\facescomponent_x$};
\fill (Ftop)    circle (1pt) node[above=2pt,faceLabel] {\small  $\facescomponent_y$};

\fill (Ex)  circle (1pt) node[above=2pt,edgeLabel] {\small  $\edgescomponent_x$};
\fill (Ey)  circle (1pt) node[right=2pt,edgeLabel] {\small  $\edgescomponent_y$};
\fill (Ez)    circle (1pt) node[above=2pt,edgeLabel] {\small  $\edgescomponent_z$};

\begin{scope}[shift={(-0.6,-0.6,-0.3)}]
  \draw[->,thick] (0,0,0) -- (1,0,0) node[anchor=north east] {$x$};
  \draw[->,thick] (0,0,0) -- (0,-1,0) node[anchor=east] {$z$};
  \draw[->,thick] (0,0,0) -- (0,0,1) node[anchor=south] {$y$};
\end{scope}

\end{tikzpicture}

%% file: figures/1D_example_diagram.tikz
\begin{tikzpicture}[
    >={Stealth[width=5pt]}, 
    every node/.style={font=\small}
]
\pgfplotsset{set layers, cell picture=true}

\def\axWsmall{4cm}
\def\axHsmall{4cm}
\def\axW{4.7cm}
\def\axH{4.7cm}
\def\xsep{3.25cm}
\def\ysep{4.5cm}





\coordinate (R3leftmost)  at (1.3cm, -2*\ysep);
\coordinate (R3rightmost) at ( \linewidth, -2*\ysep);




\begin{axis}[
    name=A3left,
    at={(R3leftmost)}, anchor=left of west,
    width=\axW, height=0.75*\axH,
    axis lines=box,
    ymin=-0.02,
    ymax=0.09,
    scaled ticks=false,
    ytick={-0.05,0,0.05},
    yticklabels={$-0.05$, $0$, $0.05$},
    title={},
    xlabel={$x$}, ylabel={},
    clip mode=individual,
]
\addplot[color=mygreen, mark=star, mark options={scale=1}] table[x index=0, y index=1] {figures/data/1D_example_nx_17_ddf_Yang3_alpha_1/u_ef_RTm.txt};
\end{axis}

\begin{axis}[
    name=A3right,
    at={(R3rightmost)}, anchor=right of east,
    width=\axW, height=0.75*\axH,
    axis lines=box,
    ymin=-0.02,
    ymax=0.09,
    scaled ticks=false,
    ytick={-0.05,0,0.05},
    yticklabels=\empty,
    xlabel={$x$}, ylabel={},
    title={},
    clip mode=individual
]
\addplot[color=myorange, mark=o, mark options={scale=0.8}] table[x index=0, y index=1] {figures/data/1D_example_nx_17_ddf_Yang3_alpha_1/f_mod_H_RTm1.txt};
\end{axis}

\coordinate (R3centerleft)  at ($(A3left.right of east) + (0.1cm,0)$);
\coordinate (R3centerright) at ($(A3right.left of west) - (0.4cm,0)$);

\begin{axis}[
    name=A3centerleft,
    at={(R3centerleft)}, anchor=left of west,
    width=\axW, height=0.75*\axH,
    axis lines=box,
    ymin=-0.02,
    ymax=0.09,
    scaled ticks=false,
    ytick={-0.05,0,0.05},
    yticklabels=\empty,
    xlabel={$x$}, ylabel={},
    title={},
    clip mode=individual,
]
\addplot[color=red, mark=diamond, mark options={scale=1.0}] table[x index=0, y index=1] {figures/data/1D_example_nx_17_ddf_Yang3_alpha_1/u_IBM_RTm.txt};
\end{axis}

\begin{axis}[
    name=A3centerright,
    at={(R3centerright)}, anchor=right of east,
    width=\axW, height=0.75*\axH,
    axis lines=box,
    ymin=-0.02,
    ymax=0.09,
    scaled ticks=false,
    ytick={-0.05,0,0.05},
    yticklabels=\empty,
    xlabel={$x$}, ylabel={},
    title={},
    clip mode=individual
]
\addplot[color=blue, mark=triangle, mark options={scale=1.0}] table[x index=0, y index=1] {figures/data/1D_example_nx_17_ddf_Yang3_alpha_1/u_mod_RTm.txt};
\end{axis}

\coordinate (R1left)   at ($0.5*(A3left.center) + 0.5*(A3centerleft.center) + (0,1.85*\ysep)$);
\coordinate (R1centerleft)   at ($(A3centerright.center) + (-0.1*\xsep,1.85*\ysep)$);
\coordinate (R1centerright)  at ($(A3centerright.center) + (0.7*\xsep,1.85*\ysep)$);

\coordinate (R0left)          at ($0.5*(A3left.center) + 0.5*(A3centerleft.center) + (0,2.7*\ysep)$);
\coordinate (R0centerright)   at ($(A3centerright.center) + (0,2.7*\ysep)$);
\coordinate (R0right)         at ($(A3right.center) + (0,2.7*\ysep)$);

\coordinate (R2left)   at ($(A3left.center) + (0,0.85*\ysep)$);
\coordinate (R2centerleft)   at ($(A3centerleft.center) + (0,0.85*\ysep)$);
\coordinate (R2centerright)   at ($(A3centerright.center) + (0,0.85*\ysep)$);
\coordinate (R2right)  at ($(A3right.center) + (0,0.85*\ysep)$);

\coordinate (R4left)   at ($(A3left.below south) + (0,-1cm)$);
\coordinate (R4center) at ($(A3centerleft.below south) + (0,-1cm)$);
\coordinate (R4right)  at ($0.5*(A3centerright.below south) + 0.5*(A3right.below south) + (0,-1.2cm)$);

\begin{axis}[
    name=A0left,
    at={(R0left)}, anchor=center,
    width=\axWsmall, height=\axHsmall,
    axis lines=box,
    ymin=-0.5,
    ymax=5.5,
    xlabel={$x$}, ylabel={},
    ylabel shift=-0.1cm,
]
\addplot[color=black, mark=square, mark options={scale=0.8}] table[x index=0, y expr=\thisrowno{1}/4] {figures/data/1D_example_nx_17_ddf_Yang3_alpha_1/f_ef_Rm.txt};
\node[draw=black,fill=white,anchor=north west] at (rel axis cs:0,1) {$\regds_{\centers}$};
\end{axis}

\begin{axis}[
    name=A0centerright,
    at={(R0centerright)}, anchor=center,
    width=\axWsmall, height=\axHsmall,
    axis lines=box,
    ymin=-0.5,
    ymax=5.5,
    xlabel={$x$}, ylabel={},
    ylabel shift=-0.1cm,
]
\addplot[color=black, mark=square, mark options={scale=0.8}] table[x index=0, y index=1] {figures/data/1D_example_nx_17_ddf_Yang3_alpha_1/R.txt};
\node[draw=black,fill=white,anchor=north west] at (rel axis cs:0,1) {$\regds_{\faces}$};
\end{axis}

\begin{axis}[
    name=A0right,
    at={(R0right)}, anchor=center,
    width=\axWsmall, height=\axHsmall,
    axis lines=box,
    xlabel={$x$}, ylabel={},
    ylabel shift=-0.5cm,
]
\addplot[color=myorange, mark=o, mark options={scale=0.8}] table[x index=0, y index=1] {figures/data/1D_example_nx_17_ddf_Yang3_alpha_1/DR.txt};
\node[draw=black,fill=white,anchor=north west] at (rel axis cs:0,1) {$\divgrid \regds_{\faces}$};
\end{axis}

\begin{axis}[
    name=A1left,
    at={(R1left)}, anchor=center,
    width=\axWsmall, height=\axHsmall,
    axis lines=box,
    xlabel={$x$}, ylabel={},
    ylabel shift=-0.1cm,
    ymin=-8,
    ymax=25,
]
\addplot[color=mygreen, mark=star, mark options={scale=1}] table[x index=0, y index=1] {figures/data/1D_example_nx_17_ddf_Yang3_alpha_1/f_ef_Rm.txt};
\addplot[color=red, mark=diamond, mark options={scale=1.0}] table[x index=0, y index=1] {figures/data/1D_example_nx_17_ddf_Yang3_alpha_1/f_IBM_Rm.txt};
\node[draw=black,fill=white,anchor=north west] at (rel axis cs:0,1) {$\regds_{\centers} \forcepoint $};
\end{axis}

\begin{axis}[
    name=A1centerleft,
    at={(R1centerleft)}, anchor=center,
    width=\axWsmall, height=\axHsmall,
    yticklabels=\empty,
    axis lines=box,
    ymin=-8,
    ymax=25,
    xlabel={$x$}, ylabel={},
]
\addplot[color=blue, mark=triangle, mark options={scale=1}] table[x index=0, y index=1] {figures/data/1D_example_nx_17_ddf_Yang3_alpha_1/f_mod_IR.txt};
\node[draw=black,fill=white,anchor=north west] at (rel axis cs:0,1) {$\spacetransform{\faces}{\centers}\regds_{\faces} \jump{\potpointnormal} $};
\end{axis}

\begin{axis}[
    name=A1centerright,
    at={(R1centerright)}, anchor=center,
    width=\axWsmall, height=\axHsmall,
    axis lines=box,
    yticklabel pos=right,
    xlabel={$x$}, ylabel={},
    ylabel shift=-0.1cm,
    ymin=-8,
    ymax=25,
    axis background/.style={fill=white}
]
\addplot[color=blue, mark=triangle, mark options={scale=1}] table[x index=0, y index=1] {figures/data/1D_example_nx_17_ddf_Yang3_alpha_1/f_mod_DR1.txt};
\node[draw=black,fill=white,anchor=north west] at (rel axis cs:0,1) {$\divgrid \regds_{\faces,1n}\jump{\potpointnormal} $};
\end{axis}

\begin{axis}[
    name=A2left,
    at={(R2left)}, anchor=center,
    width=\axW, height=\axH,
    axis lines=box,
    ytick={0,0.5},
    yticklabels={$0$, $0.5$},
    xlabel={}, ylabel={},
    xticklabels=\empty,
    title={},
    title style = {align = center, font=\small},
    ymin=-0.2,
    ymax=0.65,
    clip mode=individual
]
\addplot[color=lightgray, line join=bevel] table[x index=0, y index=1] {figures/data/1D_example_nx_17_ddf_Yang3_alpha_1/exact_solution.txt};
\addplot[color=lightgray, only marks, mark=*, mark options={scale=0.8}] table[x index=0, y index=1] {figures/data/1D_example_nx_17_ddf_Yang3_alpha_1/u_exact.txt};
\addplot[color=mygreen, mark=star, mark options={scale=1}] table[x index=0, y index=1] {figures/data/1D_example_nx_17_ddf_Yang3_alpha_1/u_ef.txt};

\pgfplotsextra{
  \draw[dashed,thick,black] (axis cs:0.7,0.09) rectangle (axis cs:1.3,0.35);
  \node[anchor=north] (lblA) at (axis cs:0.7,0.0) {\scriptsize $\epsilon_i \sim \mathcal{O}(\dx)$};
  \node[anchor=south] at (axis cs:1.5,0.5) {\scriptsize $\epsilon_i \sim \mathcal{O}(\dx^2)$};
  \draw[thick] (axis cs:0.85,0.09) -- (lblA.north);
}
\node[draw=black,fill=white,anchor=north west] at (rel axis cs:0,1) {$\mask{\potgrid}_\mathrm{AF}$};
\end{axis}

\begin{axis}[
    name=A2centerleft,
    at={(R2centerleft)}, anchor=center,
    width=\axW, height=\axH,
    ytick={0,0.5},
    yticklabels=\empty,
    xticklabels=\empty,
    xlabel={}, ylabel={},
    title={},
    title style = {align = center, font=\small},
    ymin=-0.2,
    ymax=0.65,
    clip mode=individual
]
\addplot[color=lightgray, line join=bevel] table[x index=0, y index=1] {figures/data/1D_example_nx_17_ddf_Yang3_alpha_1/exact_solution.txt};
\addplot[color=lightgray, only marks, mark=*, mark options={scale=0.8}] table[x index=0, y index=1] {figures/data/1D_example_nx_17_ddf_Yang3_alpha_1/u_exact.txt};
\addplot[color=red, mark=diamond, mark options={scale=1.0}] table[x index=0, y index=1] {figures/data/1D_example_nx_17_ddf_Yang3_alpha_1/u_IBM.txt};

\pgfplotsextra{
  \node[anchor=south] at (axis cs:1.5,0.5) {\scriptsize $\epsilon_i \sim \mathcal{O}(\dx)$};
}
\node[draw=black,fill=white,anchor=north west] at (rel axis cs:0,1) {$\mask{\potgrid}_\mathrm{IB1}$};
\end{axis}

\begin{axis}[
    name=A2centerright,
    at={(R2centerright)}, anchor=center,
    width=\axW, height=\axH,
    axis lines=box,
    ytick={0,0.5},
    yticklabels=\empty,
    xlabel={}, ylabel={},
    xticklabels=\empty,
    title={},
    title style = {align = center, font=\small},
    ymin=-0.2,
    ymax=0.65,
    clip mode=individual
]
\addplot[color=lightgray, line join=bevel] table[x index=0, y index=1] {figures/data/1D_example_nx_17_ddf_Yang3_alpha_1/exact_solution.txt};
\addplot[color=lightgray, only marks, mark=*, mark options={scale=0.8}] table[x index=0, y index=1] {figures/data/1D_example_nx_17_ddf_Yang3_alpha_1/u_exact.txt};
\addplot[color=blue, mark=triangle, mark options={scale=1}] table[x index=0, y index=1] {figures/data/1D_example_nx_17_ddf_Yang3_alpha_1/u_mod.txt};

\pgfplotsextra{
  \draw[dashed,thick,black] (axis cs:0.7,0.06) rectangle (axis cs:1.3,0.35);
  \node[anchor=north] (lblA) at (axis cs:0.7,0.0) {\scriptsize $\epsilon_i \sim \mathcal{O}(\dx)$};
  \node[anchor=south] at (axis cs:1.5,0.5) {\scriptsize $\epsilon_i \sim \mathcal{O}(\dx^2)$};
  \draw[thick] (axis cs:0.85,0.06) -- (lblA.north);
}
\node[draw=black,fill=white,anchor=north west] at (rel axis cs:0,1) {$\mask{\potgrid}_\mathrm{IB2}$};
\end{axis}

\begin{axis}[
    name=A2right,
    at={(R2right)}, anchor=center,
    width=\axW, height=\axH,
    axis lines=box,
    ytick={0,0.5,1.0},
    xlabel={}, ylabel={},
    xticklabels=\empty,
    title={},
    title style = {align = center, font=\small},
    ymin=-0.2,
    ymax=1.2,
    clip mode=individual
]
\addplot[color=myorange, mark=o, mark options={scale=0.8}] table[x index=0, y index=1] {figures/data/1D_example_nx_17_ddf_Yang3_alpha_1/HplusC.txt};
\node[draw=black,fill=white,anchor=north west] at (rel axis cs:0,1) {$\hgridcenters^+$};
\end{axis}

\node[anchor=center] (Lblleft) at (R4left)    {\small $\displaystyle\sum_i \potgrid_{(i)} \ddf_\dx(\xgrid[\centers,(i)] - x_\surf) \neq \potpointsurf$};
\node[anchor=center] (Lblcenterleft) at (R4center)   {\small $\displaystyle\sum_i \potgrid_{(i)} \ddf_\dx(\xgrid[\centers,(i)] - x_\surf) = \potpointsurf$};
\node[anchor=center, text width=0.35*\linewidth] (Lblright) at (R4right) {$\displaystyle\sum_i  \Big(  {\color{blue}{\potgrid_{(i)} \ddf_\dx(\xgrid[\centers,(i)] - x_\surf)}}$ \\ \hfill $- {\color{myorange}{\hgridcenters^+ (\xgrid[\centers,(i)] - x_\surf) \ddf_\dx(\xgrid[\centers,(i)] - x_\surf) \jump{\potpointnormal}}} \Big)  = \potpointsurf $};


\draw[->] (A0left.below south) -- (A1left.north) node[midway, above, align=center] {};

\draw[->] (A0centerright.below south) |- ($(A0centerright.below south)!0.4!(A1centerleft.north)$) -| (A1centerleft.north) node[midway, above, align=center] {};
\draw[->] (A0centerright.below south) |- ($(A0centerright.below south)!0.4!(A1centerright.north)$) -| (A1centerright.north) node[pos=1, above left, align=center] {};
\draw[->] (A0centerright.east) -- (A0right.left of west) node[midway, above, align=center] {};

\draw[->] (A1left.below south) |- ($(A1left.below south)!0.5!(A2left.north)$) -| (A2left.north) node[midway, above, align=center] (forcinglabel) {analytical forcing:\\$\forcepoint = \jump{\pot^\normalbase} = \jump{\inlder{\pot}{x}{}}$};
\draw[->] (A1left.below south) |- ($(A1left.below south)!0.5!(A2centerleft.north)$) -| (A2centerleft.north) node[midway, above, align=center] {forcing from Eq.~(\ref{eq:poisson_old_force_equation}):\\$\forcepoint = \jump{\potpointnormal}$};
\draw[->] (A1centerleft.below south) |- ($(A1centerleft.below south)!0.5!(A2centerright.north)$) -| (A2centerright.north) node[midway, above, align=center] {};
\draw[->] (A1centerright.below south) |- ($(A1centerright.below south)!0.5!(A2centerright.north)$) -| (A2centerright.north) node[midway, above, align=center] {};

\coordinate (A0rightbottomanchor) at ($(A0right.south west)!0.8!(A0right.south east)$);
\draw[->] (A0rightbottomanchor) -- (A0rightbottomanchor |- A2right.north);

\draw[->] (A2left.south) -- (A3left.north) node[midway,left, align=center] (interpolationlabel) {multiply by \\$\ddf_\dx (\xgrid[\centers,(i)] - x_\surf)$};
\draw[->] (A2centerleft.south) -- (A3centerleft.north) node[midway,left, align=center] {multiply by \\$\ddf_\dx (\xgrid[\centers,(i)] - x_\surf)$};
\draw[->] (A2centerright.south) -- (A3centerright.north) node[midway,left, align=center] {multiply by \\$\ddf_\dx (\xgrid[\centers,(i)] - x_\surf)$};
\draw[->] ($(A2right.south west)!0.8!(A2right.south east)$) -- ($(A3right.north west)!0.8!(A3right.north east)$) node[midway,left, align=center] {multiply by \\$(\xgrid[\centers,(i)] - x_\surf)\ddf_\dx (\xgrid[\centers,(i)] - x_\surf)  \jump{\potpointnormal} $};

\draw[->] (A3left.below south) -- (Lblleft.north);
\draw[->] (A3centerleft.below south) -- (Lblcenterleft.north);
\draw[->] (A3centerright.below south) |- ($(A3centerright.below south)!0.3!(Lblright.north)$) -| (Lblright.north);
\draw[->] (A3right.below south) |- ($(A3right.below south)!0.3!(Lblright.north)$) -| (Lblright.north);

\begin{pgfonlayer}{axis grid}
    \draw[gray,very thick,dashed] ($(A3centerleft.south east)!0.5!(A3centerright.south west) - (0,2cm)$) -- ($(A3centerleft.south east)!0.5!(A3centerright.south west) + (0,3.4*\ysep)$);
\end{pgfonlayer}

\coordinate (origin) at (current bounding box.west);

\coordinate (tmp) at (A0left.north);
\coordinate (topforcing) at (origin |- tmp);
\coordinate (tmp) at (forcinglabel.south);
\coordinate (bottomforcing) at (origin |- tmp);
\draw[decorate, decoration={brace, mirror}](topforcing) -- (bottomforcing) node[midway, xshift=-0.5cm, rotate=90] {IB forcing};

\coordinate (tmp) at (A2left.north);
\coordinate (topsolution) at (origin |- tmp);
\coordinate (tmp) at (A2left.south);
\coordinate (bottomsolution) at (origin |- tmp);
\draw[decorate, decoration={brace, mirror}](topsolution) -- (bottomsolution) node[midway, xshift=-0.5cm, rotate=90] {Solution};

\coordinate (tmp) at (interpolationlabel.north);
\coordinate (topinterpolation) at (origin |- tmp);
\coordinate (tmp) at (current bounding box.south);
\coordinate (bottominterpolation) at (origin |- tmp);
\draw[decorate, decoration={brace, mirror}](topinterpolation) -- (bottominterpolation) node[midway, xshift=-0.5cm, rotate=90] {Interpolation of the solution to the IB};

\node[anchor=north, align=center] () at ($(A3left.south east)!0.5!(A3centerleft.south west) + (0,3.4*\ysep)$)   {Prototypical continuous-forcing IB method};
\node[anchor=north, align=center] () at ($(A3centerright.south east)!0.5!(A3right.south west) + (0,3.4*\ysep)$) {Proposed IB method};
\end{tikzpicture}

%% file: figures/1D_poisson_error_Linf_norm.tikz
\begin{tikzpicture}
    \begin{loglogaxis}[
        xlabel=$\dx$,
        ylabel={$||\mask{\potgrid} - \pot(\xgrid[\centers])||_\infty/||\pot(\xgrid[\centers])||_\infty$},
        width=5cm,
        height=5.5cm
    ]

        \addplot[very thick, gray, solid] table[x index=0, y expr=\thisrowno{0}*5               ] {figures/data/1D_poisson_error/u_ef_Linf_error.txt};
        \addplot[very thick, gray, solid] table[x index=0, y expr=\thisrowno{0}*\thisrowno{0}*0.5] {figures/data/1D_poisson_error/u_ef_Linf_error.txt};    
        
        \addplot[color=mygreen, mark=star, mark options={scale=1}]  table[x index=0, y index=1] {figures/data/1D_poisson_error/u_ef_Linf_error.txt};
        \addplot[color=red, mark=diamond, mark options={scale=1}]   table[x index=0, y index=1] {figures/data/1D_poisson_error/u_IBM_Linf_error.txt};
        \addplot[color=blue, mark=triangle, mark options={scale=1}] table[x index=0, y index=1] {figures/data/1D_poisson_error/u_mod_Linf_error.txt};
        
        \addplot[color=mygreen, mark=star, mark options={scale=1}, dashed, mark options={solid}]  table[x index=0, y index=1] {figures/data/1D_poisson_error/u_ef_nobody_Linf_error.txt};
        \addplot[color=red, mark=diamond, mark options={scale=1}, dashed, mark options={solid}]   table[x index=0, y index=1] {figures/data/1D_poisson_error/u_IBM_nobody_Linf_error.txt};
        \addplot[color=blue, mark=triangle, mark options={scale=1}, dashed, mark options={solid}] table[x index=0, y index=1] {figures/data/1D_poisson_error/u_mod_nobody_Linf_error.txt};
        
        \node[anchor=south west, text=gray] at (axis cs:0.008,0.1) {$\sim \dx$};
        \node[anchor=north west, text=gray] at (axis cs:0.008,0.00008) {$\sim \dx^2$};
    \end{loglogaxis}

\end{tikzpicture}

%% file: figures/1D_poisson_error_L2_norm.tikz
\begin{tikzpicture}
    \begin{loglogaxis}[
        xlabel=$\dx$,
        ylabel={$||\mask{\potgrid} - \pot(\xgrid[\centers])||_2/||\pot(\xgrid[\centers])||_2$},
        width=5cm,
        height=5.5cm
    ]

        \addplot[very thick, gray, solid] table[x index=0, y expr=\thisrowno{0}*5               ] {figures/data/1D_poisson_error/u_ef_L2_error.txt};
        \addplot[very thick, gray, solid] table[x index=0, y expr=\thisrowno{0}*\thisrowno{0}*0.5] {figures/data/1D_poisson_error/u_ef_L2_error.txt};    
        
        \addplot[color=mygreen, mark=star, mark options={scale=1}]  table[x index=0, y index=1] {figures/data/1D_poisson_error/u_ef_L2_error.txt};
        \addplot[color=red, mark=diamond, mark options={scale=1}]   table[x index=0, y index=1] {figures/data/1D_poisson_error/u_IBM_L2_error.txt};
        \addplot[color=blue, mark=triangle, mark options={scale=1}] table[x index=0, y index=1] {figures/data/1D_poisson_error/u_mod_L2_error.txt};
        
        \addplot[color=mygreen, mark=star, mark options={scale=1}, dashed, mark options={solid}]  table[x index=0, y index=1] {figures/data/1D_poisson_error/u_ef_nobody_L2_error.txt};
        \addplot[color=red, mark=diamond, mark options={scale=1}, dashed, mark options={solid}]   table[x index=0, y index=1] {figures/data/1D_poisson_error/u_IBM_nobody_L2_error.txt};
        \addplot[color=blue, mark=triangle, mark options={scale=1}, dashed, mark options={solid}] table[x index=0, y index=1] {figures/data/1D_poisson_error/u_mod_nobody_L2_error.txt};
        
        \node[anchor=south west, text=gray] at (axis cs:0.008,0.1) {$\sim \dx$};
        \node[anchor=north west, text=gray] at (axis cs:0.008,0.00008) {$\sim \dx^2$};
    \end{loglogaxis}

\end{tikzpicture}

%% file: figures/1D_poisson_forcing_error.tikz
\begin{tikzpicture}
    \begin{loglogaxis}[
        xlabel=$\dx$,
        ylabel={$|\jump{\potpoint[\normalbase]} - \jump{\inlder{\pot}{x}{}}|/|\jump{\inlder{\pot}{x}{}}|$},
        width=5cm,
        height=5.5cm
    ]

        \addplot[very thick, gray, solid] table[x index=0, y expr=\thisrowno{0}*8               ] {figures/data/1D_poisson_error/u_ef_L2_error.txt};
        \addplot[very thick, gray, solid] table[x index=0, y expr=\thisrowno{0}*\thisrowno{0}*0.5] {figures/data/1D_poisson_error/u_ef_L2_error.txt};    
        
        \addplot[color=red, mark=diamond, mark options={scale=1}]   table[x index=0, y index=1] {figures/data/1D_poisson_error/f_IBM_error.txt};
        \addplot[color=blue, mark=triangle, mark options={scale=1}] table[x index=0, y index=1] {figures/data/1D_poisson_error/f_mod_error.txt};
        
        \node[anchor=south west, text=gray] at (axis cs:0.008,0.2) {$\sim \dx$};
        \node[anchor=north west, text=gray] at (axis cs:0.008,0.00008) {$\sim \dx^2$};
    \end{loglogaxis}

\end{tikzpicture}

%% file: figures/2D_poisson_mesh.tikz
\begin{tikzpicture}
  \begin{axis}[
    view={-20}{30},
    xlabel=$x / R$,
    ylabel=$y / R$,
    zlabel={$\pot(\xgridvec[\centers])$},
    height=6cm,
    width=6cm]

      
    \addplot3[
      surf,
      shader=faceted,
      line join=round,
      faceted color=black,             
      colormap={graywhite}{gray(0cm)=(0.4); gray(1cm)=(1)}
    ] table [x index=0, y index=1, z index=2] {figures/data/2D_poisson/poisson_exact_solution_dx_0.2_dsdx_1.txt};

  \end{axis}
\end{tikzpicture}

%% file: figures/2D_poisson_solution_original.tikz
\begin{tikzpicture}
    \begin{axis}[
        name=mainplot,
        xmin=0,
        xmax=2.1,
        ymin=-0.1,
        ymax=1.1,
        xlabel=$x/R$,
        ylabel={$\mask{\potgrid}(\xgrid[\centers],y=0)$},
        width=5.75cm,
        height=6cm,
    ]

    \addplot[line width=1, black, dashed] table[x index=0, y index=1] {figures/data/2D_poisson/poisson_exact_centerline_dx_0.0125_dsdx_1.0_ddf_Yang3.txt};
    
    \addplot[red!100!white, line width=1.0,mark=diamond*,mark options={scale=0.50}] table[x index=0,y index=1] {figures/data/2D_poisson/poisson_original_method_centerline_dx_0.2_dsdx_1.0_ddf_Yang3.txt};
    \addplot[red!060!white, line width=0.75,mark=diamond*,mark options={scale=0.40}] table[x index=0,y index=1] {figures/data/2D_poisson/poisson_original_method_centerline_dx_0.1_dsdx_1.0_ddf_Yang3.txt};
    \addplot[red!030!white, line width=0.50,mark=diamond*,mark options={scale=0.30}] table[x index=0,y index=1] {figures/data/2D_poisson/poisson_original_method_centerline_dx_0.05_dsdx_1.0_ddf_Yang3.txt};

    \draw[thick] (axis cs:0.88,0.88) coordinate (zoomSW) rectangle (axis cs:1.12,1.07) coordinate (zoomNE);
    \coordinate (zoomSouth) at ($(zoomSW)!0.5!(zoomNE |- zoomSW)$);
    \end{axis}

    \begin{axis}[
        name=inset,
        width=3.7cm,
        height=3.7cm,
        xmin=0.88, xmax=1.12, 
        ymin=0.88, ymax=1.07, 
        at={(mainplot.south east)},
        xshift=-7mm, yshift=2mm,
        anchor=south east,
        xtick=\empty,
        ytick=\empty,
        clip marker paths=true
    ]
    \addplot[line width=1, black, dashed] table[x index=0, y index=1] {figures/data/2D_poisson/poisson_exact_centerline_dx_0.0125_dsdx_1.0_ddf_Yang3.txt};
    
    \addplot[red!100!white, line width=1.0,mark=diamond*,mark options={scale=0.60}] table[x index=0,y index=1] {figures/data/2D_poisson/poisson_original_method_centerline_dx_0.2_dsdx_1.0_ddf_Yang3.txt};
    \addplot[red!060!white, line width=0.75,mark=diamond*,mark options={scale=0.50}] table[x index=0,y index=1] {figures/data/2D_poisson/poisson_original_method_centerline_dx_0.1_dsdx_1.0_ddf_Yang3.txt};
    \addplot[red!030!white, line width=0.50,mark=diamond*,mark options={scale=0.40}] table[x index=0,y index=1] {figures/data/2D_poisson/poisson_original_method_centerline_dx_0.05_dsdx_1.0_ddf_Yang3.txt};

    \end{axis}
    
    \draw[thick] (zoomSouth) -- (inset.north);
\end{tikzpicture}

%% file: figures/2D_poisson_solution_new.tikz
\begin{tikzpicture}
    \begin{axis}[
        name=mainplot,
        xmin=0,
        xmax=2.1,
        ymin=-0.1,
        ymax=1.1,
        xlabel=$x/R$,
        width=5.75cm,
        height=6cm,
        yticklabels=\empty,
    ]

    \addplot[line width=1, black, dashed] table[x index=0, y index=1] {figures/data/2D_poisson/poisson_exact_centerline_dx_0.0125_dsdx_1.0_ddf_Yang3.txt};

    \addplot[blue!100!white,line width=1,mark=triangle*,mark options={scale=0.50}] table[x index=0,y index=1] {figures/data/2D_poisson/poisson_new_method_centerline_dx_0.2_dsdx_1.0_ddf_Yang3.txt};
    \addplot[blue!060!white,line width=0.75,mark=triangle*,mark options={scale=0.40}] table[x index=0,y index=1] {figures/data/2D_poisson/poisson_new_method_centerline_dx_0.1_dsdx_1.0_ddf_Yang3.txt};
    \addplot[blue!030!white,line width=0.5,mark=triangle*,mark options={scale=0.30}] table[x index=0,y index=1] {figures/data/2D_poisson/poisson_new_method_centerline_dx_0.05_dsdx_1.0_ddf_Yang3.txt};

    \draw[thick] (axis cs:0.88,0.88) coordinate (zoomSW) rectangle (axis cs:1.12,1.07) coordinate (zoomNE);
    \coordinate (zoomSouth) at ($(zoomSW)!0.5!(zoomNE |- zoomSW)$);
    \end{axis}

    \begin{axis}[
        name=inset,
        width=3.7cm,
        height=3.7cm,
        xmin=0.88, xmax=1.12, 
        ymin=0.88, ymax=1.07, 
        at={(mainplot.south east)},
        xshift=-7mm, yshift=2mm,
        anchor=south east,
        xtick=\empty,
        ytick=\empty,
        clip marker paths=true
    ]

    \addplot[line width=1, black, dashed] table[x index=0, y index=1] {figures/data/2D_poisson/poisson_exact_centerline_dx_0.0125_dsdx_1.0_ddf_Yang3.txt};

    \addplot[blue!100!white, line width=1,mark=triangle*,mark options={scale=0.60}] table[x index=0,y index=1] {figures/data/2D_poisson/poisson_new_method_centerline_dx_0.2_dsdx_1.0_ddf_Yang3.txt};
    \addplot[blue!060!white,line width=0.75,mark=triangle*,mark options={scale=0.50}] table[x index=0,y index=1] {figures/data/2D_poisson/poisson_new_method_centerline_dx_0.1_dsdx_1.0_ddf_Yang3.txt};
    \addplot[blue!030!white, line width=0.5,mark=triangle*,mark options={scale=0.30}] table[x index=0,y index=1] {figures/data/2D_poisson/poisson_new_method_centerline_dx_0.05_dsdx_1.0_ddf_Yang3.txt};

    \end{axis}
    \draw[thick] (zoomSouth) -- (inset.north);

\end{tikzpicture}

%% file: figures/2D_poisson_error_Linf_norm_dsdx_variation.tikz
\begin{tikzpicture}
    \begin{loglogaxis}[
        xlabel=$\dx/R$,
        ylabel={$||\mask{\potgrid} - \pot(\xgridvec[\centers])||_\infty/||\pot(\xgridvec[\centers])||_\infty$},
        width=5cm,
        height=5.5cm
    ]

        \addplot[very thick, gray, solid] table[x index=0, y expr=\thisrowno{5} * 2] {figures/data/2D_poisson/poisson_original_method_error_dsdx_1.3_ddf_Yang3.txt};
        
        \addplot[very thick, gray, solid] table[x index=0, y expr=\thisrowno{6} / 2] {figures/data/2D_poisson/poisson_new_method_error_dsdx_1.3_ddf_Yang3.txt};    

        \node[anchor=south east, text=gray] at (axis cs:0.17,0.12) {$\sim \dx$};
        \node[anchor=north west, text=gray] at (axis cs:0.07,0.002) {$\sim \dx^2$};

        \foreach \method/\marker/\color in {
        original/diamond/red,
        new/triangle/blue}{
            \edef\temp{
                \noexpand\addplot[
                    solid,
                    \color,
                    mark=\marker, line width=1.0,
                    mark options={solid}
                ] table [x index=0, y index=1] {figures/data/2D_poisson/poisson_\method_method_error_dsdx_1.3_ddf_Yang3.txt};
                \noexpand\addplot[
                    dashed,
                    \color,
                    mark=\marker, line width=1.0,
                    mark options={solid}
                ] table [x index=0, y index=3] {figures/data/2D_poisson/poisson_\method_method_error_dsdx_1.3_ddf_Yang3.txt};
                \noexpand\addplot[
                    solid,
                    \color!060!white, line width=0.7,
                    mark=\marker,
                    mark options={solid}
                ] table [x index=0, y index=1] {figures/data/2D_poisson/poisson_\method_method_error_dsdx_0.7_ddf_Yang3.txt};
                \noexpand\addplot[
                    dashed,
                    \color!060!white, line width=0.7,
                    mark=\marker,
                    mark options={solid}
                ] table [x index=0, y index=3] {figures/data/2D_poisson/poisson_\method_method_error_dsdx_0.7_ddf_Yang3.txt};
                \noexpand\addplot[
                    solid,
                    \color!040!white, line width=0.3,
                    mark=\marker,
                    mark options={solid}
                ] table [x index=0, y index=1] {figures/data/2D_poisson/poisson_\method_method_error_dsdx_0.1_ddf_Yang3.txt};
                \noexpand\addplot[
                    dashed,
                    \color!040!white, line width=0.3,
                    mark=\marker,
                    mark options={solid}
                ] table [x index=0, y index=3] {figures/data/2D_poisson/poisson_\method_method_error_dsdx_0.1_ddf_Yang3.txt};
            }
            \temp
        }
        
    \end{loglogaxis}

\end{tikzpicture}

%% file: figures/2D_poisson_error_L2_norm_dsdx_variation.tikz
\begin{tikzpicture}
    \begin{loglogaxis}[
        xlabel=$\dx/R$,
        ylabel={$||\mask{\potgrid} - \pot(\xgridvec[\centers])||_2/||\pot(\xgridvec[\centers])||_2$},
        width=5cm,
        height=5.5cm
    ]

        \addplot[very thick, gray, solid] table[x index=0, y expr=\thisrowno{7} * 2] {figures/data/2D_poisson/poisson_original_method_error_dsdx_1.3_ddf_Yang3.txt};
    
        \addplot[very thick, gray, solid] table[x index=0, y expr=\thisrowno{8} / 2] {figures/data/2D_poisson/poisson_new_method_error_dsdx_1.3_ddf_Yang3.txt};    

        \node[anchor=south east, text=gray] at (axis cs:0.03,0.04) {$\sim \dx$};
        \node[anchor=north west, text=gray] at (axis cs:0.07,0.002) {$\sim \dx^2$};

        \foreach \method/\marker/\color in {
        original/diamond/red,
        new/triangle/blue}{
            \edef\temp{
                \noexpand\addplot[
                    solid,
                    \color,
                    mark=\marker, line width=1.0,
                    mark options={solid}
                ] table [x index=0, y index=2] {figures/data/2D_poisson/poisson_\method_method_error_dsdx_1.3_ddf_Yang3.txt};
                \noexpand\addplot[
                    dashed,
                    \color,
                    mark=\marker, line width=1.0,
                    mark options={solid}
                ] table [x index=0, y index=4] {figures/data/2D_poisson/poisson_\method_method_error_dsdx_1.3_ddf_Yang3.txt};
                \noexpand\addplot[
                    solid,
                    \color!060!white, line width=0.7,
                    mark=\marker,
                    mark options={solid}
                ] table [x index=0, y index=2] {figures/data/2D_poisson/poisson_\method_method_error_dsdx_0.7_ddf_Yang3.txt};
                \noexpand\addplot[
                    dashed,
                    \color!060!white, line width=0.7,
                    mark=\marker,
                    mark options={solid}
                ] table [x index=0, y index=4] {figures/data/2D_poisson/poisson_\method_method_error_dsdx_0.7_ddf_Yang3.txt};
                \noexpand\addplot[
                    solid,
                    \color!040!white, line width=0.3,
                    mark=\marker,
                    mark options={solid}
                ] table [x index=0, y index=2] {figures/data/2D_poisson/poisson_\method_method_error_dsdx_0.1_ddf_Yang3.txt};
                \noexpand\addplot[
                    dashed,
                    \color!040!white, line width=0.3,
                    mark=\marker,
                    mark options={solid}
                ] table [x index=0, y index=4] {figures/data/2D_poisson/poisson_\method_method_error_dsdx_0.1_ddf_Yang3.txt};
            }
            \temp
        }
        
    \end{loglogaxis}

\end{tikzpicture}

%% file: figures/2D_poisson_forcing_error_Linf_norm_dsdx_variation.tikz
\begin{tikzpicture}
    \begin{loglogaxis}[
        ymin=0.000002,
        ymax=300,
        xlabel=$\dx/R$,
        ylabel={$||\jump{\potpointnormal} - \jump{\pot^\normalbase}||_\infty/||\jump{\pot^\normalbase}||_\infty$},
        width=5cm,
        height=5.5cm
    ]
     
        \addplot[very thick, gray, solid] table[x index=0, y expr=\thisrowno{2} * 2] {figures/data/2D_poisson/poisson_original_method_forcing_error_dsdx_1.3_ddf_Yang3.txt};    
  
        \addplot[very thick, gray, solid] table[x index=0, y expr=\thisrowno{3} / 2] {figures/data/2D_poisson/poisson_new_method_forcing_error_dsdx_1.3_ddf_Yang3.txt};    
    
        \node[anchor=south west, text=gray] at (axis cs:0.07,0.3) {$\sim \dx$};
        \node[anchor=north west, text=gray] at (axis cs:0.07,0.001) {$\sim \dx^2$};

        \foreach \method/\marker/\color in {
        original/diamond/red,
        new/triangle/blue}{
            \edef\temp{
                \noexpand\addplot[
                    solid,
                    \color,
                    mark=\marker, line width=1.0,
                    mark options={solid}
                ] table [x index=0, y index=1] {figures/data/2D_poisson/poisson_\method_method_forcing_error_dsdx_1.3_ddf_Yang3.txt};
                \noexpand\addplot[
                    solid,
                    \color!060!white, line width=0.7,
                    mark=\marker,
                    mark options={solid}
                ] table [x index=0, y index=1] {figures/data/2D_poisson/poisson_\method_method_forcing_error_dsdx_0.7_ddf_Yang3.txt};
                \noexpand\addplot[
                    solid,
                    \color!040!white, line width=0.3,
                    mark=\marker,
                    mark options={solid}
                ] table [x index=0, y index=1] {figures/data/2D_poisson/poisson_\method_method_forcing_error_dsdx_0.1_ddf_Yang3.txt};
            }
            \temp
        }
        
    \end{loglogaxis}

\end{tikzpicture}

%% file: figures/2D_poisson_condition_number.tikz
\begin{tikzpicture} 
    \begin{loglogaxis}[
        ymin=1,
        ymax=1.1e20,
        xlabel=$\ds/\dx$,
        ylabel={Condition number of $\schur$ and $\tilde{\schur}$},
        width=5cm,
        height=6cm
    ]

    \foreach \method/\linestyle/\color in {
        original/dotted/red,
        new/dashed/blue
    }{
        \edef\temp{
            \noexpand\addplot[
                \linestyle,
                \color,
                mark=o,
                mark options={solid}
            ] table [x index=0, y index=1] {figures/data/2D_poisson/poisson_\method_method_condition_number_dx_0.2_ddf_Yang3.txt};
            \noexpand\addplot[
                \linestyle,
                \color,
                mark=square,
                mark options={solid}
            ] table [x index=0, y index=1] {figures/data/2D_poisson/poisson_\method_method_condition_number_dx_0.1_ddf_Yang3.txt};
            \noexpand\addplot[
                \linestyle,
                \color,
                mark=x,
                mark options={solid}
            ] table [x index=0, y index=1] {figures/data/2D_poisson/poisson_\method_method_condition_number_dx_0.05_ddf_Yang3.txt};
        }
        \temp
    }

    \end{loglogaxis}

\end{tikzpicture}

%% file: figures/2D_poisson_forcing_original.tikz
\begin{tikzpicture}
    \begin{axis}[
        xmin=-0.2,
        xmax=6.48,
        ymin=-2.8,
        ymax=2.8,
        xlabel=$\theta \vphantom{\ds/\dx}$,
        ylabel=$\jump{\potpointnormal}$,
        width=6.25cm,
        height=6cm,
        clip marker paths=true,
        line join=round,
    ]

    \addplot[black, line width=1.9, dashed] table[x index=0, y index=1] {figures/data/2D_poisson/poisson_exact_forcing.txt};

    \addplot[red!100!white, line width=1.00,  mark options={scale=0.50, line join=miter}] table[x index=0,y index=1] {figures/data/2D_poisson/poisson_original_method_forcing_dx_0.1_dsdx_1.3_ddf_Yang3.txt};
    \addplot[red!060!white, opacity=0.75, line width=0.50,  mark options={scale=0.40, line join=miter}] table[x index=0,y index=1] {figures/data/2D_poisson/poisson_original_method_forcing_dx_0.1_dsdx_0.7_ddf_Yang3.txt};

    \end{axis}

\end{tikzpicture}

%% file: figures/2D_poisson_forcing_new.tikz
\begin{tikzpicture}
    \begin{axis}[
        xmin=-0.2,
        xmax=6.48,
        ymin=-2.8,
        ymax=2.8,
        xlabel=$\theta \vphantom{\ds/\dx}$,
        yticklabels=\empty,
        width=6.25cm,
        height=6cm,
        clip marker paths=true,
        line join=round
    ]

    \addplot[black, line width=1.9, dashed] table[x index=0, y index=1] {figures/data/2D_poisson/poisson_exact_forcing.txt};

    \addplot[blue!100!white, line width=1.00,mark options={scale=0.50, line join=miter}] table[x index=0,y index=1] {figures/data/2D_poisson/poisson_new_method_forcing_dx_0.1_dsdx_1.3_ddf_Yang3.txt};
    \addplot[blue!060!white, opacity=0.75, line width=0.50,mark options={scale=0.40, line join=miter}] table[x index=0,y index=1] {figures/data/2D_poisson/poisson_new_method_forcing_dx_0.1_dsdx_0.7_ddf_Yang3.txt};
    \addplot[blue!030!white, opacity=0.5, line width=0.25,mark options={scale=0.30, line join=miter}] table[x index=0,y index=1] {figures/data/2D_poisson/poisson_new_method_forcing_dx_0.1_dsdx_0.1_ddf_Yang3.txt};

    \end{axis}

\end{tikzpicture}

%% file: figures/2D_double_cylinder_diagram.tikz
\begin{tikzpicture}
\def\Rone{1.0}
\def\Rtwo{2.0}
\def\Roffset{0.5}
\def\ut{1.2}
\pgfmathsetmacro{\rmid}{0.5*(\Rone+\Rtwo)}

\begin{axis}[
    width = 6.0cm,
    height = 6.0cm,
    xlabel = $x$, 
    ylabel = $y$
]
    
    \addplot[domain=0:360, samples=200, lightgray, thick] 
        ({\Rone*cos(x)}, {\Rone*sin(x)});
    
    \addplot[domain=0:360, samples=200, lightgray, thick] 
        ({\Rtwo*cos(x)}, {\Rtwo*sin(x)});
    
    \node at (axis cs:-0.9,-0.9) {$\omega$};
    \addplot[
        domain=160:215,
        samples=100,
        ->,
        -{stealth},
    ]
    ({1.25*cos(x)}, {1.25*sin(x)});
    
    \node at (axis cs:0.707*\Rone-0.707*\Roffset,-0.707*\Rone+0.707*\Roffset) {$\domain^+$};
    \node at (axis cs:0.707*\rmid,-0.707*\rmid) {$\domain^-$};
    \node at (axis cs:0.707*\Rtwo+0.707*\Roffset,-0.707*\Rtwo-0.707*\Roffset) {$\domain^+$};
    \node at (axis cs:0.8,1.4) {$\vel_\theta(r)$};
    
    \foreach \r in {0.3, 0.5, 0.7, 0.9, 1.1, 1.3, 1.5, 1.7} {
        \pgfmathsetmacro{\v}{
            (\r <= \Rone) 
            * (\ut / \Rone * \r) 
            + (\r > \Rone) 
            * (\ut / (1 - (\Rone/\Rtwo)^2) * (\Rone/\r - \r*\Rone/\Rtwo^2))
        }
        \addplot[-stealth] coordinates {(\r,0) (\r,\v)};
    }
    
    \addplot[domain=0:\Rone, samples=100] 
        ({x}, {\ut/\Rone*x});
    
    \addplot[domain=\Rone:\Rtwo, samples=100] 
        ({x}, {\ut/(1 - (\Rone/\Rtwo)^2)*(\Rone/x - x*\Rone/\Rtwo^2)});

    \coordinate (ANCHOR1) at (axis cs: -0.707*\Rone, 0.707*\Rone);
    \coordinate (CTRL11)  at (axis cs: -0.707*\Rone - 0.1, 0.707*\Rone + 0.1);
    \coordinate (CTRL12)  at (axis cs: -0.707*\Rone - 0.2, 0.707*\Rone + 0.1);
    \coordinate (LABEL1)  at (axis cs: -0.707*\Rone - 0.2, 0.707*\Rone + 0.26);
    \draw [] (ANCHOR1) .. controls (CTRL11)  and (CTRL12) .. (LABEL1) node[above] {$\surf_1$};
    
    \coordinate (ANCHOR2) at (axis cs: -0.707*\Rtwo, 0.707*\Rtwo);
    \coordinate (CTRL21)  at (axis cs: -0.707*\Rtwo - 0.1, 0.707*\Rtwo + 0.1);
    \coordinate (CTRL22)  at (axis cs: -0.707*\Rtwo - 0.2, 0.707*\Rtwo + 0.1);
    \coordinate (LABEL2)  at (axis cs: -0.707*\Rtwo - 0.2, 0.707*\Rtwo + 0.26);
    \draw [] (ANCHOR2) .. controls (CTRL21)  and (CTRL22) .. (LABEL2) node[above] {$\surf_2$};

\end{axis}

\end{tikzpicture}

%% file: figures/2D_double_cylinder_solution_original.tikz
\begin{tikzpicture}
    \begin{axis}[
        name=mainplot,
        xmin=0,
        xmax=2.3,
        ymin=-0.1,
        ymax=1.1,
        xlabel=$x$,
        ylabel={$\mask{\velgridcomponent}_y(\xgrid[\facescomponent_y],y=0)$},
        width=5.75cm,
        height=6cm,
    ]

    \addplot[black, line width=1.0, dashed] table[x expr=\thisrowno{0}/0.75, y expr=\thisrowno{1}/0.15] {figures/data/2D_double_cylinder/double_cylinder_exact_profile.txt};
    

    \addplot[red!100!white, line width=1, mark=diamond*, mark options={scale=0.50, line join=miter}] table[x expr=\thisrowno{0}/0.75,y expr=\thisrowno{1}/0.15] {figures/data/2D_double_cylinder/ns_original_method_v_centerline_dx_0.125_dsdx_1.0_ddf_Yang3.txt};
    \addplot[red!060!white, line width=0.75, mark=diamond*,  mark options={scale=0.40, line join=miter}] table[x expr=\thisrowno{0}/0.75,y expr=\thisrowno{1}/0.15] {figures/data/2D_double_cylinder/ns_original_method_v_centerline_dx_0.0625_dsdx_1.0_ddf_Yang3.txt};
    \addplot[red!030!white, line width=0.5, mark=diamond*,  mark options={scale=0.30, line join=miter}] table[x expr=\thisrowno{0}/0.75,y expr=\thisrowno{1}/0.15] {figures/data/2D_double_cylinder/ns_original_method_v_centerline_dx_0.03125_dsdx_1.0_ddf_Yang3.txt};
 
    \draw[thick] (axis cs:1.88,-0.05) coordinate (zoomSW) rectangle (axis cs:2.22,0.06) coordinate (zoomNE);
    \coordinate (zoomWest) at ($(zoomSW)!0.5!(zoomSW |- zoomNE)$);
    \end{axis}

    \begin{axis}[
        name=inset,
        width=3.3cm,
        height=3.2cm,
        xmin=1.88, xmax=2.22, 
        ymin=-0.05, ymax=0.06, 
        at={(mainplot.south)},
        xshift=-3mm, yshift=2mm,
        anchor=south,
        xtick=\empty,
        ytick=\empty,
        clip marker paths=true
    ]

    \addplot[black, line width=1.0, dashed] table[x expr=\thisrowno{0}/0.75, y expr=\thisrowno{1}/0.15] {figures/data/2D_double_cylinder/double_cylinder_exact_profile.txt};


    \addplot[red!100!white, line width=1, mark=diamond*, mark options={scale=0.60, line join=miter}] table[x expr=\thisrowno{0}/0.75,y expr=\thisrowno{1}/0.15] {figures/data/2D_double_cylinder/ns_original_method_v_centerline_dx_0.125_dsdx_1.0_ddf_Yang3.txt};
    \addplot[red!060!white, line width=0.75, mark=diamond*,  mark options={scale=0.50, line join=miter}] table[x expr=\thisrowno{0}/0.75,y expr=\thisrowno{1}/0.15] {figures/data/2D_double_cylinder/ns_original_method_v_centerline_dx_0.0625_dsdx_1.0_ddf_Yang3.txt};
    \addplot[red!030!white, line width=0.5, mark=diamond*,  mark options={scale=0.40, line join=miter}] table[x expr=\thisrowno{0}/0.75,y expr=\thisrowno{1}/0.15] {figures/data/2D_double_cylinder/ns_original_method_v_centerline_dx_0.03125_dsdx_1.0_ddf_Yang3.txt};
 
    \end{axis}
    \draw[thick] (zoomWest) -- (inset.east);
\end{tikzpicture}

%% file: figures/2D_double_cylinder_solution_new.tikz
\begin{tikzpicture}
    \begin{axis}[
        name=mainplot,
        xmin=0,
        xmax=2.3,
        ymin=-0.1,
        ymax=1.1,
        xlabel=$x$,
        width=5.75cm,
        height=6cm,
        yticklabels=\empty,
    ]

    \addplot[black, line width=1.0, dashed] table[x expr=\thisrowno{0}/0.75, y expr=\thisrowno{1}/0.15] {figures/data/2D_double_cylinder/double_cylinder_exact_profile.txt};

    \addplot[blue!100!white, line width=1, mark=triangle*, mark options={scale=0.50, line join=miter}] table[x expr=\thisrowno{0}/0.75,y expr=\thisrowno{1}/0.15] {figures/data/2D_double_cylinder/ns_new_method_v_centerline_dx_0.125_dsdx_1.0_ddf_Yang3.txt};
    \addplot[blue!060!white, line width=0.75, mark=triangle*, mark options={scale=0.40, line join=miter}] table[x expr=\thisrowno{0}/0.75,y expr=\thisrowno{1}/0.15] {figures/data/2D_double_cylinder/ns_new_method_v_centerline_dx_0.0625_dsdx_1.0_ddf_Yang3.txt};
    \addplot[blue!030!white, line width=0.5, mark=triangle*, mark options={scale=0.30, line join=miter}] table[x expr=\thisrowno{0}/0.75,y expr=\thisrowno{1}/0.15] {figures/data/2D_double_cylinder/ns_new_method_v_centerline_dx_0.03125_dsdx_1.0_ddf_Yang3.txt};

    \draw[thick] (axis cs:1.88,-0.05) coordinate (zoomSW) rectangle (axis cs:2.22,0.06) coordinate (zoomNE);
    \coordinate (zoomWest) at ($(zoomSW)!0.5!(zoomSW |- zoomNE)$);
    \end{axis}

    \begin{axis}[
        name=inset,
        width=3.3cm,
        height=3.2cm,
        xmin=1.88, xmax=2.22, 
        ymin=-0.05, ymax=0.06, 
        at={(mainplot.south)},
        xshift=-3mm, yshift=2mm,
        anchor=south,
        xtick=\empty,
        ytick=\empty,
        clip marker paths=true
    ]

    \addplot[black, line width=1.0, dashed] table[x expr=\thisrowno{0}/0.75, y expr=\thisrowno{1}/0.15] {figures/data/2D_double_cylinder/double_cylinder_exact_profile.txt};
    

    \addplot[blue!100!white, line width=1, mark=triangle*, mark options={scale=0.60, line join=miter}] table[x expr=\thisrowno{0}/0.75,y expr=\thisrowno{1}/0.15] {figures/data/2D_double_cylinder/ns_new_method_v_centerline_dx_0.125_dsdx_1.0_ddf_Yang3.txt};
    \addplot[blue!060!white, line width=0.75, mark=triangle*, mark options={scale=0.50, line join=miter}] table[x expr=\thisrowno{0}/0.75,y expr=\thisrowno{1}/0.15] {figures/data/2D_double_cylinder/ns_new_method_v_centerline_dx_0.0625_dsdx_1.0_ddf_Yang3.txt};
    \addplot[blue!030!white, line width=0.5, mark=triangle*, mark options={scale=0.40, line join=miter}] table[x expr=\thisrowno{0}/0.75,y expr=\thisrowno{1}/0.15] {figures/data/2D_double_cylinder/ns_new_method_v_centerline_dx_0.03125_dsdx_1.0_ddf_Yang3.txt};

    \end{axis}
    \draw[thick] (zoomWest) -- (inset.east);
\end{tikzpicture}

%% file: figures/2D_double_cylinder_error_Linf_norm.tikz
\begin{tikzpicture}
    \begin{loglogaxis}[
        xlabel=$\dx$,
        ylabel={$||\mask{\velgrid} - \vel(\xgridvec[\faces])||_\infty/||\vel(\xgridvec[\faces])||_\infty$},
        width=6cm,
        height=5.5cm
    ]

        \addplot[very thick, gray, solid] table[x expr=\thisrowno{0}/0.75, y expr=\thisrowno{5} * 3.1] {figures/data/2D_double_cylinder/ns_original_method_error_dsdx_1.0_ddf_Yang3.txt};
        
        \addplot[very thick, gray, solid] table[x expr=\thisrowno{0}/0.75, y expr=\thisrowno{6} / 1.6] {figures/data/2D_double_cylinder/ns_new_method_error_dsdx_1.0_ddf_Yang3.txt};    
        
        \addplot[mark=triangle, red, solid, mark options={solid, scale=1.0}] table[x expr=\thisrowno{0}/0.75, y index=1] {figures/data/2D_double_cylinder/ns_original_method_error_dsdx_1.0_ddf_Yang3.txt};
        
        \addplot[mark=triangle, red, dashed, mark options={solid, scale=1.0}] table[x expr=\thisrowno{0}/0.75, y index=3] {figures/data/2D_double_cylinder/ns_original_method_error_dsdx_1.0_ddf_Yang3.txt};   
        \addplot[mark=triangle, blue, solid, mark options={solid, scale=1.0}] table[x expr=\thisrowno{0}/0.75, y index=1] {figures/data/2D_double_cylinder/ns_new_method_error_dsdx_1.0_ddf_Yang3.txt};
    
        \addplot[mark=triangle, blue, dashed, mark options={solid, scale=1.0}] table[x expr=\thisrowno{0}/0.75, y index=3] {figures/data/2D_double_cylinder/ns_new_method_error_dsdx_1.0_ddf_Yang3.txt};

        \node[anchor=south east, text=gray] at (axis cs:0.04,0.05) {$\sim \dx$};
        \node[anchor=north west, text=gray] at (axis cs:0.05,0.0015) {$\sim \dx^2$};
    \end{loglogaxis}

\end{tikzpicture}

%% file: figures/2D_double_cylinder_error_L2_norm.tikz
\begin{tikzpicture}   
    \begin{loglogaxis}[
        xlabel=$\dx$,
        ylabel={$||\mask{\velgrid} - \vel(\xgridvec[\faces])||_2/||\vel(\xgridvec[\faces])||_2$},
        width=6cm,
        height=5.5cm
    ]

        \addplot[very thick, gray, solid] table[x expr=\thisrowno{0}/0.75, y expr=\thisrowno{7} * 1.4] {figures/data/2D_double_cylinder/ns_original_method_error_dsdx_1.0_ddf_Yang3.txt};
        
        \addplot[very thick, gray, solid] table[x expr=\thisrowno{0}/0.75, y expr=\thisrowno{8} / 1.9] {figures/data/2D_double_cylinder/ns_new_method_error_dsdx_1.0_ddf_Yang3.txt};    
        
        \addplot[mark=triangle, red, solid, mark options={solid, scale=1.0}] table[x expr=\thisrowno{0}/0.75, y index=2] {figures/data/2D_double_cylinder/ns_original_method_error_dsdx_1.0_ddf_Yang3.txt};
        
        \addplot[mark=triangle, red, dashed, mark options={solid, scale=1.0}] table[x expr=\thisrowno{0}/0.75, y index=4] {figures/data/2D_double_cylinder/ns_original_method_error_dsdx_1.0_ddf_Yang3.txt};      
        \addplot[mark=triangle, blue, solid, mark options={solid, scale=1.0}] table[x expr=\thisrowno{0}/0.75, y index=2] {figures/data/2D_double_cylinder/ns_new_method_error_dsdx_1.0_ddf_Yang3.txt};
    
        \addplot[mark=triangle, blue, dashed, mark options={solid, scale=1.0}] table[x expr=\thisrowno{0}/0.75, y index=4] {figures/data/2D_double_cylinder/ns_new_method_error_dsdx_1.0_ddf_Yang3.txt};

        \node[anchor=south east, text=gray] at (axis cs:0.04,0.05) {$\sim \dx$};
        \node[anchor=north west, text=gray] at (axis cs:0.05,0.0015) {$\sim \dx^2$};
    \end{loglogaxis}

\end{tikzpicture}

%% file: figures/2D_double_cylinder_condition_number.tikz
\begin{tikzpicture}
    \begin{loglogaxis}[
        xlabel=$\ds/\dx$,
        ylabel={Condition number of $\schur$ and $\tilde{\schur}$},
        width=4cm,
        height=6cm
    ]

    \foreach \method/\linestyle/\color in {
        original/dotted/red,
        new/dashed/blue
    }{
        \edef\temp{
            \noexpand\addplot[
                \color,
                \linestyle,
                mark=o,
                mark options={solid}
            ] table [x index=0, y index=1] {figures/data/2D_double_cylinder/ns_\method_method_condition_number_dx_0.125_ddf_Yang3.txt};
            \noexpand\addplot[
                \color,
                \linestyle,
                mark=square,
                mark options={solid}
            ] table [x index=0, y index=1] {figures/data/2D_double_cylinder/ns_\method_method_condition_number_dx_0.0625_ddf_Yang3.txt};
        }
        \temp
    }

    \end{loglogaxis}

\end{tikzpicture}

%% file: figures/2D_double_cylinder_forcing_original.tikz
\begin{tikzpicture}
    \begin{axis}[
        width=6.25cm,
        height=6cm,
        xmin=-0.2,
        xmax=6.48,
        ymin=-5.0,
        ymax=5.0,
        ylabel={$\forcepoint_{x, \surf_1} \Re$},
        xlabel=$\theta \vphantom{\ds/\dx}$,
        clip marker paths=true,
        line join=round
    ]

    \addplot[black, line width=1.0, dashed] table[x index=0, y expr=\thisrowno{1}/0.2] {figures/data/2D_double_cylinder/double_cylinder_exact_inner_x_forcing.txt};
    

    \addplot[red!100!white, line width=1.5, mark options={scale=0.50, line join=miter}] table[x index=0,y expr=\thisrowno{1}/0.2] {figures/data/2D_double_cylinder/ns_original_method_inner_x_forcing_dx_0.0625_dsdx_1.3_ddf_Yang3.txt};
    \addplot[red!060!white, line width=0.750,  mark options={scale=0.40, line join=miter}] table[x index=0,y expr=\thisrowno{1}/0.2] {figures/data/2D_double_cylinder/ns_original_method_inner_x_forcing_dx_0.0625_dsdx_1.0_ddf_Yang3.txt};
    
    \end{axis}

\end{tikzpicture}

%% file: figures/2D_double_cylinder_forcing_new.tikz
\begin{tikzpicture}
    \begin{axis}[
        width=6.25cm,
        height=6cm,
        xmin=-0.2,
        xmax=6.48,
        ymin=-5.0,
        ymax=5.0,
        ylabel={$\jump[\surf_1]{\velcomponentpointnormal{x}}$},
        xlabel=$\theta \vphantom{\ds/\dx}$,
        clip marker paths=true,
        line join=round
    ]

    \addplot[black, line width=1.0, dashed] table[x index=0, y expr=\thisrowno{1}/0.2] {figures/data/2D_double_cylinder/double_cylinder_exact_inner_x_forcing.txt};
    
    
    \addplot[blue!100!white, line width=1.5, mark options={scale=0.50, line join=miter}] table[x index=0,y expr=\thisrowno{1}/0.2] {figures/data/2D_double_cylinder/ns_new_method_inner_x_forcing_dx_0.0625_dsdx_1.3_ddf_Yang3.txt};
    \addplot[blue!060!white, line width=0.750,  mark options={scale=0.40, line join=miter}] table[x index=0,y expr=\thisrowno{1}/0.2] {figures/data/2D_double_cylinder/ns_new_method_inner_x_forcing_dx_0.0625_dsdx_1.0_ddf_Yang3.txt};
    \addplot[blue!030!white, line width=0.25,  mark options={scale=0.30, line join=miter}] table[x index=0,y expr=\thisrowno{1}/0.2] {figures/data/2D_double_cylinder/ns_new_method_inner_x_forcing_dx_0.0625_dsdx_0.7_ddf_Yang3.txt};
    \end{axis}

\end{tikzpicture}